\documentclass{amsart}
\usepackage[nohug]{diagrams}
\usepackage{graphicx}
\usepackage[mathscr]{eucal}
\usepackage{amssymb, amsmath, amsthm}
\usepackage{amsfonts}
\usepackage{latexsym}
\usepackage{tabularx,array}
\usepackage{enumerate}
\usepackage[all,arc]{xy}
\usepackage{mathrsfs}
\usepackage{latexsym}
\usepackage{color}
\usepackage{mathtools}
\usepackage{leftidx}
\usepackage{bm}
\usepackage{url}
\newfont{\msam}{msam10}

\newtheorem{theorem}[]{Theorem}
\newtheorem{proposition}[]{Proposition}
\newtheorem{corollary}[]{Corollary}
\newtheorem{lemma}[]{Lemma}

\theoremstyle{definition}
\newtheorem{definition}[]{Definition}
\newtheorem{defn}[theorem]{Definition}
\newtheorem{remark}[]{Remark}
\newtheorem{example}[]{Example}
\newtheorem{conj}[]{Conjecture}

\let\nc\newcommand

\def\bthm{\begin{theorem}}
\def\ethm{\end{theorem}}
\def\blemma{\begin{lemma}}
\def\elemma{\end{lemma}}
\def\bproof{\begin{proof}}
\def\eproof{\end{proof}}
\def\bprop{\begin{proposition}}
\def\eprop{\end{proposition}}
\def\bcor{\begin{corollary}}
\def\ecor{\end{corollary}}
\def\bconj{\begin{conj}}
\def\econj{\end{conj}}
\nc{\la}{\label}
\def\O{\mathcal{O}}
\def\Z{\mathbb{Z}}
\def\N{\mathbb{N}}
\def\Q{\mathbb{Q}}
\def\c{\mathbb{C}}
\def\M{\mathcal{M}}

\def\L {\boldsymbol{L}}
\def\R{\mathbb{R}}

\def\Alg{\mathtt{Alg}}
\def\sAlg{\mathtt{sAlg}}

\def\Mod{\mathtt{Mod}}

\def\Bimod{\mathtt{Bimod}}
\def\cAlg{\mathtt{Comm\,Alg}}
\def\scAlg{\mathtt{sComm\,Alg}}
\def\gcAlg{\mathtt{grComm\,Alg}}
\def\Sets{\mathtt{Set}}

\def\D{\mathcal{D}}
\def\C{\mathcal{C}}
\def\W{\mathcal{W}}

\def\M{\mathcal{M}}
\def\N{\mathcal{N}}
\def\E{\mathcal{E}}

\def\Ho{{\mathtt{Ho}}}

\nc{\ocolim}{{\rm ocolim}}
\nc{\Ob}{{\rm Ob}}
\nc{\Hom}{{\rm{Hom}}}
\nc{\Homcont}{{\mathcal{H}om}}
\nc{\HOM}{\underline{\rm{Hom}}}
\nc{\DER}{\underline{\rm{Der}}}
\nc{\END}{\underline{\rm{End}}}
\nc{\bSym}{\mathbf{Sym}}
\nc{\Ext}{{\rm{Ext}}}
\nc{\Rep}{{\rm{Rep}}}
\nc{\DRep}{{\rm{DRep}}}
\nc{\ODRep}{{\mathcal O}{\rm{DRep}}}
\nc{\NCRep}{\widetilde{\rm{Rep}}}
\nc{\RAct}{{\rm{RAct}}}
\nc{\bs}{\backslash}
\nc{\ob}{{\tt{Obs}}}
\nc{\CE}{\mathcal{C}}
\nc{\TP}{{T\!P}}
\nc{\un}{\underline{n}}
\nc{\um}{\underline{m}}
\nc{\rn}{\langle n \rangle}
\nc{\nn}{{{\natural} {\natural}}}
\nc{\n}{{{\natural}}}
\nc{\A}{\mathbb A}
\nc{\B}{{\mathrm{B}}}
\nc{\Ba}{\overline{\mathrm{B}}}
\nc{\bC}{\overline{C}}
\nc{\bOmega}{\boldsymbol{\Omega}}
\nc{\bB}{\boldsymbol{B}}
\nc{\EXT}{\underline{\rm{Ext}}}
\nc{\TOR}{\underline{\rm{Tor}}}
\def\H{\mathrm H}

\def\HR{\mathrm{HR}}

\nc{\End}{{\rm{End}}}
\nc{\GL}{{\rm{GL}}}
\nc{\gl}{{\mathfrak{gl}}}
\nc{\rgl}{\overline{{\mathfrak{gl}}}}
\nc{\g}{{\mathfrak{g}}}
\nc{\h}{{\mathfrak{h}}}
\nc{\PGL}{{\rm{PGL}}}
\nc{\SL}{{\rm{SL}}}
\nc{\sll}{\mathfrak{sl}}
\nc{\cn}{ \mbox{\rm c\^{o}ne} }
\nc{\PSL}{{\rm{PSL}}}
\nc{\ad}{{\rm{ad}}}
\nc{\Ad}{{\rm{Ad}}}
\nc{\dlim}{\varinjlim}
\nc{\plim}{\varprojlim}
\nc{\colim}{{{\rm colim}}}

\newcommand{\bL}{\mathbb{L}}
\newcommand{\bR}{\boldsymbol{R}}

\newcommand{\HH}{{\rm{HH}}}

\newcommand{\Tor}{{\rm{Tor}}}
\newcommand{\Spec}{{\rm{Spec}}}

\newcommand{\Sym}{\Lambda}
\newcommand{\Aut}{{\rm{Aut}}}

\newcommand{\id}{{\rm{Id}}}

\newcommand{\Ker}{{\rm{Ker}}}

\newcommand{\diag}{{\rm{diag}}}

\newcommand{\im}{{\rm{Im}}}

\newcommand{\into}{\,\hookrightarrow\,}

\newcommand{\onto}{\,\twoheadrightarrow\,}
\newcommand{\sonto}{\,\stackrel{\sim}{\twoheadrightarrow}\,}
\def\Cc{\mathscr{C}}
\def\cb{\boldsymbol{\Omega}}

\def\bs{\backslash}

\def\fin{\mathfrak{F}}
\def\ve{\mathtt{Vect}}

\def\sGr{\mathtt{sGr}}
\def\Mon{\mathtt{Mon}}
\def\sMon{\mathtt{sMon}}
\def\Gr{\mathtt{Gr}}
\def\fgr{\mathtt{FGr}}
\def\ffgr{\mathfrak{G}}
\def\set{\mathtt{Set}}
\def\sset{\mathtt{sSet}}
\def\sSp{\mathtt{sSp}}
\def\LsSp{{\mathfrak L}\mathtt{sSp}}
\def\lgr{\mathbb{G}}
\def\lin{\mathrm{lin}}

\newcommand{\rar}{\xrightarrow{}}

\newcommand\frgr[1]{\langle #1 \rangle}
\nc{\env}{\mathrm{End}(V)}
\nc{\FT}{\mathcal{C}}

\numberwithin{equation}{section}
\numberwithin{theorem}{section}
\numberwithin{lemma}{section}
\numberwithin{proposition}{section}
\numberwithin{definition}{section}
\numberwithin{corollary}{section}
\numberwithin{example}{section}
\numberwithin{remark}{section}

\newcommand{\rH}{\overline{\mathrm{H}}}
\def\cO{\mathcal O}
\def\bdf{\begin{defn}}
\def\edf{\end{defn}}

\def\brm{\begin{remark}}
\def\erm{\end{remark}}

\theoremstyle{definition}
\def\bdf{\begin{definition}}
\def\edf{\end{definition}}

\newcommand{\bD}{{\mathbb D}}
\newcommand{\bS}{{\mathbb S}}
\newcommand{\bP}{{\mathbb P}}
\newcommand{\bG}{{\mathbb G}}
\newcommand{\bF}{{\mathbb F}}
\newcommand{\bT}{{\mathbb T}}
\newcommand{\cA}{{\mathcal A}}
\newcommand{\m}{{\mathfrak m}}
\newcommand{\bHH}{{\mathbb H}{\mathbb H}}
\newcommand{\cH}{{\mathcal H}}
\newcommand{\ucH}{\underline{\mathcal H}}
\newcommand{\uN}{\underline{N}}

\newcommand{\mH}{{\mathcal H}}
\newcommand{\umH}{\underline{\mathcal H}}

\DeclareMathOperator*{\Moplus}{\text{\raisebox{0.25ex}{\scalebox{0.75}{$\bigoplus$}}}}
\DeclareMathOperator*{\Motimes}{\text{\raisebox{0.25ex}{\scalebox{0.75}{$\bigotimes$}}}}
\DeclareMathOperator*{\Mvee}{\text{\raisebox{0.25ex}{\scalebox{0.8}{$\bigvee$}}}}
\def\arbreBA{\vcenter{\xymatrix@R=2pt@C=2pt{
&&&&\\
&&&*{}\ar@{-}[ul] & \\
&&*{}\ar@{-}[uurr] \ar@{-}[uull] \ar@{-}[d]     &&\\
&&&&
}}}

\def\arbreAB{\vcenter{\xymatrix@R=2pt@C=2pt{
&&&&\\
&*{}\ar@{-}[ur] &&& \\
&&*{}\ar@{-}[uurr] \ar@{-}[uull] \ar@{-}[d]     &&\\
&&&&
}}}

\def\arbreABC{\vcenter{\xymatrix@R=1pt@C=1pt{
&&&&&&\\
&*{}\ar@{-}[ur] &&&&& \\
&&*{}\ar@{-}[uurr] &&&&\\
&&&*{}\ar@{-}[uuurrr] \ar@{-}[uuulll] \ar@{-}[d] &&&\\
&&&&&&
}}}

\def\arbreBAC{\vcenter{\xymatrix@R=1pt@C=1pt{
&&&&&&\\
&&&*{}\ar@{-}[ul] &&& \\
&&*{}\ar@{-}[uurr] &&&&\\
&&&*{}\ar@{-}[uuurrr] \ar@{-}[uuulll] \ar@{-}[d] &&&\\
&&&&&&
}}}

\def\arbreACB{\vcenter{\xymatrix@R=1pt@C=1pt{
&&&&&&\\
&*{}\ar@{-}[ur] &&&&& \\
&&&&*{}\ar@{-}[uull] &&\\
&&&*{}\ar@{-}[uuurrr] \ar@{-}[uuulll] \ar@{-}[d] &&&\\
&&&&&&
}}}

\def\arbreBCA{\vcenter{\xymatrix@R=1pt@C=1pt{
&&&&&&\\
&&&&&*{}\ar@{-}[ul] & \\
&&*{}\ar@{-}[uurr] &&&&\\
&&&*{}\ar@{-}[uuurrr] \ar@{-}[uuulll] \ar@{-}[d] &&&\\
&&&&&&
}}}

\def\arbreCAB{\vcenter{\xymatrix@R=1pt@C=1pt{
&&&&&&\\
&&&*{}\ar@{-}[ur] &&& \\
&&&&*{}\ar@{-}[uull] &&\\
&&&*{}\ar@{-}[uuurrr] \ar@{-}[uuulll] \ar@{-}[d] &&&\\
&&&&&&
}}}

\def\arbreCBA{\vcenter{\xymatrix@R=1pt@C=1pt{
&&&&&&\\
&&&&&*{}\ar@{-}[ul] & \\
&&&&*{}\ar@{-}[uull] &&\\
&&&*{}\ar@{-}[uuurrr] \ar@{-}[uuulll] \ar@{-}[d] &&&\\
&&&&&&
}}}

\def\arbreACA{\vcenter{\xymatrix@R=1pt@C=1pt{
&&&&&&\\
&*{}\ar@{-}[ur] &&&&*{}\ar@{-}[ul] & \\
&&&&&&\\
&&&*{}\ar@{-}[uuurrr] \ar@{-}[uuulll] \ar@{-}[d] &&&\\
&&&&&&
}}}

\begin{document}

\title{Representation Homology of Topological Spaces}
\author{Yuri Berest}
\address{Department of Mathematics,
Cornell University, Ithaca, NY 14853-4201, USA}
\email{berest@math.cornell.edu}
\author{Ajay C. Ramadoss}
\address{Department of Mathematics,
Indiana University,
Bloomington, IN 47405, USA}
\email{ajcramad@indiana.edu}
\author{Wai-Kit Yeung}
\address{Department of Mathematics,
Indiana University,
Bloomington, IN 47405, USA}
\email{yeungw@iu.edu}
\begin{abstract}
In this paper, we introduce and study representation homology of topological spaces, which is a natural homological extension of representation varieties of fundamental groups. We give an elementary construction of representation homology parallel to the Loday-Pirashvili construction of higher Hochschild homology; in fact, we establish a direct geometric relation between the two theories by proving that the representation homology of the suspension of a (pointed connected) space is isomorphic to its higher Hochschild homology. We also construct some natural maps and spectral sequences relating representation homology to other homology theories associated with spaces (such as Pontryagin algebras, $\bS^1$-equivariant homology of the free loop space and stable homology of automorphism groups of f.g. free groups). We compute representation homology explicitly (in terms of known invariants) in a number of interesting cases, including spheres, suspensions, complex projective spaces, Riemann surfaces and some 3-dimensional manifolds, such as link complements in $\R^3$ and the lens spaces $ L(p,q) $. In the case of link complements, we identify the representation homology in terms of ordinary Hochschild homology, which gives a new algebraic invariant of links in $\R^3$.

\end{abstract}
\maketitle
\tableofcontents
\section{Introduction}

Representation homology is an algebraic homology theory associated with derived representation schemes, which are natural (`derived') extensions of classical representation varieties. The subject may be plainly viewed as part of derived algebraic geometry (see, e.g., \cite{To}); however, somewhat surprisingly, there are more elementary  constructions. What makes representation homology interesting (it seems) are the relations between these different constructions and interpretations coming from different areas of mathematics.

In the present paper, we give two (equivalent) definitions of representation homology of topological spaces: one in terms of (non-abelian) derived functors on simplicial groups and the other in terms of classical homological algebra in functor categories. The first definition is inspired by our earlier work on representation homology of algebras (see  \cite{BKR, BR, BFPRW}), while the second by the Loday-Pirashvili approach to higher Hochschild homology \cite{P1}. Both definitions are conceptually very simple and accessible to computations: in this paper, we will use them in a complementary way to establish basic properties of representation homology and  do some examples; in our subsequent paper \cite{BRY2}, we will look at applications. 
\subsection{Representation varieties and representation homology}
We begin with some motivation for studying representation homology.
Let $G$ be a finite-dimensional affine algebraic group defined over a field $k$ of characteristic zero. For any (discrete) group $ \Gamma $, the set of all representations of $ \Gamma $ in $G$ has a natural structure
of an affine $k$-scheme called the {\it representation scheme}
$\,\Rep_G(\Gamma) $. Representation schemes and associated varieties  play an important role in many areas of mathematics, most notably in representation theory and low-dimensional topology. In representation theory, the fundamental problem is to understand the structure of representations of $ \Gamma $ in $G$. One can approach this problem geometrically by studying the natural (adjoint) action of the group $G$ on the variety $\,\Rep_G(\Gamma) $.
When $ k $ is algebraically closed and $ \Gamma $ is finitely generated,
the equivariant geometry of $ \Rep_G(\Gamma) $ is closely related to the representation theory of
$ \Gamma $:
the equivalence classes of representations of $ \Gamma $ in $G$ are in bijection with
the $G$-orbits in $ \Rep_G(\Gamma) $,
and the geometry of $G$-orbits determines the algebraic structure of representations.
This relation has been extensively studied since the late 70's, and the representation varieties have become a standard tool in representation theory of groups (see, for example, \cite{LM85, Sik12}).

In topology, one is usually interested in global algebro-geometric invariants of spaces defined in terms of representation varieties of fundamental groups. For example, if $ K $ is a knot in $ \bS^3 $,  many classical invariants of $K$ arise from its character variety $ \chi_G(K) := \Rep_G[\pi_1(X_K)]/\!/G \,$, which is the  (categorical) quotient of the representation variety of the fundamental group of the knot complement $ X_K := \bS^3\! \setminus\! K $.  These invariants include, in particular, the classical
Alexander polynomial $ \Delta_K(t) $  (in the simplest case when $ G = \c^* $, see, e.g., \cite{McM02}), the so-called $A$-polynomial $ A_K(m,l) $ (see \cite{CCGLS94}), the Casson invariant \cite{Cu01} and
the famous Chern-Simons invariant \cite{KK93}, all of which
are defined for $ G = \SL_2(\c) $. In fact, for $ G = \SL_2(\c) $, the entire character variety, or rather its coordinate ring $\, \O[\chi_G(K)]\,  $, has a purely topological interpretation
as a Kauffman bracket skein module of  $ X_K $ (see \cite{PS00}).

Despite being useful tools, the representation varieties have some intrinsic deficiencies. First of all, these varieties are usually very singular, which makes it hard to understand their geometry. Thus, in representation theory, one faces the problem of resolving singularities of $ \Rep_G(\Gamma) $. In topology, the use of representation varieties is mostly limited to (compact orientable) surfaces, hyperbolic $3$-manifolds and knot complements in $ \bS^3 $, all of which are known to be aspherical spaces. The homotopy type of such a space is completely determined by the isomorphism type of its fundamental group, which makes representation varieties of these groups very strong and efficient  invariants. For more general spaces, however, one needs to take into account a higher homotopy information, and looking at representation varieties of fundamental groups (or even, higher homotopy groups) is not enough.

A natural way to remedy these problems is to replace the representation functor $ \Rep_G \,$ with its (non-abelian) derived functor  $ \DRep_G $ much in the same way as one replaces non-exact additive functors  in classical homological algebra (such as `$ \otimes $' and `$ \Hom $') with corresponding derived
functors (`$\,\otimes^{{\boldsymbol{L}}}\,$' and `$ \mathbf{R}\Hom $'). Geometrically,
passing from the representation scheme $ \Rep_G(\Gamma) $ to the derived representation scheme $ \DRep_G(\Gamma) $ amounts to desingularizing
$ \Rep_G(\Gamma) $, while topologically, this yields a new homology theory of spaces that captures a good deal of
homotopy information and refines the classical representation varieties of fundamental groups in an interesting and nontrivial way.

To explain this idea in more precise terms, we recall that the representation scheme $ \Rep_G(\Gamma) $ is defined as the functor on the category of commutative $k$-algebras:
\begin{equation}
\la{rep0}
\Rep_G(\Gamma):\ \cAlg_{k} \to \Sets\ ,\quad A \mapsto \Hom_{\Gr}(\Gamma,\, G(A))\ ,
\end{equation}
assigning to a $k$-algebra $A$ the set of families of representations of $\Gamma $ in $G$
parametrized by the $k$-scheme $ \Spec(A) $. It is well known that the functor \eqref{rep0} is representable, and we denote the corresponding commutative algebra by $ \Gamma_G = \O[\Rep_G(\Gamma)] $: this is the coordinate ring of the affine $k$-scheme $\Rep_G(\Gamma)$.
Varying $ \Gamma $ (while keeping $ G $ fixed),  we can now regard $ \Gamma_G $  as a functor on the category of groups:
\begin{equation}
\la{rep1}
(\,\mbox{--}\,)_G :\ \Gr \to \cAlg_k\ , \quad \Gamma \mapsto \Gamma_G \ , \
\end{equation}
which we call  the {\it representation functor} in $ G $.
The functor \eqref{rep1} extends naturally to the category $ \sGr $ of simplicial groups, taking values in the category $ \scAlg_k $ of simplicial commutative algebras. Both categories  $ \sGr $ and $ \scAlg_k $
carry standard (simplicial) model structures, with weak equivalences being
the weak homotopy equivalences of underlying simplicial sets. The functor $\,(\,\mbox{--}\,)_G:\,\sGr \to \scAlg_k\,$ is not homotopy invariant:  in general, it does not preserve weak equivalences and hence
does not descend to a functor between the homotopy categories $ \Ho(\sGr) $ and $ \Ho(\scAlg_k) $. However, it is easy to check that $(\,\mbox{--}\,)_G$ takes weak equivalences between cofibrant objects in $ \sGr $ to weak equivalences in $ \scAlg_k $  (see Lemma~\ref{sgfun}). Hence, by standard homotopical algebra, it has a (total) left derived functor
\begin{equation}
\la{rep2}
\L(\,\mbox{--}\,)_G:\ \Ho(\sGr) \to \Ho(\scAlg_k)\ .
\end{equation}
We call \eqref{rep2} the {\it derived representation functor} in $G$. Heuristically, $\L(\,\mbox{--}\,)_G$ may be thought of as the ``best possible"  approximation of the representation functor \eqref{rep1} at the level of homotopy categories. When applied to a simplicial group $ \Gamma $,  the functor \eqref{rep2} is represented by a simplicial commutative algebra which we denote by $ \O[\DRep_G(\Gamma)] $. The derived representation scheme
$ \DRep_G(\Gamma) $  is then defined formally as the `\textbf{Spec}' of $\,\O[\DRep_G(\Gamma)\,$, i.e. the simplicial
algebra $\,\O[\DRep_G(\Gamma)]\,$ viewed as an object of the opposite category $ \Ho(\scAlg_k)^{\rm op} $.
The homotopy groups of $ \O[\DRep_G(\Gamma)] $  depend
only on $ \Gamma $ and $G$, with $ \pi_0 \O[\DRep_G(\Gamma)] $ being canonically isomorphic to $ \pi_0(\Gamma)_G $. In particular, if $ \Gamma $ is a discrete simplicial group, then $ \pi_0 \O[\DRep_G(\Gamma)] \cong \Gamma_G $. Extending our terminology from \cite{BKR, BFPRW}, we will refer to $ \pi_* \O[\DRep_G(\Gamma)] $ as the {\it representation homology} of $ \Gamma $ in $G$ and denote it $ \HR_\ast(\Gamma, G) $.
We should mention that representation homology of associative and Lie algebras was introduced and studied in \cite{BKR, BR, BFPRW}. The idea of deriving the representation functor was motivated by noncommutative geometry, where the representation functor plays an important role  (see \cite{G2,KR} and also \cite{BFR}).

Next, we recall that the model category $ \sGr $ of simplicial groups is Quillen
equivalent to the category of reduced simplicial sets, $ \sset_0 $, which is, in turn, Quillen equivalent to the category $ {\tt Top}_{0,*} $ of pointed connected topological spaces. These classical equivalences are given by two pairs of adjoint functors:
$$
\lgr : \sset_0 \rightleftarrows \sGr :\overline{W}\,, \qquad |\,\mbox{--}\,| : \sset_0 \rightleftarrows {\tt Top}_{0,*} :\, \overline{S}\,,
$$
the construction of which will be briefly reviewed in Section~\ref{klg}. Here, we only recall that $ \lgr $ is
the Kan loop group functor that assigns to a reduced simplicial set $ X \in \sset_0 $ a semi-free
simplicial group $\,\lgr X\,$, which is a simplicial model of the based loop space $ \Omega |X| $ (see \cite{Kan1}). The  Kan loop group functor preserves weak equivalences and hence induces a functor
between the homotopy categories: $\,\lgr : \Ho(\sset_0) \to \Ho(\sGr) \,$. Combining this last functor with  \eqref{rep2}, we set $\,\O[\DRep_G(X)] := \L(\lgr X)_G \,$ and define  the {\it representation homology} of $ X \in \sset_0 $ by
\begin{equation}
\la{repX}
\HR_*(X, G) :=  \pi_* \O[\DRep_G(X)]\ .
\end{equation}
By definition, $ \HR_*(X, G) $ is a graded commutative algebra that depends only on
the homotopy type of $X$ and hence is a homotopy invariant of the corresponding space
$ |X| $. In degree zero, we have $\,\HR_0(X,G) \cong (\pi_1(X))_G = \O[\Rep_G(\pi_1(X))]\,$, where $ \pi_1(X)  $ is
the fundamental group of $X$. To avoid confusion, we emphasize that $ \HR_*(X, G) \not\cong
\HR_*(\pi_1(X), G) $  in general; however, if $ \Gamma $ is a discrete group and $X$ is a $ K(\Gamma, 1)$-space ({\sl e.g.},
$ X = \B\Gamma $), then we do have a natural isomorphism $\,\HR_*(X, G) \cong \HR_\ast(\Gamma, G) \,$, so there is no ambiguity in our notation.

The goal of the present paper is threefold. First, we establish basic properties of the derived representation functor \eqref{rep2}.
Second, we give an elementary construction of representation homology in terms of classical (abelian) homological algebra. Our construction is analogous to Pirashvili's construction of higher order Hochschild homology, and it provides a natural interpretation of representation homology as functor homology. This opens up the way to efficient
computations and places representation homology in one row with other classical invariants such as Hochschild and cyclic homology. Third, we construct some spectral sequences and natural maps relating representation homology to other homology theories associated with spaces (including the Pontryagin algebra
$ \H_\ast(\Omega X) $, higher Hochschild homology and stable homology of the automorphism groups of f.g. free groups $ \bF_n $).
We also compute representation homology explicitly in a number of interesting cases, including the spheres $ \bS^n $, suspensions $ \Sigma X $, co-$H$-spaces, closed surfaces of arbitrary genus and some classical 3-dimensional spaces, such as the link complements in $ \R^3 $ and the lens spaces $ L(p,q) $.  In our subsequent paper, \cite{BRY2}, we will extend these computations to arbitrary simply connected topological spaces by expressing the representation homology of a 1-connected space of finite rational type in terms of its Quillen and Sullivan models, and give
some applications to representation theory.

\subsection{Main results}
We now proceed with a summary of the main results of the paper. Recall that an affine algebraic group $G$ is defined by its functor of points, which is a group-valued representable functor
on commutative algebras. This functor extends in the natural way  to simplicial commutative algebras:
\begin{equation}
\la{Gs}
G:\ \scAlg_k \to \sGr\ ,\quad A_\ast \mapsto G(A_\ast)\ .
\end{equation}
By definition, the representation functor \eqref{rep1} is left adjoint to the functor of points
of $G$, hence its simplicial extension is left adjoint to \eqref{Gs}. Thus, for any affine algebraic group, we have the adjunction
\begin{equation}
\la{sAdj}
(\,\mbox{--}\,)_G\,: \,\sGr\, \rightleftarrows \, \scAlg_k\, :\, G\ .
\end{equation}
Our first main result reads
\bthm
\la{Thm0}
The functor \eqref{Gs} has a total right derived functor $ \bR G: \Ho(\scAlg_k) \to \Ho(\sGr) $,  which is right adjoint to the derived representation functor \eqref{rep2}: thus,
\eqref{sAdj} induces the derived adjunction
\begin{equation}
\la{derAdj}
\L(\,\mbox{--}\,)_G\,: \,\Ho(\sGr)\, \rightleftarrows \, \Ho(\scAlg_k)\, :\, \bR G
\end{equation}
\ethm

Note that the categories $\sGr$ and $\scAlg_k$ have natural (simplicial) model structures, and the above result would be immediate from the well-known Adjunction Theorem of \cite{Q1} 
if \eqref{sAdj} were Quillen functors.
However, it is easy to see  that the functors \eqref{sAdj} do {\it not} form a Quillen pair of model categories, nor even do they form a deformable adjunction
of homotopical categories in the sense of \cite{DHKS}. Theorem~\ref{Thm0} is therefore an interesting and fairly nontrivial result, which is -- to the best of our knowledge -- new.

In the special case when $ G = \GL_n $, the derived functor $ \bR G $ can be described explicitly as the composite of two well-known functors:
\begin{equation}
\la{RGL}
\bR \GL_n \cong \L l \circ \widehat{\GL}_n\ .
\end{equation}
The functor $ \widehat{\GL}_n: \scAlg_k \to \sMon $ takes values in the category of
simplicial monoids, assigning to a simplicial algebra $ A_\ast $ the simplicial monoid
of $ (n \times n)$-matrices over $A_*$  `invertible up to homotopy':
more precisely,  $ \widehat{\GL}_n(A_\ast) $ is defined by the pull-back diagram in the category
$ \sMon $:
\begin{equation}
\la{waldfun}
\begin{diagram}[small]
\widehat{\GL}_n(A_\ast) & \rTo^{} &  {\GL}_n(\pi_0 A_\ast)    \\
\dTo^{}                 &         & \dTo^{}\\
M_n(A_\ast)      &  \rTo^{}       & M_n(\pi_0 A_\ast)
 \end{diagram}
\end{equation}
This functor was originally introduced by F. Waldhausen \cite{Wa} to define a homotopy invariant version of algebraic $K$-theory of simplicial rings. The second functor in \eqref{RGL} is the total left
derived functor $ \L l: \Ho(\sMon) \to \Ho(\sGr) $ of the group completion (localization) of simplicial monoids: it can be viewed as a special case of the classical Dwyer-Kan localization of simplicial categories studied in \cite{DK}. Formula \eqref{RGL} is rather unusual as it expresses
a right derived functor in terms of a left derived one.

For an arbitrary algebraic group $G$, we construct an explicit model for $ \bR G $ using the recent work of S. Galatius and A. Venkatesh \cite{GV}. In this model, instead of simplicial monoids, we factor
$ \bR G  $ through the reduced simplicial spaces (or reduced Segal precategories) in the sense of J.
Bergner \cite{Be07}. This construction leads to a more general definition of representation homology
which applies to simplicial spaces and does not use the Kan loop group equivalence (see Section~\ref{S34}).

The second main result of this paper is an interpretation of representation homology in terms of
classical homological algebra in functor categories. To this end
we consider the category $ \ffgr $ of finitely generated free groups with objects $ \rn := \bF\langle x_1, x_2, \ldots, x_n \rangle $ (for $ n \ge 0 $) and morphisms being the arbitrary group homomorphisms.
This category is a PROP (i.e., a small permutative category) with monoidal structure $ \rn \boxtimes \langle m \rangle = \langle n+m\rangle $,
and it is known that the category of $k$-algebras over $\ffgr $ (i.e., the category of strict monoidal functors from $ \ffgr $ to the category $ \ve_k $ of $k$-vector spaces) is equivalent to the category of commutative Hopf $k$-algebras\footnote{A direct combinatorial proof of this fact can be found, e.g., in \cite{Ha}.}. Under this equivalence, a commutative Hopf algebra $ {\mathcal H} $ corresponds to the functor $ \umH: \ffgr \to \ve_k $, $\,\rn \mapsto \mH^{\otimes n}\,$, which actually takes its values in the category of commutative algebras.
Note that any functor $ F: \ffgr \to \cAlg_k $ extends naturally (by taking the left Kan extension along the inclusion $ \ffgr \into \fgr$) to the category of all (based) free groups, $ \fgr $, whose objects are the free groups $ \bF\langle S \rangle $ given with a prescribed generating set $S$ and morphisms are the arbitrary group homomorphisms; we denote this extension
by $ \tilde{F}: \fgr \to \cAlg_k $.

Now, mimicking the Pirashvili construction of higher Hochschild homology
({\it cf.} \cite{P1} and Section~\ref{S22} below), for a reduced simplicial set $ X \in \sset_0 $
and a commutative Hopf algebra $ \mH $, we consider the composition of functors
$$
\Delta^{\rm op} \xrightarrow{\,\lgr X\,} \fgr \xrightarrow{\,\tilde{\umH}\,} \cAlg_k\ ,
$$
where $ \lgr X $ is the Kan loop group construction of $X$ and $\tilde{\umH} $ is the (left Kan) extension of the strict monoidal functor $  \umH: \ffgr \to \ve_k $ corresponding to $ \mH$. This defines a simplicial commutative  algebra $ \umH(\lgr X) $,  whose homotopy groups we denote by
\begin{equation}
\la{HRH1}
\HR_\ast(X, \mH) := \pi_\ast \umH(\lgr X) = \H_\ast[N(\umH(\lgr X))] \ .
\end{equation}
It turns out that this definition is {\it equivalent} to our original definition of representation homology \eqref{repX} given in terms of the  derived representation functor $ \L(\,\mbox{--}\,)_G $.  Precisely ({\it cf.} Proposition~\ref{rephom0}),
\bprop
\la{IProp1}
Let $ G $ be an affine group scheme over $ k $ with coordinate ring $ \mH = \O(G) $. Then,
for any $ X \in \sset_0 $, there is a natural isomorphism of graded commutative algebras
\begin{equation*}
\la{eqhr}
\HR_\ast(X, \O(G)) \cong \HR_\ast(X, G)\ .
\end{equation*}
\eprop
Thanks to Proposition~\ref{IProp1}, we may (and will) use the notation $ \HR_\ast(X, G) $ and $ \HR_\ast(X, \mH) $  interchangeably, without causing confusion. Although its proof is almost immediate, Proposition \ref{IProp1} has a number of important implications. First, we state
the following theorem, which is the main result of Section~\ref{section_rephom} (see Theorem~\ref{rephomdtp2}).
\bthm
\la{IThm1} For any $ X \in \sset_0 $, there is a natural first quadrant spectral sequence
\begin{equation}
\la{ssrh}
E^2_{pq}\,=\, \mathrm{Tor}^{\ffgr}_p(\underline{\H}_{\,q}(\Omega X;k), \umH)\ \underset{p}{\Longrightarrow}\  \mathrm{HR}_{n}(X, \mH)\
\end{equation}
converging to the representation homology of $X$.
\ethm
The spectral sequence \eqref{ssrh} relates the representation homology $\HR_*(X, \mH) $ of a space $X$ to its Pontryagin algebra  $ \H_\ast(\Omega X; k) $. To describe the $ E_2$-term of \eqref{ssrh} we recall that $ \H_\ast(\Omega X; k) $
has a natural structure of a graded cocommutative Hopf algebra with coproduct induced by the Alexander-Whitney diagonal
and the product by the Eilenberg-Zilber map. For each $ q \in \Z $, the assignment $\,\rn \mapsto [\H^{\otimes n}]_q\,$, where $ [\H^{\otimes n}]_q $ is the $q$-th graded component of the $n$-th tensor power of $ \H = \H_\ast(\Omega X; k) $, defines a functor
$ \underline{\H}_{\,q}: \ffgr^{\rm op} \to \ve_k \,$, which is the first argument of the ``Tor'' in \eqref{ssrh}. The ``Tor'' itself is the (abelian)
derived functor of the  tensor product $ \, \otimes_{\ffgr} \,$ between  covariant and contravariant  $\ve_k$-valued functors over the (small) category $ \ffgr $. The spectral sequence \eqref{ssrh} is a counterpart of Pirashvili's fundamental spectral sequence for higher Hochschild homology ({\it cf.} \cite[Theorem~2.4]{P1}); however, in the case of
representation homology it takes a more geometric form.

Theorem~\ref{IThm1} has several interesting implications. First of all, it shows that the representation homology $ \HR_\ast(X, \mH) $ is stable under Pontryagin equivalences\footnote{that is, maps of spaces $ X \to Y $ inducing isomorphisms  of Pontryagin algebras
$\, \H_\ast(\Omega X; k) \stackrel{\sim}{\to} \H_\ast(\Omega Y; k) \,$.}, and hence, if $X$ is simply-connected, $ \HR_\ast(X, \mH) $ is actually a {\it rational}\, homotopy invariant of $X$
(see Proposition~\ref{rephomhtpinv}).
Next, if $X$ is a $ K(\Gamma,1)$-space, the spectral sequence \eqref{ssrh} degenerates giving an isomorphism ({\it cf.} Corollary~\ref{bgamma})
\begin{equation}
\la{torgamma}
\HR_{\ast}(\Gamma, G) \cong \Tor_{\ast}^{\ffgr}(k[\Gamma], \O(G))\ ,
\end{equation}
where $ k[\Gamma] $ is the group algebra of $ \Gamma $ equipped with
the natural (cocommutative) Hopf algebra structure. The isomorphism \eqref{torgamma} shows that the representation homology has a natural `Tor' interpretation, similar to the classical (Connes) interpretation of the Hochschild and cyclic homology (see \cite[Chap.~6]{L}). It is also interesting to compare \eqref{torgamma} with another natural isomorphism
\begin{equation}
\la{torgamma1}
\H_{\ast+1}(\Gamma, k) \cong \Tor_{\ast}^{\ffgr}(k[\Gamma], \lin_k)\ ,
\end{equation}
which provides a `Tor' interpretation (over $ \ffgr$) for the ordinary homology of $ \Gamma $
as a discrete group. Here $ \lin_k $ stands for the linearization functor
$ \ffgr \to \ve_k $ that takes the free group $ \rn $ to the vector space $ k^n $
(see \eqref{linfun}). Note that if $ G  = {\mathbb G}_a $ is the additive group over $k$,
then we have a natural isomorphism of functors
$ \O({\mathbb G}_a) \cong {\rm Sym}_k(\lin_k) $, which implies
$\, \HR_\ast(\Gamma, {\mathbb G_a}) \cong \Sym_k[\,\H_{\ast+1}(\Gamma, k)] \,$. More generally,
for any (pointed connected) space $X$, we have an isomorphism of graded commutative
algebras (see Example~\ref{addgr})
\begin{equation}
\la{hrvsh}
\HR_\ast(X, {\mathbb G_a}) \cong  \Sym_k[\,\H_{\ast+1}(X, k)]\ ,
\end{equation}
where $ \Sym_k[\,\H_{\ast+1}(X, k)] $ is the symmetric algebra of the graded
vector space $ \H_{\ast+1}(X, k) = \oplus_{i \ge 0}\, \H_{i+1}(X, k) $. Thus, we may
think of representation homology as a generalization of the ordinary (singular) homology
of spaces.

In Section~\ref{sect5}, we show that representation homology can be also
viewed as a generalization of higher Hochschild homology of spaces.
The main result  of this section reads
({\it cf.} Theorem \ref{repvshh1} and Theorem~\ref{repvshh}):
\bthm
\la{MT51}
Let $ \mH $ be a commutative Hopf algebra.

$(a)$ For any simplicial set $\,X \in \sset\,$, there is a natural isomorphism
\begin{equation}
\la{sihrhh}
\HR_\ast(\Sigma(X_{+}), \mH)\,\cong\,\HH_\ast(X,\mH)\ ,
\end{equation}
where $ X_{+} = X \sqcup \{\ast\} $ is a pointed simplicial set obtained
from $X$ by adjoining functorially a basepoint, and $\, \Sigma \,$ is
the $($reduced$)$ suspension functor on the category of pointed simplicial sets.

$(b)$ For any pointed simplicial set $\, X \in \sset_* \,$, there is a natural isomorphism
\begin{equation}
\la{sihrhh1}
\HR_\ast(\Sigma X, \mH)\,\cong\,\HH_\ast(X,\mH; k)\ ,
\end{equation}
where $ \HH_\ast(X,\mH; k) $ is the Pirashvili-Hochschild homology of the commutative algebra $ \mH $ with
coefficients in $ k $ viewed as an $\mH$-module via the Hopf algebra counit $ \varepsilon: \mH \to k $.
\ethm
The proof of Theorem~\ref{MT51}  is based on Milnor's classical  $FK$-construction \cite{Mi} that gives a simple
simplicial group model for the space $ \Omega \Sigma |X| $.

Theorem~\ref{MT51}  has strong implications: in particular, it allows one to compute the representation homology of suspensions in a completely explicit way. It is known that $ \Sigma X $ for any pointed connected space $X$ is {\it rationally} homotopy equivalent to a bouquet of spheres of dimension $ \ge 2\,$. Since representation homology depends only on the rational homotopy type of a space, the isomorphism \eqref{sihrhh1} together with Pirashvili's computations \cite{P1} of higher Hochschild homology of spheres, implies ({\it cf.} Proposition~\ref{rephomsusp})
\bprop
\la{Psusp}
For any pointed connected space $X$ of finite type, there is an isomorphism
\begin{equation*}
\la{psusp}
\HR_{\ast}(\Sigma X, G)\,\cong\, \Sym_k [\,\rH_\ast(X; \g^{\ast})]\ .
\end{equation*}
where $\, \Sym_k [\,\rH_\ast(X; \g^{\ast})] \,$ is the graded symmetric algebra of the reduced $($singular$)$ homology of $X$ with coefficients in the dual Lie algebra of the group $G$.
\eprop
\noindent
By induction, Proposition~\ref{Psusp} implies
$\, \HR_{\ast}(\Sigma^n X, G) \cong  \Sym_k \left( \rH_\ast(X; \g^{\ast})[n-1] \right)\,$
for all $\, n \ge 1 $. In particular, for $ \bS^n \cong \Sigma \bS^{n-1} $, we have
\begin{equation}
\la{IHRsphers}
 \HR_\ast(\bS^n, G) \,\cong \, \Sym_k (\g^*[n-1])\ ,\ n \ge 2 \ .
\end{equation}

In Section~\ref{sec7}, we compute representation homology of some classical non-simply connected spaces.
Our examples include closed surfaces of arbitrary genus (both orientable and non-orientable) as well as some
three-dimensional spaces (the link complements in $ \R^3 $ and $ \bS^3 $, the lens spaces $ L(p,q) $ and a general
closed orientable $3$-manifold). The representation homology of surfaces and link complements is expressed
in terms of classical Hochschild homology of $ \O(G) $ and related commutative algebras. For example, for the link
complements in $ \R^3 $, we prove ({\it cf.} Theorem~\ref{hrlink})
\bthm
\la{Ihrlink}
Let $ L $ be a link in $ \R^3 $ obtained as the Alexander closure of a braid $ \beta \in B_n $.
Then the representation homology of the complement of its $($regular$)$ neighborhood in $ \R^3 $ is given by
\begin{equation}
\la{IHRRL}
\HR_\ast(\R^3\!\setminus\! L, G) \,\cong \,\HH_\ast(\O(G^n),\,\O(G^n)_{\beta})\ .
\end{equation}
\ethm
The right-hand side of \eqref{IHRRL} is the (ordinary) Hochschild homology of the associative algebra
$ \O(G^n) $ with bimodule coefficients. The bimodule $ \O(G^n)_{\beta} $ is isomorphic
to $ \O(G^n) = \O(G)^{\otimes n} $ as a left module, while the right action of $ \O(G^n)$  is twisted by an
element $ \beta $ viewed as an automorphism of $ \O(G)^{\otimes n} $ via the Artin representation of the braid group $B_n$.

Theorem~\ref{Ihrlink} shows that the Hochschild homology groups $\,\HH_\ast(\O(G^n),\,\O(G^n)_{\beta})\,$
are algebraic invariants of links in $\R^3$, which, to the best of our knowledge, have not appeared
in the earlier literature. We should mention, however, that the representation homology of link complements bears a striking resemblance to {\it knot contact homology}, which is a new geometric homology theory of knots and links in $ \R^3 $ defined in \cite{Ng1} and studied extensively in recent years (see Remark after Theorem~\ref{hrlink}). We will discuss the relation between representation homology and knot contact homology in our subsequent paper.

In Section~\ref{sec7}, we also discuss a multiplicative version of the
derived Harish-Chandra Conjecture proposed in \cite{BFPRW}. If $G$ is a connected reductive
group with a maximal torus $ T \subset G $ and $W$ is the associated Weyl group, then for any space $X$, there is a natural map
\begin{equation}
\la{IHCmap}
\HR_{\ast}(X,G)^G \to \HR_\ast(X, T)^W\ ,
\end{equation}
which we call the {\it derived Harish-Chandra homomorphism} ({\it cf.} \cite[Section~7]{BFPRW}).
In view of \eqref{IHRsphers},  by the classical Chevalley Restriction Theorem \cite{Ch}, the map  \eqref{IHCmap} is an isomorphism for any odd-dimensional sphere $ X= \bS^{2p+1} $.
%
%
We conjecture that \eqref{IHCmap} is also an isomorphism for the two-dimensional torus
$ \bT^2 = \bS^1 \times \bS^1 $, which gives the following explicit formula for the
representation homology of $  \bT^2  $ (see Section~\ref{S7.1.1}, Conjecture~\ref{conj2}):
\begin{equation*}
 \HR_*(\bT^2, G)^G\, \cong\, [\O(T \times T) \otimes \Lambda^*_k(\h^*)]^W\ ,
\end{equation*}
where $ \h $ is the Lie algebra of $T$ (i.e., a Cartan subalgebra of $ \g $) and $ \h^* $ is its linear dual.
As for the Drinfeld homomorphism, it would be interesting to find more examples of spaces, for which the map
\eqref{IHCmap} is an isomorphism, and/or give an abstract characterization of all such spaces.

In the last section of the paper, we give another interpretation of representation homology
as the Hochschild-Mitchell homology of a certain bifunctor on the category of finitely generated free groups
$ \mathfrak{G} $. Such an interpretation is useful for several reasons. First, it allows us to define
representation {\it cohomology} in a natural way (by simply replacing the Hochschild-Mitchell homology with
the Hochschild-Mitchell cohomology of the same bifunctor). Second, it suggests that it is
natural to extend the definition of representation (co)homology by taking the Hochschild-Mitchell
(co)homology of $ \mathfrak{G} $ with coefficients in an arbitrary bifunctor $ D\,$: i.e., $\, \HR(D) :=
\HH(\mathfrak{G}, D) $. Third and most important, it exhibits a close analogy with  {\it topological} Hochschild
homology, which is known to be isomorphic to the Hochschild-Mitchell homology of the category  $ \mathfrak{G}_{\rm ab} $
of finitely generated free {\it abelian} groups (see \cite{PW}). Motivated by this analogy, we construct functorial trace maps
$$
{\rm DTr}_{n}^{\mathfrak G}(D):\
\H_\ast(\Aut(\bF_n),\,D_{n})\, \to\,   \HR_\ast(D)\ ,\quad \forall\,n \ge 1\ ,
$$
relating homology of the automorphism groups of f.g. free groups with appropriate coefficients to
representation homology. These maps are compatible with natural inclusions $ \Aut(\bF_n) \into \Aut(\bF_{n+1}) $
and hence have an stable limit as $ n \to \infty \,$. The corresponding stable map
$\, {\rm DTr}_{\infty}^{\mathfrak G}(D):\,\H_\ast(\Aut_{\infty},
\,D_{\infty})\, \to\,   \HR_\ast(D)\,$ can be viewed as a non-abelian analogue of the
classical Dennis trace relating topological Hochschild homology to stable homology of general linear groups.
We conjecture that the map $ {\rm DTr}_{\infty}^{\mathfrak G}(D) $ is actually
an isomorphism, whenever $ D $ is a {\it polynomial} bifunctor ({\it cf.}
Conjecture~\ref{conj3}). This is a non-abelian analogue of
a theorem of Scorichenko \cite{FFPS}.

\subsection{Relation to derived algebraic geometry} The derived representation schemes $ \DRep_G(X) $ are basic objects of derived algebraic geometry. To the best of our
knowledge, the first construction of this kind -- the derived moduli space $ \bR {\rm Loc}_G(X) $ of $G$-local systems over a finite, pointed, connected CW complex $X$ -- was proposed by M. Kapranov in \cite{K}. He defined $ \bR {\rm Loc}_G(X) $ using a simplicial DG scheme $ \bR BG $ which played the role of a canonical `injective resolution' of the classifying space $ BG $ of the algebraic group $G$ in the category of simplicial DG schemes. A more refined construction $ \textbf{Map}(X,BG) $ -- called the derived mapping stack of flat $G$-bundles on $X$ -- was developed by B.~To\"en and G. Vezzosi in \cite{TV08} (see also \cite{PTVV13}), using  local homotopy theory of simplicial presheaves on the category of (derived) affine schemes. For a detailed comparison of these two constructions with our construction of $ \DRep_G(X) $,
we refer the reader to the appendix of \cite{BRY}, where we showed that -- despite different frameworks -- all three constructions are essentially equivalent.

We would like to conclude this introduction by mentioning some interesting topological generalizations of higher Hochschild homology that appeared in recent years, such as factorization homology (see, e.g., \cite{G, GTZ}) and higher topological Hochschild homology  \cite{BCD}. Our results show that representation homology, while closely related to Hochschild homology, is a richer and somewhat more geometric theory that blends topology and representation theory in a very natural way. It would therefore be  interesting to see if representation homology admits topological refinements similar to those of Hochschild homology.
\subsection{Appendix} The paper contains an appendix, where we collect basic facts and prove 
some new results in abstract homotopy theory concerning derived functors. The main result of the
appendix --- Theorem~\ref{ThA2} --- arises from our attempt to abstract the situation of 
Theorem~\ref{ThAdj}: it is a version of Quillen's Derived Adjunction Theorem for homotopical categories. This theorem as well as Theorem~\ref{ThA3} and Lemma~\ref{LeA1} 
are of independent interest.

\subsection{Outline of the paper} The paper is organized as follows.
In Section~\ref{S2}, we introduce notation and recall some basic facts
about simplicial sets and spaces. In Section~\ref{S3}, we study basic
properties of the derived representation functor
and define representation homology. In Section~\ref{section_rephom},
we give our second construction of representation homology in terms
of functor homology and derive its implications. In Section~\ref{section_rephom},
we establish the isomorphism between the representation homology of
suspensions and  higher Hochschild homology. In Section~\ref{sec7},
we give examples computing representation homology explicitly for some geometrically
interesting spaces. In Section~\ref{sec8}, we identify representation homology
in terms of Hochschild-Mitchell homology  and construct a non-abelian analogue of the Dennis trace map relating representation homology to the stable homology of automorphism groups of finitely generated free groups. The paper ends with an appendix where we recall basic definitions and prove a few results from abstract homotopy theory used in Section~\ref{S3}.

\subsection*{Acknowledgements}
We are very grateful to  K. Hess, M. Khovanov, A. Lindenstrauss, B. Richter and V. Turaev
for interesting questions and comments. We are particularly indepted to C. Berger for his explanations concerning the Kan loop group construction and
for providing us with reference \cite{Berg}. We also thank J. Davis, S. Gadgil, M. Kassabov, M. Mandell and S. Patotski for useful discussions.
The second author expresses his gratitude to the Department of Mathematics, Indian Institute of Science, Bengaluru for conducive working conditions during his visit in the summer of $2016$. Research of the first two authors was partially supported by the Simons Foundation Collaboration Grant 066274-00002B as well as NSF grants DMS 1702323 and DMS 1702372.

\section{Preliminaries}
\la{S2}

In this section, we introduce notation and recall some basic definitions related to
simplicial sets. Standard references for this material are \cite{M}, \cite{GJ} and \cite[Chapter 8]{W}.

\subsection{Simplicial objects} \la{simpobj}
Let $\Delta$ denote the simplicial category. Recall that the objects of $\Delta$ are the finite ordered sets $[n]\,:=\,\{0, 1, \ldots, n\}$, $n \geq 0$, and the morphisms are the  (weakly) order preserving maps $[n] \to [m] $.  A {\it simplicial object} in a category $\mathscr{C}$ is a contravariant functor from $\Delta$ to $\mathscr{C}$: i.e.,
$\,\Delta^{\mathrm{op}} \rar \mathscr{C} \,\text{.}$
The simplicial objects in $\mathscr{C}$ form a category, with morphisms being the natural transformations of functors. We denote this category by $s\mathscr{C}$. If $X\,\in\, {\rm Ob}(s\mathscr{C})$, we write $X_n := X([n])$.

The category $\Delta$ is generated by two distinguished classes of morphisms $\{ \delta^i\}^{n \geq 1}_{0 \leq i \leq n}$ and $\{\sigma^j\}^{n \geq 0}_{0 \leq j \leq n}$, whose images under $X\,\in\,s\mathscr{C}$ are called the {\it face} and {\it degeneracy maps} of $X$, respectively. The map $\delta^i\,:\, [n-1] \rar [n]$ is the (unique) injection that does not contain ``$\,i\,$'' in its image; the corresponding face map is denoted by $d_i\,:=\,X(\delta^i)\,:\,X_n \rar X_{n-1}$. Similarly, for $n \geq 0$, the map $\sigma^i\,:\,[n+1] \rar [n]\,$ is the (unique) surjection in $\Delta$ that takes
value ``$\,i\,$'' twice. The image of $\sigma^i$ under $X$ is the degeneracy map $s_i\,:=\,X(\sigma^i)\,:\,X_n \rar X_{n+1}$. The face and degeneracy maps of a simplicial object satisfy the following {\it simplicial relations}\,:
\begin{align}
 & d_id_j \,=\, d_{j-1}d_i \,\,\,\,\text{ if } i<j \, \nonumber \\
 & d_is_j\,=\,s_{j-1}d_i\,\,\,\, \text{ if } i<j\,  \nonumber\\
 & d_is_j\,=\,s_jd_{i-1} \,\,\,\,\text{ if } i>j+1 \,  \la{sids} \\
 & s_is_j\,=\,s_{j+1}s_i \,\,\,\,\text{ if } i \leq j\,  \nonumber\\
 & d_is_j\,=\,\id \qquad\ \text{ if } i=j\,,\,j+1\,\text{.}  \nonumber
\end{align}
Thus, a simplicial object in $ s\mathscr{C}$ is determined by a family $X = \{X_n\}_{n \geq0}$ of objects in $\mathscr{C}$ together with morphisms $d_i\,:\,X_n \rar X_{n-1}$ and $s_j\,:\,X_n \rar X_{n+1}$ satisfying the relations~\eqref{sids}.
The object $X_n$ is usually called the ``set'' of {\it $n$-simplices} of $X$, and the $0$-simplices are usually
called the {\it vertices} of $X$.

We let $\sset $ denote the category of simplicial sets ({\it i.e.}, simplicial objects in the category $\mathtt{Set}$).  A simplicial set $ X $ is called {\it reduced} if it has a single vertex, i.e. $X_0 = \{\ast \}\,$.
The full subcategory of $ \sset $ consisting of reduced simplicial sets will be denoted $\sset_0$. A simplicial set $X $ is called {\it pointed} if there are distinguished simplices $ x_n \in X_n $, one in each degree, such that $x_n\,=\,s_0(x_{n-1})$ for all $n \geq 1$. The sequence
$ (x_0, x_1, x_2, \ldots) \in \prod_{n \ge 0} X_n $ is called a basepoint of $ X $.
The category of pointed simplicial sets will be denoted $\sset_{\ast}$. Note that $\sset_0 $ can also be viewed as a full subcategory of $\sset_{\ast} $ as every reduced simplicial set has a canonical (unique) basepoint.

Given $X\,\in\,\sset$, the set of {\it nondegenerate} $n$-simplices of $X$ is defined to be
$$ \overline{X}_n\,:=\,X_n \setminus \bigcup_{i=0}^{n-1} s_i(X_{n-1})\,\text{.}$$
Every element of $X_n$ can be uniquely expressed in terms of the nondegenerate elements of $X$ (see~\cite[Lemma 11]{GM} for a precise statement). In particular, a simplicial set can be defined by specifiying its nondegenerate simplices together with the restriction of each face map to the set of nondegenerate simplices.

\vspace{1ex}

\noindent  We give a few basic examples of simplicial
sets that will be used in this paper.

\subsubsection{Discrete simplicial objects}

To any object $A\,\in\,\mathscr{C}$ one can associate a simplicial object $ A_{\ast}\,\in\,s\mathscr{C}$, with $A_n=A$ and $d_i$, $ s_j$ being the identity map of $A$ for all $n,i,j$. This gives a fully faithful embedding $\mathscr{C} \hookrightarrow s\mathscr{C}$. The objects of $s\mathscr{C}$ arising this way are called {\it discrete simplicial objects}.

\subsubsection{Geometric simplices}
\la{geomsimp}

The $n$-dimensional geometric simplex is the topological space
 $$\Delta^n\,:=\, \{ (x_0,\ldots,x_n)\,\in\,\mathbb R^{n+1}\,|\, \sum_{i=0}^n x_i=1 \,,\, x_i \geq 0\,\}\,\text{.}$$
Let $e_i$ denote the vertex of $\Delta^n$ with $i$-th coordinate $1$. For any morphism $f\,:\,[m] \rar [n]$ in $\Delta$, there is a (unique) linear map $\R^{m+1} \rar \R^{n+1}$ sending $e_i$ to $e_{f(i)}$, that restricts to a map of topological spaces $ f^*\,:\,\Delta^m \rar \Delta^n$. The collection $\Delta^{\ast}\,:=\,\{\Delta^n\}_{n \geq 0}$ forms a {\it cosimplicial space}, {\it i.e.}, a (covariant) functor $\Delta \rar \mathtt{Top}$, where  $\mathtt{Top}$ denotes the category of (compactly generated weakly Hausdorff) topological spaces. This functor is faithful: it gives a topological realization of the simplicial category, which was historically the first definition of $\Delta $.

\subsubsection{Standard simplices} \la{sstd}

Let $\, Y: \Delta \into \sset \,$ denote the Yoneda embedding. The functor $ Y $ assigns to $ [n] $  a simplicial set  $\Delta[n]_{\ast} $ called the {\it standard $n$-simplex}. Explicitly, $\Delta[n]_{\ast}$ is given by
 $$\Delta[n]_k\,:=\,\Hom_{\Delta}([k],[n])\,\cong\, \{(n_0,\ldots,n_k)\,|\,0 \leq n_0 \leq \ldots \leq n_k \leq n\} \,,$$
where a function $f:\,[k] \rar [n]$ is identified with the sequence of its values $(f(0),\ldots,f(k))$. Under this identification, the nondegenerate simplices correspond to {\it strictly} increasing functions, and the face and degeneracy maps in  $\Delta[n]_{\ast}$ are given by
$$
d_i(n_0,\ldots,n_k)\,=\,(n_0,\ldots,\hat{n}_i,\ldots,n_k)\,,\,\,\,\, s_j(n_0,\ldots,n_k)\,=\,(n_0,\ldots, n_j, n_j, \ldots, n_k)\ .
 $$
By Yoneda Lemma, for any simplicial set $X$, there is a natural bijection
$$
\Hom_{\sset}(\Delta[n]_{\ast}, X)\,\cong\, X_n\,,
$$
which shows that $ \Delta[n]_{\ast}$ (co)represents the functor: $\sset \rar \set\,,\,\,X \mapsto X_n$.

\subsubsection{Simplicial spheres} \la{sscirc}
The Yoneda functor $ Y: \Delta \rar \sset$ can be also regarded as a cosimplicial object in the category of simplicial sets. In particular, for any $n \geq 1$, there are $ n+1 $ coface maps $\,d^i :\,\Delta[n-1]_{\ast} \rar \Delta[n]_{\ast} $, $\, 0 \leq i \leq n$. Using these maps, we define the {\it boundary} of $\Delta[n]_{\ast}$ to be the simplicial subset
 $$\partial \Delta[n]_{\ast}\,:=\, \bigcup_{ 0 \leq i \leq n} d^i(\Delta[n-1]_{\ast}) \,\subset \,\Delta[n]_{\ast}\ , $$
The {\it simplicial $n$-sphere} is then defined to be the corresponding quotient set $\mathbb S^n_{\ast}\,:=\,\Delta[n]_{\ast}/\partial \Delta[n]_{\ast}$. It is easy to see that the only nondegenerate simplices in $\mathbb S^n_{\ast}$ are in degree $0$ and $n$, with $ \overline{\mathbb S}^n_0\,=\,\{\ast\}$ and $\overline{\mathbb S}^n_n\,=\,\{S\}$,
 where $S$ is the image of the map $\id\,\in\,\Delta[n]_n$ in $\mathbb S^n_n$. Note that $d_i(S)=s_0^{n-1}(\ast)$ for all $i$. Thus, the simplicial structure of $\mathbb S^n_{\ast}$ reflects the standard CW decomposition of the $n$-sphere $\mathbb S^n$ with one cell in dimension $0$ and one cell in dimension $n$.

The simplicial $1$-sphere $ \mathbb S^1_{\ast} $ is called the {\it simplicial circle}. By Example~\ref{sstd}, we have
$ \Delta[1]_k\,\,\cong\,\{(\underbrace{0,\ldots,0}_{\text{$i$ }},\underbrace{1,\ldots,1}_{\text{$k+1-i$ }})\,|\, i=0,1,\ldots,k+1\}$ and $ \partial \Delta[1]_k\,=\,\{(0,\ldots 0),(1,\ldots,1)\}$. Hence, $ \mathbb S^1_{\ast} $ is
given explicitly by
 $$\mathbb S^1_k \,\cong\, \{(\underbrace{0,\ldots,0}_{\text{$i$ }},\underbrace{1,\ldots,1}_{\text{$k+1-i$ }})\,|\, i=1,\ldots, k+1\}\,,$$
with $(0,\ldots,0) $ corresponding to the basepoint $ \ast $.


There is an important functor $ |\,\mbox{--}\,|:\,\sset \rar \mathtt{Top}$ assigning to each simplicial set $ X $ a topological space $ |X| $ called
the {\it geometric realization} of $X$. Explicitly, the space $ |X| $ is defined by
$$ |X| := \bigsqcup_{n \geq 0} (X_n \times \Delta^n)/\sim \,, $$
where each set $X_n$ is equipped with discrete topology and the equivalence relation is given by
\begin{align*}
 & (d_ix, p) \sim (x, d^ip) \text{ for } (x,p)\,\in\, X_n \times \Delta^{n-1}\\
 & (s_jx,p) \sim (x,s^jp) \text{ for } (x,p) \,\in \, X_{n-1} \times \Delta^{n}\,\text{.}
\end{align*}
More formally (see, e.g.,~\cite[Section 1.3]{R}), the functor $ |\,\mbox{--}\,|:\,\sset \rar \mathtt{Top}$ can be defined as the (left) Kan extension $\, |\,\mbox{--}\,| = {{\rm Lan}}_Y(\Delta^{\ast}) $ of the geometric simplex  $\,\Delta^{\ast} $ along the Yoneda embedding $ Y:
\Delta \to \sset $. It follows from this definition that $ |\Delta[n]_*| \cong \Delta^n $ for all $ n \ge 0 $,
and in general, $\, |X| \cong\mathrm{colim}\, \Delta^n $, where the colimit is taken over all morphisms of $ \Delta[n]_{\ast} \rar X\,,\,n \geq 0 $.
If $X \in \sset $ is a simplicial set and $ x_0 \in X_0 $,
we write $ \pi_n(X,x_0) $ for the {\it $n$-th homotopy group} of $X$ at $x_0$,
which is, by definition, the $n$-th homotopy group $ \pi_n(|X|, x_0) $ of the 
geometric realization of $X$.

The category $\sset$ has a standard model structure, where the  weak equivalences are the morphisms inducing weak homotopy equivalences of the corresponding geometric realizations.
The cofibrations are levelwise injective maps and the fibrations are the Kan fibrations (see~\cite[\S 7]{M}). This structure gives a model structure on $\sset_0 $.

Let $(X, \ast)$ be a pointed topological space. The {\it (total) singular complex} of $X$ is a simplicial set  $ S_{\ast}(X)$ defined by $ S_n(X)\,:=\,\Hom_{\mathtt{Top}}( \Delta^n, X)$. The {\it Eilenberg subcomplex} of $ S_{\ast}(X)$ is
$$ \overline{S}_n(X) := \{ f\,:\,\Delta^n \rar X\,:\, f(v_i)\,=\,\ast\,\,\, \text{ for all vertices } \,v_i \,\in\,\Delta^n\,\}\,\text{.}$$
If $X$ is connected, the natural inclusion $\overline{S}_{\ast}(X) \into {S}_{\ast}(X)$ is a weak equivalence of simplicial sets. Further, if we restrict $\overline{S}$ to the category  $\mathtt{Top}_{0,\ast}$ of connected pointed spaces, we get the pair of adjoint functors
\begin{equation} \la{qeq}
|\,\mbox{--}\,| \,: \mathtt{sSet}_0 \,\rightleftarrows \,\mathtt{Top}_{0,\ast}:\, \overline{S} \,\text{,}
\end{equation}
which induce mutually inverse equivalences of the homotopy categories: $\Ho(\mathtt{sSet}_0)\simeq \Ho(\mathtt{Top}_{0,\ast}) $.
This equivalence justifies the following standard convention which we will follow throughout the paper.

\vspace{2ex}

\noindent
\textbf{Convention.} We shall not notationally distinguish between a reduced simplicial set $X$ and its geometric realization $|X|$. Nor shall we distinguish notationally between a topological space and a (reduced) simplicial model of that space.

\subsection{The Kan loop group construction}
\la{klg}
We will briefly review the classical construction of Kan \cite{Kan1} which provides a functorial simplicial group model of the based loop space $ \Omega X $. For details and proofs
we refer the reader to \cite[Chapter~VI]{M} and \cite[Chapter~V]{GJ}).
Let $\sGr$ denote the category of simplicial groups. It has a standard model structure, where the weak equivalences and fibrations of simplicial  groups are the weak equivalences and fibrations of the underlying simplicial sets. We note that, unlike $\sset$, the model category $\sGr$ is fibrant: by a classical theorem of Moore,  every simplicial group is a Kan complex (see \cite[Theorem~17.1]{M}).

\vspace*{1ex}

\begin{definition}
\la{Defsf}
A simplicial group $\Gamma = \{\Gamma_n\}_{n \ge 0} $ is called  \emph{semi-free} if there is a sequence of subsets $\, B_n \subset \Gamma_n \,$, one in each degree, such that $ \Gamma_n $ is freely
generated by $ B_n $, and the set $\, B =  \bigcup_{n \ge 0}\,B_n $ is closed under degeneracies of $ \Gamma $, {\it i.e.}, $\, s_j(B_{n-1}) \subseteq B_n $ for all $ 0 \leq j \leq n-1 $ and $ n \ge 1 $. The subset
$ \overline{B}_n := B_n  \backslash \bigcup_{i=0}^{n-1}s_i(B_{n-1})$
is called the set of {\it nondegenerate generators} of $\Gamma $ of degree $n$.
\end{definition}

\vspace*{1ex}

One can show that every element in $\overline{B}_n$ is nondegenerate (when considered as an element
of the underlying simplicial set), and a semi-free simplicial group is determined by
specifying the sets of nondegenerate generators $\overline{B}_n$ and the face elements of these generators.

Semi-free simplicial groups are cofibrant objects in the model category $\sGr$. The Kan loop group construction provides an important class of semi-free simplicial groups that arise naturally from reduced simplicial sets.
To be precise, the Kan construction defines a pair of adjoint functors:
\begin{equation}
\la{kanl}
\lgr :\, \sset_0\, \rightleftarrows \,\sGr\,:\overline{W}
\end{equation}
where $ \lgr $ is called the {\it Kan loop group functor} and $\overline{W} $ is the {\it classifying simplicial complex}. The functor $ \lgr $ preserves weak equivalences and cofibrations, while $ \overline{W} $ preserves weak equivalences and fibrations (see \cite[Proposition~V.6.3]{GJ}).
Hence, \eqref{kanl} is a Quillen pair, which is actually a Quillen equivalence: i.e., the functors $ \lgr $ and
$ \overline{W} $ induce mutually inverse equivalences between the homotopy categories of
$ \sset_0 $ and $ \sGr $ (see \cite[Corollary~V.6.4]{GJ}).
Combining this with the classical Quillen equivalence \eqref{qeq} between topological spaces and simplicial sets:
$$
 \mathtt{Top}_{0,\ast} \xrightarrow{\ \overline{S}\ } \sset_0  \xrightarrow{\ \lgr\ } \sGr
$$
we get equivalences of the homotopy categories:
$$
\Ho(\mathtt{Top}_{0,\ast})\cong  \Ho(\sset_0) \cong \Ho(\sGr)\ .
$$

For further use, we recall the explicit construction of the functor $\lgr{}$. Given a reduced simplicial set $ X = \{X_n\}_{n \ge 0} $, the set of $n$-simplices of $ \lgr{X} $ is defined by
$$\lgr{X}_n\,=\,\, \frgr{X_{n+1}}/\langle s_0(x)=1\, ,\ \forall \,x\,\in\, X_n \rangle\,\cong \,\frgr{B_n}\ ,
$$
where $ B_n := X_{n+1}\!\setminus s_0(X_n) $ and the isomorphism is induced by the inclusion $ B_n \into X_{n+1} $.
The degeneracy maps $ s_j^{\lgr{X}}: \lgr{X}_{n} \to \lgr{X}_{n+1} $
are induced by the degeneracy maps $ s_{j+1}: X_{n+1} \to X_{n+2} $ of the simplicial set $X$, and the face maps $ d_i^{\lgr{X}}: \lgr{X}_{n} \to \lgr{X}_{n-1} $ are given by
$$
d_0^{\lgr{X}}(x) := (d_1x) \cdot (d_0x)^{-1}\quad \mbox{and}\quad d_i^{\lgr{X}}(x) := d_{i+1}(x)\,,\ \forall\, i>0\ .
$$
Conversely, given a simplicial group $\Gamma = \{\Gamma_n\}_{n \ge 0}$, the simplicial set $\overline{W}\Gamma$ is defined by $\overline{W}\Gamma_0\,:=\,\{\ast\}$  and $\, \overline{W}\Gamma_n\,:=\,\Gamma_{n-1} \times \Gamma_{n-2} \times \ldots \times \Gamma_0$ for $n \geq 0$. The degeneracy and face maps of $\overline{W}\Gamma$  are given explicitly in~\cite[\S 21]{M}.
We note that when restricted to {\it discrete} simplicial groups, the functor $ \overline{W} $ coincides with the usual nerve construction, {\it i.e.}, $ \overline{W}\Gamma =  {\rm B} \Gamma $ for any discrete group $\Gamma$.
\bprop
\la{Kanprop}
The Kan loop group $\lgr{X} $ of any reduced simplicial set $X$ is semi-free.
More precisely, for each $n>0$, the composite map $\tau :  X_n \to \frgr{X_n}  \onto \lgr{X}_{n-1}$ is injective when restricted to the subset
$\overline{X}_n \subset X_n$, and the image $\tau( \overline{X}_n ) \subset \lgr{X}_{n-1}$
forms the set of nondegenerate generators
$\overline{B}_{n-1} = \tau( \overline{X}_n )$ in degree $(n-1) $ of the semi-free basis $\{B_n\}_{n\geq 0}$ of $\,\lgr{X}$.
\eprop
The following fundamental theorem clarifies the meaning of the Kan loop group construction.

\bthm[Kan \cite{Kan1}]
\la{kanstheorem}
For any reduced simplicial set $X$, there is a weak homotopy equivalence
%
$$|\lgr{X}| \simeq \Omega|X|\,,$$
where $\Omega|X|$ is the $($Moore$)$ based loop space of $|X|$.
\ethm
A detailed proof of Theorem~\ref{kanstheorem} can be found in \cite[\S\,26]{M}.
Its significance becomes clear from the following considerations.
Given any path-connected CW complex $ Y $ one can choose a pointed connected simplicial set $X'$ such that
$ |X'| \simeq Y $. If $ X $ is the path-connected component of $ X' $ containing the basepoint, then
$X $ is a reduced simplicial set such that $\, |X| \simeq |X'| \simeq Y \,$ because $Y$ is connected.
Hence, applying the Kan loop group construction to $X$, we get $ |\lgr{X}| \simeq \Omega Y $.
Thus, $\lgr{X}$ is a semi-free simplicial group model of the based loop space of $Y$. In this way,
the based loop space of any path-connected CW complex admits a simplicial group model.

\section{Representation homology}
\la{S3}

In this section, we define representation homology as the homotopy groups of the (non-abelian) derived representation functor associated with an affine algebraic group. We establish the existence and basic properties of this functor as well as indicate some generalizations.
Our construction follows the approach of our earlier papers \cite{BKR, BFPRW, BR} where
we studied the  representation homology of associative and Lie algebras.

\subsection{Definition of representation homology}
\la{rephomsect}
Fix an affine algebraic group scheme $G$ over a field $k$ of characteristic $0$.  Recall that
$G$ is given by a representable functor on the category of commutative $k$-algebras with values in the category of groups:
\begin{equation}
\la{grsc}
G:\,\cAlg_k \rar \mathtt{Gr}\,,\,\,\,\, A \mapsto G(A)\ .
\end{equation}
A commutative algebra that represents  \eqref{grsc} is called
the coordinate ring of $ G $ and denoted $ \O(G) $. This algebra is equipped with a coproduct $ \Delta: \O(G) \to \O(G) \otimes \O(G) $, $\,f \mapsto f^{(1)} \otimes f^{(2)}\,$, which is dual to the multiplication in $ G $ and makes $ \O(G) $ a commutative Hopf algebra.

\vspace*{1ex}

\noindent
{\bf Lemma/Definition.}\
The functor \eqref{grsc} has a left adjoint
\begin{equation}
\la{lgrsc}
(\,\mbox{--}\,)_G\,:\,\mathtt{Gr} \to \cAlg_k\ ,\quad \Gamma \mapsto \Gamma_G\ .
\end{equation}
which we call the {\it representation functor} in $ G $.

\bproof
Given a group $ \Gamma \in \Gr $ define the  algebra $ \Gamma_G $ by
the following canonical presentation
$$
\Gamma_G = {\rm Sym}_k (\, k[\Gamma] \otimes_k \O(G)\,)/I ,
$$
where the ideal $I$ of relations is generated by
\begin{eqnarray}\la{relAA}
&& \gamma \otimes f_1 f_2 - (\gamma \otimes f_1) \cdot (\gamma \otimes f_2)\ ,  \nonumber\\
&&\gamma_1 \gamma_2 \otimes f - (\gamma_1 \otimes f^{(1)})\cdot (\gamma_2 \otimes f^{(2)})\ ,
\\
&& e_\Gamma \otimes f - f(e_G)\cdot 1\ ,\quad \gamma \otimes 1 - 1  \nonumber
\end{eqnarray}
for all $\, \gamma, \,\gamma_1,\, \gamma_2 \in \Gamma \,$ and $\,f, \,f_1, f_2 \in \O(G)\,$.
If $A \in \cAlg_k $ is a commutative algebra, a group homomorphism $ \varphi: \Gamma \to G(A) =
\Hom(\O(G),A) $ determines a linear map $ k[\Gamma] \otimes \O(G) \to A $, which, in turn, induces --- modulo the relations \eqref{relAA} --- an algebra homomorphism $ \varphi^{\#}: \Gamma_G \to A $. It is straightforward to check
that $ \varphi \mapsto \varphi^{\#} $ gives the required bijection
$\,\Hom_{\Gr}(\Gamma, G(A)) \cong \Hom_{\cAlg_k}(\Gamma_G, A)\,$.
\eproof
We remark that, for a fixed group $ \Gamma $, the algebra $ \Gamma_G $ represents the functor
$$
\Rep_G(\Gamma): \cAlg_k \to \set \ ,\quad A \mapsto \Hom_{\Gr}(\Gamma, G(A))\ ,
$$
which is the functor of points of an affine $k$-scheme $ \Rep_G(\Gamma) $
parametrizing the representations of $ \Gamma $ in $G$; hence, geometrically,
we can think of $ \Gamma_G $ as the coordinate ring $ \O[\Rep_G (\Gamma)] $ of $ \Rep_G(\Gamma) $.

Next, we embed the category of groups into the category
$\sGr $ of simplicial groups and extend the functor \eqref{lgrsc} to $\sGr $ in the natural way, assigning to a simplicial group $ \Gamma_\ast: \Delta^{\rm op} \to \Gr $ the simplicial commutative algebra
$\,(\Gamma_\ast)_G: \Delta^{\rm op} \to \Gr \to \cAlg_k\,$. We will keep the notation $ (\,\mbox{--}\,)_G $ for this extended representation functor:
\begin{equation}
\la{srep}
(\,\mbox{--}\,)_{G}:\, \sGr \to \mathtt{s}\cAlg_k\ .
\end{equation}
Both categories $ \sGr $ and $ \scAlg_k $ have natural (simplicial) model structures, with weak
equivalence being the weak homotopy equivalence of the underlying simplicial sets.  However,
the representation functor \eqref{srep} is not homotopy invariant --- it does not preserve weak equivalences --- hence, in order to work
in a homotopical context we should replace or approximate \eqref{srep} with a derived functor (see \cite{Q1}, \cite{DS95}). The existence of this derived functor is easy to establish.
\blemma
\la{sgfun}
The functor \eqref{srep} maps the weak equivalences between cofibrant objects in $ \sGr $
to weak equivalences in $ \scAlg_k $, and hence has a total left derived functor
\begin{equation}
\la{Lrep}
\L(\,\mbox{--}\,)_{G} :\,\Ho(\mathtt{sGr}) \rar \Ho(\mathtt{s}\cAlg_k)\,\text{.}
\end{equation}
\elemma
\begin{proof}
Suppose that $f: \Gamma \rightarrow \Gamma' $ is a weak equivalence between cofibrant simplicial groups.
Since $\mathsf{sGr}$ is a fibrant model category, $ \Gamma $ and $ \Gamma' $ are both fibrant-cofibrant objects.
By Whitehead's Theorem, the map $f$ has then a homotopy inverse $ g: \Gamma' \rightarrow \Gamma $, such that
$\,fg \sim \id\,$ and $\,gf \sim \id\,$. Now, any homotopy between fibrant-cofibrant objects can be realized
using a good cylinder object in $ \sGr $. Since $\mathsf{sGr}$ is a simplicial model category,
there is a natural choice of good cylinder objects for $\Gamma $ and $ \Gamma' $: namely, $\, \Gamma \sqcup \Gamma \rightarrow \Gamma \times \Delta[1] \rightarrow \Gamma\,$, and similarly for $ \Gamma' \,$.
For such cylinder objects, the simplicial homotopies (see \cite[Def.~5.1]{M}) can be defined by explicit combinatorial relations which are preserved by the functor $ (\,\mbox{--}\,)_{G} $. Thus, we conclude that $\, g_G: \Gamma_G' \to \Gamma_G $ is a homotopy inverse of $\, f_G: \Gamma_G \rightarrow \Gamma_G' \,$ in $ \mathtt{s}\cAlg_k $ and hence $ f_G $ and $ g_G $ are mutually inverse isomorphisms in $ \Ho(\mathtt{s}\cAlg_k) $. The existence of the derived functor \eqref{Lrep} follows now from \cite[Prop.~9.3]{DS95}.
\end{proof}
Now, for a fixed simplicial group $\Gamma \in \sGr $, we formally define the {\it derived representation scheme} $\, \DRep_{G}(\Gamma)\, $ as $\, \Spec\,\L(\Gamma)_{G}\, $, i.e. the simplicial
algebra $ \L(\Gamma)_{G} $ viewed as an object of the opposite category $\Ho(\scAlg_k)^{\mathrm{op}}$. We call the  homotopy groups of $ \L(\Gamma)_{G} $ the {\it representation homology of $\Gamma$ in $G\,$} and write
$$
\HR_\ast(\Gamma,G)  := \pi_\ast \L(\Gamma)_{G} \ .
$$

By comparing the universal mapping properties, it is easy to check that the  functor \eqref{srep} commutes with $ \pi_0 $; hence, for any $\Gamma \in \sGr $, there is a natural isomorphism in $ \cAlg_k $:
\begin{equation}
\la{commpi}
\HR_0(\Gamma,G)\,\cong\, [\pi_0(\Gamma)]_{G}
\end{equation}
In particular, if $\Gamma \in \Gr $ is a constant simplicial group, we have
$\,\HR_0(\Gamma,G) \cong \Gamma_{G} \,$, which
justifies our notation and terminology for $ \DRep_{G}(\Gamma) $.

Next, recall the fundamental theorem of Kan \cite{Kan1} that identifies the homotopy types of simplicial groups with those of pointed connected spaces. To be precise, the Kan Theorem asserts that the category of simplicial groups is Quillen equivalent to the category $ \sset_0 $ of reduced simplicial sets, which is, in turn, Quillen equivalent to the category $ \mathtt{Top}_{0,*} $ of pointed connected  spaces. As a result, we have natural equivalences of homotopy categories
\begin{equation}\la{tsg}
\Ho(\mathtt{Top}_{0,\ast}) \,\cong\, \Ho(\sset_0)\,\cong\, \Ho(\sGr)\ .
\end{equation}
This leads us to the main definition.
\begin{definition}
\la{DRepX}
For a space $X \in \mathtt{Top}_{0,*} $, we define the {\it derived representation scheme}
$\,\DRep_G(X)\,$ to be $\,\DRep_G({\Gamma}X)\,$, where $ {\Gamma}X$  is a(ny) simplicial group model\footnote{that is, a simplicial group that corresponds to $X$ under the Kan equivalence \eqref{tsg}.} of $X$. The {\it representation homology of $X$ in $G$} is then defined by
\begin{equation}
\la{hrX}
\HR_\ast(X,G)   := \pi_\ast \L({\Gamma}X)_{G} \ .
\end{equation}
\end{definition}
By definition,  $ \HR_*(X,G) $ is a graded commutative algebra, with
$ \HR_0(X,G) $  naturally isomorphic to $ [\pi_1(X)]_G = \cO[\Rep_G(\pi_1(X))]$,
the coordinate ring of the representation scheme $ \Rep_G[\pi_1(X)] $.
The last isomorphism is the composition of \eqref{commpi} with the natural isomorphism $ \pi_0(\Gamma X) \cong \pi_1(X) $.

For a reduced simplicial set $ X \in \sset_0 $, the Kan loop group $ \lgr{X} $
provides a canonical (functorial) simplicial group model for $ |X| $. Since this
simplicial group is semi-free (see Section~\ref{klg}), we have
\begin{equation}
\la{explform}
\HR_\ast(X,G) \cong \pi_\ast \,(\lgr{X})_{G}
\end{equation}
This formula can be used to compute representation homology in some simple cases.
\begin{example}
\la{addgr}
Let $ \lgr_a $ be the additive group over $k$, i.e. the affine algebraic group defined by
the functor $ \lgr_a: \cAlg_k \to \Gr \,$,$\,A \mapsto (A, +)\,$, where $ (A, +) $ denotes
the underlying abelian group of the algebra $A$. It is easy to see that,
for any $ \Gamma \in \Gr $, there is a natural bijection
$\, \Hom_{\Gr}(\Gamma, \lgr_a(A)) \cong \Hom_{\cAlg_k}({\rm Sym}_k (\lin_k \Gamma), A) \,$,
where $\, \lin_k(\Gamma) := \Gamma_{\rm ab} \otimes_{\Z} k \,$.
Hence, the representation functor in $ \lgr_a $ is given by the composition
$\,(\,\mbox{--}\,)_{\lgr_a} = {\rm Sym}_k \circ \lin_k:\,\Gr \to \ve_k \to \cAlg_k \,$.
Using formula \eqref{explform}, for an arbitrary $ X \in \sset_0 $, we can now compute
\begin{eqnarray*}
\HR_\ast(X,\lgr_a)
&\cong &  \pi_\ast\,\Sym_k\,[(\lgr{X})_{\rm ab} \otimes_\Z k] \\
& \cong & \Sym_k\left[\pi_\ast((\lgr{X})_{\rm ab} \otimes_\Z k)\right]\\
& \cong & \Sym_k\left[\pi_\ast(\lgr{X})_{\rm ab} \otimes_\Z k\right]\\
& \cong & \Sym_k\left[\H_{\ast+1}(X, \Z) \otimes_\Z k\right]\\
& \cong & \Sym_k\left[\H_{\ast+1}(X, k)\right]\ ,
\end{eqnarray*}
where $ \Sym_k $ is the graded symmetric algebra functor over $ k $ and
$ \H_{\ast+1}(X, \Z) := \oplus_{i \ge 0}\, \H_{i+1}(X, \Z) $ is the singular
homology of $X$. Note that, besides \eqref{explform}, we used here the classical
isomorphism of Kan: $\, \pi_\ast(\lgr{X})_{\rm ab} \cong \H_{\ast+1}(X, \Z)\,$
(see, e.g., \cite[Theorem~26.9]{M}) and the well-known fact that the functor
$ \Sym_k $ commutes with homology when $k$ has characteristic zero
(see, e.g., \cite[Part~I, Prop. 4.5]{Q}). This example shows that
representation homology may be viewed as a generalization of the ordinary singular
homology.
\end{example}

\subsection{The derived representation adjunction}
\la{DRA}
By definition, the representation functor \eqref{lgrsc} is left adjoint to
(the functor of points of) the algebraic group $G$. This adjunction extends automatically
to simplicial categories:
\begin{equation}
\la{sadj}
(\,\mbox{--}\,)_G\,: \,\sGr\, \rightleftarrows \, \scAlg_k\, :\, G\ ,
\end{equation}
and the natural question is whether \eqref{sadj} induces an adjunction between derived functors
on the corresponding homotopy categories. The (affirmative) answer to this question would be
immediate from Quillen's fundamental theorem \cite[I.4.5, Theorem~3]{Q1} if \eqref{sadj} were a pair of
Quillen functors between model categories. However, this is not the case. By definition, any left Quillen
functor preserves cofibrations, which means, in particular, that it maps cofibrant objects in
one model category to cofibrant ones in the other. Unfortunately, the representation functor
\eqref{srep} lacks this property even in simplest cases. Take, for example,
$G=\mathbb{G}_m$, the multiplicative group, and apply \eqref{srep} to the free group on
one generator $\Gamma= \mathbb{F}_1$, which is obviously a cofibrant object in $ \sGr$.
The result is $\Gamma_G\,\cong\,k[x,x^{-1}]$, which is not a cofibrant simplicial algebra in $\mathtt{s}\cAlg_k$. Another problem is that the right adjoint functor in \eqref{sadj} is not homotopical and
hence should be replaced by a right derived functor $ \bR G $. But the existence of $ \bR G $
is not clear because the standard (projective) model structure on $ \scAlg_k $ is fibrant.
Nevertheless, somewhat surprisingly, we still have
\bthm
\la{ThAdj}
The algebraic group functor $ G:\, \scAlg_k \to \sGr $ has a total right derived functor,
 which is right adjoint to the derived representation functor \eqref{Lrep}:
\begin{equation*}
\L(\,\mbox{--}\,)_G\,: \,\Ho(\sGr)\, \rightleftarrows \, \Ho(\scAlg_k)\, :\, \bR G
\end{equation*}
Moreover, both $\L(\,\mbox{--}\,)_G$ and $ \bR G$ are absolute derived functors in the sense of Deligne-Maltsiniotis \cite{Mal}.
\ethm
The first statement of Theorem~\ref{ThAdj} shows that the simplicial adjunction \eqref{sadj}
behaves like a Quillen adjunction, and the last statement shows that the corresponding derived functors
are as `good' (well-behaved) as derived functors of Quillen functors. In particular, like the
derived functor of a left Quillen functor, the derived representation functor has the following
important property which plays a crucial role in computations of representation homology in Section~\ref{sec7}.
\begin{theorem}
\la{hopushout}
The derived representation functor \eqref{Lrep} preserves arbitrary $($small$)$ homotopy colimits.
\end{theorem}

The main idea behind our proof of Theorem~\ref{ThAdj} is to `forget' the model structure on $ \sGr $, thinking of this category simply as a homotopical category
in the sense of \cite{DHKS}, and then `approximate' it with another model category, which is
`almost' Quillen equivalent to $ \sGr $ in the sense of \cite{CS}.
For reader's convenience, we recall basic definitions and the necessary results from \cite{DHKS}
and \cite{CS} in the Appendix, where we also prove abstract versions of Theorem~\ref{ThAdj}
and Theorem~\ref{hopushout} (see Theorem~\ref{ThA2} and Theorem~\ref{ThA3}, respectively). The proof of Theorem~\ref{ThAdj}
will consist of verifying the conditions of Theorem~\ref{ThA2}; we will divide it
into two cases: $\, G = \GL_n$ and the general case: $G$ is an arbitrary
algebraic group. For $ \GL_n $, we will provide detailed arguments, while in the general case, we will only sketch the proof leaving technical details for our subsequent paper.

We begin with the following observation refining the result of Lemma~\ref{sgfun}.
\bprop
\la{sgfun2}
The representation functor \eqref{srep} is left deformable on $ \sGr $,
and hence its total left derived functor \eqref{Lrep} is an absolute derived functor
in the sense of \cite{Mal}.
\eprop
\bproof
Write $ \C $ for $ \sGr $ viewed as a homotopical category, and let $ \C_Q $
denote its full subcategory consisting of semi-free simplicial groups (see Definition~\ref{Defsf}).
Then $ \C_Q \into \C $ is a left deformation retract of $ \C $, with retraction functor
$ Q: \C \to \C $ being the composition $ Q := \mathbb{G} \,\overline{W} $ and
the morphism $ q: Q \to \id_{\C} $ given by the counit of the Kan loop group adjunction
\eqref{kanl}. Indeed, by Kan's theorem, the morphism $q$ is a natural weak equivalence  (see \cite[Prop. V.6.3.]{GJ}), and its image is contained in $ \C_Q $. The proof of Lemma~\ref{sgfun} shows that
\eqref{srep} is homotopical on $ \C_Q $ and hence, by definition, left deformable.
The result now follows from Proposition~\ref{leftdef}.
\eproof
\subsection{Proof of Theorem~\ref{ThAdj} for $ G = \GL_n $}
\la{S33}
Let $ \sMon $ denote the category of simplicial monoids equipped
with the standard (projective) model structure. Consider the natural adjunction
\begin{equation}
\la{sMonGr}
l: \sMon \, \rightleftarrows \, \sGr: r
\end{equation}
where $r$ is the inclusion functor, and $ l $ is the group completion (localization) functor.
The next observation is a consequence of a known theorem of Dwyer and Kan \cite{DK}.
\blemma
\la{DKlemma}
The adjunction \eqref{sMonGr} is a left model approximation of $ \sGr $ in the sense of \cite{CS}.
\elemma
\bproof
We need to verify the three conditions of Definition~\ref{DefA1} (see Appendix). Condition $ (1) $ is obvious,
$ (2) $ follows from the fact that $ l $ is a left Quillen functor (when $ \sMon $ and $ \sGr $
are regarded as model categories). It suffices only to check $(3)$. For this, observe that
any map from a simplicial monoid to a simplicial group, say $ f: M \to r (\Gamma) \,$, can be factored
 as $\,M \xrightarrow{\eta_M} r l(M) \xrightarrow{r f^{\#}} r (\Gamma)\,$, where
$ \eta_M $ is the group completion (localization) map and $ f^{\#}: l(M) \to \Gamma $ is the map adjoint to $f$ under \eqref{sMonGr}. If $ f $ is a weak equivalence in $ \sMon $,
then $ M $ is a group-like simplicial monoid (i.e., $\pi_0(M) \cong \pi_0(\Gamma) $ is a group), and
hence, if $ M $ is also cofibrant, by \cite[Proposition~10.4]{DK}, $ \eta_M $ is a weak equivalence. By 2-of-3 property, the map $ r f^{\#}: r l(M) \to r (\Gamma) $
is then a weak equivalence, and since $r$ reflects weak equivalences,
$\, f^{\#}: l(M) \to \Gamma $ is a weak equivalence as well.
\eproof
We will apply Theorem~\ref{ThA2} to the representation adjunction  \eqref{sadj} with $ G = \GL_n $ using  the left model approximation \eqref{sMonGr}. Note that this model approximation is good for the representation functor \eqref{srep}  (for any algebraic group $G$), since the group completion functor maps cofibrant objects in $ \sMon $, which
are (retracts of) semi-free simplicial monoids, to (retracts of) semi-free simplical groups, on which the functor \eqref{srep} is homotopical by Lemma~\ref{sgfun}.

From now on, we assume that $ G = \GL_n $ and write $ F :=  (\,\mbox{---}\,)_{\GL_n}\,: \,\sGr \to \scAlg_k $ for the corresponding representation functor.
Let  $ \sMon_0 $ denote the full subcategory of $ \sMon $ consisting of group-like simplicial monoids ---  by the Dwyer-Kan theorem \cite{DK}, the essential image of the functor $ \bar{r}: \Ho(\sGr) \into \Ho(\sMon) $
is precisely $ \Ho(\sMon_0) $.

The inclusion functor $i: \sMon_0 \into \sMon $ has
a right adjoint $ \hat{U}: \sMon \to \sMon_0 $ defined by the pull-back diagram in $ \sMon \,$:
\begin{equation*}
\begin{diagram}[small]
\hat{U}(M_\ast) & \rTo^{} & U[\pi_0 (M_\ast)] \\
\dTo^{}         &         & \dTo^{}\\
M_\ast &  \rTo^{}& \pi_0(M_\ast)
 \end{diagram}
\end{equation*}
where $ U: \Mon \to \Gr $ is the functor assigning to a monoid its subgroup of units.
To construct the functors $ \hat{F} $ and $ \hat{G} $ we start with the natural adjunction
$\,
(k[\,\mbox{--}\,])_n: \,\sMon\, \rightleftarrows \, \scAlg_k\,: \, M_{n}
\,$
and compose it with $\, i: \sMon_0 \rightleftarrows \sMon: \hat{U} \,$, i.e. define
\begin{equation}
\la{hatadj}
\hat{F} := (k[\,\mbox{--}\,])_n \circ i: \ \sMon_0\, \rightleftarrows \, \scAlg_k\,: \ \hat{G} :=
\hat{U} \circ M_{n}
\end{equation}

Note that the right adjoint $ \hat{G} $ in \eqref{hatadj} is precisely the Waldhausen functor $ \widehat{\GL}_n $ defined in the Introduction (see \eqref{waldfun}). In particular,  it is a homotopical functor  that takes its
values in group-like simplicial monoids: thus, we have $ \bR \hat{G} = \hat{G}$ and $ \im(\hat{G}) \subseteq \im(\bar{r})$, i.e. condition $(iii)$
of Theorem~\ref{ThA2} holds.

The left adjoint $ \hat{F} $ is obtained by restricting to $ \sMon_0 $ the functor $ (k[\,\mbox{--}\,])_n $, which is left Quillen on $\sMon\,$: hence, $ \hat{F} $ is homotopical on cofibrant objects in $ \sMon_0 $, and therefore Theorem~\ref{ThA2}$(i)$ holds.

Next, factor $ r:= r_0 \circ i:  \sGr \into \sMon_0 \into \sMon \,$
and observe that $ \hat{F} \circ r_0  =
(k[\,\mbox{--}\,])_n \circ r $ is left adjoint to
$ U \circ M_n = \GL_n $. Hence, there is a canonical isomorphism of functors
$\, \hat{F} \circ r_0 \cong F \,$. It remains only to show that
$\,\sGr \xrightarrow{r_0} \sMon_0 \xrightarrow{\hat{F}} \scAlg_k \,$
is a left deformable pair. For this, in the notation of Proposition~\ref{comdef}, we take $\, \C_Q $ to be the full subcategory
of semi-free simplicial groups in $ \C := \sGr $ and $ \D_Q $  the full subcategory
of $ \D := \sMon_0 $ consisting of monoids $ M $ such that $ k[M] $ is a simplicial $k$-algebra that is degreewise a direct limit of formally smooth $k$-algebras having semifree DG resolutions with finitely many generators in each homological degree. Both $ \C_Q $ and $ \D_Q $ are left deformation
retracts of the corresponding homotopical categories: $ \C_Q $ contains the image of the deformation functor $ Q =  {\mathbb G} \, \overline{W} $ associated with the Kan
loop group adjunction (see Proposition~\ref{sgfun2}), while $ \D_Q $ contains the image of  $ Q_0: \sMon_0 \to \sMon_0 $,  which is the restriction of  the cofibrant reprelacement functor $ Q $ on $ \sMon $. Since $ r_0 $ is homotopical and $ r_0(\C_Q) \subseteq \D_Q $, we need only to check that $ \hat{F} $ is homotopical on $ \D_Q $.  Since  $ \hat{F} = (\,\mbox{--}\,)_n \circ k[\,\mbox{--}\,]  $ and $k[\,\mbox{--}\,]$ is homotopical on $ \sMon_0 $, it suffices to check that $ \, (\,\mbox{--}\,)_n:\, \sAlg_k  \to \scAlg_k $ is homotopical on simplicial $k$-algebras that are degreewise (direct limits of) formally smooth $k$-algebras having semifree DG resolutions with finitely many generators in each homological degree. Now, this last fact follows from \cite[Theorem~21]{BFR}, saying that such associative algebras are adapted for the representation functor
$\, (\,\mbox{--}\,)_n: \Alg_k \to \cAlg_k \,$ (in the sense that
$\L(A)_n \cong A_n $ for such $A$'s) and the well-known abstract result
from homotopical algebra saying that the simplicial objects, which are  {\it degreewise} adapted for a functor $ F $, are actually adapted for $F$
(see, e.g., \cite[Theorem 9.2.2]{WY}).

Summing up, we showed that all three conditions of Theorem~\ref{ThA2} hold for the adjoint pair \eqref{hatadj}, except that the left adjoint $ \hat{F} $ is not defined on the entire model category
$ \sMon $ but rather on its full subcategory $ \sMon_0 $. However, this last subcategory is closed under weak equivalences and coincides with  $ \im(\bar{r}) $, hence
the result of Theorem~\ref{ThA2} still holds (see Remark~2 after the proof of Theorem~\ref{ThA2}).
This completes the proof of Theorem~\ref{ThAdj} for $ G = \GL_n $.

Theorem~\ref{ThA2} gives an explicit formula for the total derived functor $ \bR G $:
namely,
\begin{equation}
\la{func}
\bR \GL_n =  \L l \circ \widehat{\GL}_n \ .
\end{equation}
This allows us to compute the homotopy groups $\, \bR^i \GL_n := \pi_i \bR \GL_n \,$ for all $\, i \ge 0 $.
\bprop
For any $ A_\ast \in \scAlg_k $,
\begin{equation*}
\bR^i \GL_n(A_*) \,\cong \,\left\{\!\!
\begin{array}{lll}
\GL_n[\pi_0(A_*)] & \mbox{for} & i = 0 \\*[1ex]
{\rm M}_n[\pi_i(A_*)] & \mbox{for} & i \ge 1
\end{array}
\right.
\end{equation*}
\eprop
\bproof
Let $ Q \widehat{\GL}_n(A_\ast) \xrightarrow{\sim} \widehat{\GL}_n(A_\ast) $ be a cofibrant resolution of
$ \widehat{\GL}_n(A_\ast) $ in the (model) category $ \sMon $. By \eqref{func}, we have
$\,\bR \GL_n(A_\ast) \cong l Q \widehat{\GL}_n(A_\ast)\,$. On the other hand,
$ Q \widehat{\GL}_n(A_\ast) $ is a group-like simplicial monoid, since so is $ \widehat{\GL}_n(A_\ast) $. Hence,
by the Dwyer-Kan Theorem, the group completion map
$\, Q \widehat{\GL}_n(A_\ast) \xrightarrow{\sim} l Q \widehat{\GL}_n(A_\ast) \,$ is a weak equivalence. Thus, we have a zig-zag of weak equivalences
$$
\widehat{\GL}_n(A_\ast) \xleftarrow{\sim} Q \widehat{\GL}_n(A_\ast)
\xrightarrow{\sim} l Q \widehat{\GL}_n(A_\ast)\ ,
$$
from which we conclude that $ \bR^i \GL_n(A_\ast) \cong \pi_i \widehat{\GL}_n(A_\ast) $.
The result now is immediate from the definition of $ \widehat{\GL}_n $ (see \cite{Wa}).
\eproof
\subsection{(Sketch of) Proof of  Theorem~\ref{ThAdj} in the general case}
\la{S34}
For a general algebraic group $G$, we will use a different model approximation of $ \sGr $
given by reduced simplicial spaces. By a {\it simplicial space} we mean a bisimplicial
set of which we think as a functor $ X_*: \Delta^{\rm op} \to \sset $, $\,[n] \mapsto X_n\,$,
with simplicial components $ X_n $ viewed as `vertical' simplicial sets. We call $ X_* $ {\it reduced} if $X_0 = \Delta[0] $ is the one-point (discrete) simplicial set. We write
$ \sSp = \sset^{\Delta^{\rm op}} $ for the category of all simplicial spaces and $ \sSp_* $ for
its full subcategory consisting of reduced ones. The category $ \sSp $ is known to carry several
interesting model structures. We will use two of these: the projective model structure\footnote{This is a special case of the Bousfield-Kan model structure on the category $ \sset^I $ of diagrams
of simplicial sets for $ I = \Delta^{\rm op}$ (see, e.g., \cite[Sect. IV.3.1]{GJ}).}  in which the weak equivalences and fibrations are the levelwise weak equivalences (resp., fibrations) of simplicial sets and its (left Bousfield) localization  with respect to Segal maps introduced in \cite{Re}. We denote
the projective model structure simply by $ \sSp$ and its localization by $ \LsSp $. As shown in \cite{Be07}, both model structures `restrict' to reduced simplicial spaces,
and we denote the corresponding model categories by $ \sSp_* $ and $ \LsSp_* $, respectively.

The reduced simplicial spaces are related to simplicial groups by the pair of adjoint functors
\begin{equation}
\la{sSpGr}
\underline{\pi}_1 : \LsSp_* \, \rightleftarrows \, \sGr: \underline{N}
\end{equation}
where $ \underline{N} $ is the nerve functor applied degreewise to
components of simplicial groups: i.e., for $ \Gamma_* \in \Ob(\sGr) $,
$$
\underline{N}_{\bullet}(\Gamma_*):\,\Delta^{\rm op} \to \sset\ , \quad
[n] \mapsto N_n(\Gamma_*) = \Gamma_*^n \ ,
$$
and $ \underline{\pi}_1 $ is the fundamental group functor applied degreewise to bisimplicial sets: i.e., for $ X = \{X_{p,q}\}_{p,q \ge 0} $,
$$
\underline{\pi}_1(X): \Delta^{\rm op} \to \Gr\ ,\quad [q] \mapsto \pi_1(X_{\ast, q})\ .
$$
The fact that $ \underline{\pi}_1 $ is  left adjoint to $ \underline{N} $
follows from the well-known fact that the fundamental group functor
$ \pi_1: \sset_0 \to \Gr $ on reduced simplicial sets is left adjoint to
the simplicial nerve $ N: \Gr \to \sset_0 $ on the category of groups.
Now, in place of Lemma~\ref{DKlemma}, we have
\blemma
\la{DKlemma2}
The adjunction \eqref{sSpGr} is a left model approximation of $ \sGr $.
\elemma
\bproof
This follows from a theorem of Bergner (see \cite[Theorem~1.6]{Be07}) which
asserts that the model category $ \LsSp_* $ is Quillen equivalent to the category
of simplicial monoids, $ \sMon $, equipped with the standard (projective) model structure. In fact, one can check that
\eqref{sSpGr} factors as $\,\LsSp_*  \,\rightleftarrows\, \sMon \,\rightleftarrows  \,\sGr\,$, where the first adjunction is Bergner's Quillen equivalence, with its right adjoint being
a homotopical functor, and the second adjunction is \eqref{sMonGr}, which is, by Lemma~\ref{DKlemma}, a left model approximation of $ \sGr $. It follows that  \eqref{sSpGr} is a left model approximation of $ \sGr $ as well.
\eproof

Next, to define the functor $ \hat{G}: \scAlg_k \to \LsSp_* $ we will use a construction of Galatius
and Venkatesh (see \cite[Section~5]{GV}). We start with the cosimplicial commutative algebra
$$
\O(N_\bullet G):\,\Delta \to \cAlg_k\ ,\quad [n] \mapsto \O(N_n G) = \O(G)^{\otimes n}\ .
$$
Taking the cofibrant replacement  of $ \O(N_\bullet G) $
in $ \scAlg_k$ in each cosimplicial degree, we get a cosimplicial simplicial algebra $ c\,\O(N_\bullet G) \in \scAlg_k^{\Delta} $ and then define $ {\hat G} $ by
\begin{equation}
\la{Ghat}
\hat{G}: \scAlg_k \to \LsSp_\ast\ ,\quad A_* \mapsto {\rm Map}\,(c\,\O(N_\bullet G), \,A_\ast)\ ,
\end{equation}
where `Map' stands for the standard (simplicial) function complex in $ \scAlg_k $. Note that $ {\hat G}(A_\ast)_0 := {\rm Map}\,(c\,\O(N_0 G), \,A_\ast) =
{\rm Map}\,(k, \,A_\ast) \cong \Delta[0] $ for any $ A_\ast $, so
$ \hat{G} $ indeed takes its values in the category of {\it reduced} simplicial spaces.

By formal properties of function complexes, the functor $ \hat{G} $ has a left adjoint given by
\begin{equation}
\la{Fhat}
\hat{F}:\ \LsSp_\ast \to \scAlg_k \ ,\quad X \mapsto X \otimes_{\Delta}
c\,\O(N_\bullet G)\ ,
\end{equation}
where $ \,\otimes_{\Delta}\, $ is the functor tensor product
over the category $ \Delta $ in the simplicial category $ \scAlg_k $
(see, e.g., \cite[(4.1.1)]{R}).
\bprop
\la{adjfg}
The adjoint functors
$\,\hat{F}: \LsSp_\ast  \rightleftarrows  \scAlg_k : \hat{G} \,$
form a Quillen pair.
\eprop
\bproof[Sketch of proof]
One proves this in two steps. First, one checks that the functors $ (\hat{F}, \hat{G}) $
form a Quillen pair $\,\hat{F}: \sSp_\ast  \rightleftarrows  \scAlg_k : \hat{G} \,$ for the projective model structure on $ \sSp_* $. Then, one shows that this Quillen pair `localizes' to a Quillen pair on $ \LsSp_* $
by checking that the left derived functor $ \bL F: \sSp_* \to \Ho(\scAlg_k) $
maps the Segal morphisms in $ \sSp_* $ to isomorphisms in $\Ho(\scAlg_k)$
\eproof
We have now defined all ingredients of Theorem~\ref{ThA2}. To show that
this theorem applies to the representation adjunction \eqref{sadj} we need to verify its assumptions $(i)$, $(ii)$ and $(iii)$.  Condition $(i)$ ---
$ (\hat{F}, \hat{G}) $ being a deformable adjunction ---
is immediate from Proposition~\ref{adjfg}. Condition $(iii)$ is not difficult to check since  $ \hat{G} $ is a homotopical functor (and therefore $ \bR \hat{G} = \hat{G} $). The main work is to verify condition $(ii)$: in particular, to prove that $ \L F \cong \L \hat{F} \circ \underline{N}\,$.
The details of this verification will appear in our subsequent paper. Here we only mention that, as an intermediate step, we
prove the following lemma which provides an alternative way to define
representation homology of spaces, without using the Kan loop group
construction.
\blemma
For any $ X \in \sset_0 $, there is a natural isomorphism in $ \Ho(\scAlg_k)$:
$$
\L F(\lgr X) \cong \L \hat{F}(X^t)\ ,
$$
where $\,(\,\mbox{--}\,)^t: \sset_0 \into \sSp_* \,$ is the `transpose'
inclusion functor identifying a simplicial set $ X  $
with the simplicial space $ X^t = \{(X^t)_n\}_{n\ge 0} $
with discrete components $ (X^t)_n = X_n $.
\elemma

We conclude by pointing out that, once the conditions of Theorem~\ref{ThA2} are verified and Theorem~\ref{ThAdj} is proved, Theorem~\ref{hopushout} follows immediately from
Theorem~\ref{ThA3}, since in our situation both $ \C =  \sGr $ and $ \D = \scAlg_k $ carry model category structures and hence the colimits on these categories exist and are left deformable by results of \cite{CS} (see Theorem~\ref{ThA4}$(3)$).

\section{Functor homology interpretation}
\la{section_rephom}
In this section, we give our second definition of representation homology
parallel to Pirashvili's definition of higher Hochschild homology \cite{P1}.
We begin by reviewing the construction of \cite{P1}.

\subsection{Higher Hochschild homology}
\la{S22}
Let $ \fin_{\ast} $ denote the category of finite pointed sets with objects $ [n] = \{0,1, \ldots, n\} $, $ n\ge 0 $, and morphisms $ f: [n] \to [m] $ being arbitrary set maps such that $ f(0) = 0 $. Let $F :\,\fin_{\ast} \rar \ve_k $ be a covariant functor. We extend $ F $ to the category $\set_{\ast}$ of {\it all} pointed sets in a natural way, using the left Kan extension along the inclusion $\fin_{\ast} \hookrightarrow \set_{\ast}$. We keep the notation $ F $ for the extended functor: explicitly, $ F: \set_\ast \to \ve_k $ is given by
$\, F(X) = \colim\, F([n])\,$, where the colimit is taken over all pointed inclusions $ [n] \into X $.

Given a pointed simplicial set $X \in \sset_* $, we define a simplicial $k$-vector space $F(X)$ as the composition of functors
\begin{equation} \la{defpointed}
F(X):\
\Delta^{\mathrm{op}}  \xrightarrow{X}  \set_{\ast} \xrightarrow{F} \ve_k \ .
\end{equation}
We denote the homotopy groups of $F(X)$ by $ \pi_{\ast}F(X) $ and recall that $ \pi_{\ast}F(X) := \H_{\ast}[N(F(X))] $, where $N$ is the Dold-Kan normalization functor.

Now, any commutative $k$-algebra $A$ and an $A$-module $M$ (viewed as a symmetric bimodule) give rise to a functor $ \fin_{\ast} \rar \ve_k$ that assigns to the set $ [n] $ the vector space $ M \otimes
A^{\otimes \,n}$ and to a pointed map $f: [n] \to [m] $, the action of $f$ on $ M \otimes A^{\otimes \,n} $ given by
$$
f_*(a_0 \otimes a_1 \otimes \ldots \otimes a_n) := b_0 \otimes b_1 \otimes \ldots \otimes b_m\ ,
$$
where $\, b_j := \prod_{i \in f^{-1}(j)} a_i \,$ for $ j =0,1, \ldots m\,$. Following \cite{P1}, we denote this
functor by $ \mathcal L(A,M) $, and for a pointed simplicial set $ X \in \sset_* $, define
$$ \HH_{\ast}(X,A,M)\,:=\, \pi_{\ast}\mathcal L(A,M)(X)\,\text{.}$$
Thus, $\,  \HH_{\ast}(X,A,M) $ is the homology of the complex
$ {\rm C}_\ast(X, A, M) := N[\mathcal L(A,M)(X)] $, which we call the {\it Pirashvili-Hochschild complex} of $A$ with coefficients in $M$  associated to $X$.

\begin{example}
Let $X=\mathbb S^1_\ast $ be the simplicial circle. Recall that the set of $n$-simplices $ \mathbb S^1_n$ can be
identified with the set of monotone sequences of $0$'s and $1$'s of length $n+1$ modulo the identification
$ (0, 0, \ldots, 0) \sim (1,1,\ldots, 1) $ (see Section~\ref{sscirc}). For a nonzero sequence $ x \in  \mathbb S^1_n $,
let $n(x)$ denote the position of the first $1$. The map $ x \mapsto n(x)-1 $ identifies $\mathbb S^1_n$ with $[n]$.
Under this identification, the degeneracy map $\, s_i\,:\,[n] \rar [n+1] $ corresponds to the unique monotone injection skiping $ i+1 $ in its image and the face map $ d_i\,:\,[n] \rar [n-1]$ is given by $d_i(j)=j$ for $j<i$, $d_i(i)=i$ for $i<n$, $d_n(n)=0$ and $d_i(j)=j-1$ for $j>i$. From this description of $ \mathbb S^1_\ast $, it is easy to see that the Pirashvili complex $ {\rm C}_\ast(\mathbb S^1, A, M) $ for $ \mathbb S^1 $ is precisely the classical Hochschild complex $ {\mathrm C}_\ast(A,M) $. Thus, $\HH_{\ast}(\mathbb S^1, A, M) = \HH_{\ast}(A, M)$ for any commutative algebra $A$ and
$A$-module $M$. In a similar way, one can exlicitly describe the Pirashvili complex $ {\rm C}_\ast(\mathbb S^n, A, M) $ for the $n$-dimensional simplicial sphere $ \mathbb S^n_\ast $. The corresponding
homology groups $ \HH_{\ast}(\mathbb S^n, A, M) $ are denoted $ \HH^{[n]}_\ast(A, M) $  and called
{\it the Hochschild homology of $(A,M)$ of order $n$}.
\end{example}

In the present paper, we will mostly deal with two cases: $M=A$ and $M=k$, where in the last case
the module structure on $ k $ comes from an augmentation $ A \to k $. To simplify the notation
we will write $\HH_{\ast}(X,A)$ for $\HH_{\ast}(X,A, A)$ and regard $ X \mapsto \HH_{\ast}(X,A) $ as
a functor on the category of (pointed) simplicial sets assuming $A$  to be fixed.
We will refer to this functor as a {\it higher Hochschild homology} of spaces.

There is another, more conceptual way to define higher Hochschild homology, using homological algebra of functor categories over PROPs. Recall that a PROP is a permutative category $ ({\mathcal P}, \boxtimes) $ whose set of objects is indexed by (or identified with) the natural numbers $ \N $ and whose monoidal structure $ \boxtimes $ is given by addition in $ \N $ (see \cite{Mac}). A $k$-algebra over a PROP  $ {\mathcal P} $ is a strict symmetric monoidal functor from $ {\mathcal P} $ to the tensor category $ \ve_k $.

To define Hochschild homology we take $ {\mathcal P} $ to be a category $ \fin $ of finite sets with monoidal structure given by disjoint union. More precisely, we let $ \fin $ denote the full subcategory of $\set$ whose objects are the sets $\, \un := \{1,2,\ldots, n\} \,$ for $ n \ge 0 $ (where, by convention, $ \underline{0} = \varnothing $) and morphisms are arbitrary set maps.
The monoidal structure on $ \fin $ is given by $ \un \boxtimes \um = \underline{n+m} $. It is well known and easy to prove (see, e.g., \cite[Section~2]{P2}) that the category of $k$-algebras over $\fin $ is equivalent to the category $ \cAlg_k $, the equivalence being given by the functor $ A \mapsto  [(\,\mbox{--}\,\otimes A):\, \un \mapsto A^{\otimes n}] $. We will write $ {\underline A} $ for the algebra over $\fin$ corresponding to the commutative algebra $ A \in \cAlg_k $.

Now, let $\fin$-$\mathtt{Mod}$  (resp., $\mathtt{Mod}$-$\fin$) denote the category of all covariant (resp., contravariant) functors from $\fin$ to the category of vector spaces.  The notation suggests that one should think of
the objects of $\fin$-$\mathtt{Mod}$ and $\mathtt{Mod}$-$\fin$ as left and right $\fin$-modules, respectively. These categories are both abelian with enough projective and injective objects. Furthermore, they are related by a bifunctor
$$ \mbox{--} \,\otimes_{\fin} \mbox{--}\,:\, \mathtt{Mod}\text{-}\fin \times \fin\text{-}\mathtt{Mod} \rar \ve_k $$
that is right exact with respect to each argument, preserves sums and is left balanced (see, e.g., \cite[Sect.~1.5]{P1}). Explicitly, for a right $\fin$-module ${\mathcal N} $ and a left $\fin$-module
${\mathcal M} $,
\begin{equation}
\la{tenspr}
{\mathcal N} \otimes_{\fin} {\mathcal M}\,=\, \big[ \Moplus_{n \ge 0} {\mathcal N}(\un) \otimes_k {\mathcal M}(\un) \big]/R \,,
\end{equation}
where $R$ is the subspace spanned by the vectors of the form ${\mathcal N}(f) x\otimes y - x \otimes {\mathcal M}(f)y\ $
with $ x \in {\mathcal N}(\un)$ and $y \in {\mathcal M}(\um)$ and $ f $ running over all maps in
$ \Hom_{\set}(\um,\un) $.


Next, we consider the functor
\begin{equation}\la{hyo}
 h \,:\, \fin \rar \mathtt{Mod}\text{-}\fin\,,\,\,\,\, \un \mapsto\,  k[\Hom_{\fin}(\mbox{--}, \un)]\ ,
 \end{equation}
where $k[S]$ denotes the vector space generated by a set $S$,
%
and extend \eqref{hyo} to the category simplicial sets in two steps. First, we define a functor $ \set \to \mathtt{Mod}\text{-}\fin $ by taking the left Kan extension of  \eqref{hyo}  along the natural  inclusion $\fin \into \mathtt{Set}$, and then we extend this  degreewise to simplicial sets. Abusing notation, we will continue to denote the resulting functor by
$\, h:\,\mathtt{sSet} \rar \mathtt{s}\mathtt{Mod}\text{-}\fin $. Composing $ h $ with the normalization functor
$ \underline{N}:\, \mathtt{s}\mathtt{Mod}\text{-}\fin \to \mathtt{Ch}_{\ge 0}(\mathtt{Mod}\text{-}\fin) $ assigns to every
simplicial set $ X $ a chain complex  and hence an object in the derived category
$ {\mathscr D}(\mathtt{Mod}\text{-}\fin) $ that we denote  by $ \underline{N}(h(X)) $.

Now, recall that any commutative algebra $A$ defines an algebra ${\underline A} $ over the PROP $ \fin $ that can be viewed as an object
 in $ \fin$-$\mathtt{Mod} $. With this interpretation of $A$,
we have the following result.

\bthm
\la{dtp}
For any $ X \in \sset $ and $ A \in \cAlg_k $, there is a natural isomorphism
$$
\HH_{\ast}(X,A)\,\cong\, \H_{\ast}[\underline{N}(h(X)) \otimes^{\L}_{\fin} \underline A\,]
$$
\ethm
Although Theorem~\ref{dtp} is not explicitly stated in \cite{P1}, it can be deduced from results of this paper.
We do not give a proof of Theorem~\ref{dtp} here as in the next section, we prove the analogous theorem for representation homology
(see Theorem~\ref{rephomdtp}).

\subsection{Representation homology as functor homology}
\la{rhs}
We now define the representation homology of a (reduced) simplicial set by mimicking  Pirashvili's
definition of higher Hochshild homology. Our starting point is the known fact  that the category of commutative Hopf algebras over a field $k$ is equivalent to the category of $k$-algebras of the PROP of finitely generated free groups (see, e.g., \cite[Sect.~5]{P2} and \cite{Ha} for a detailed proof).
To be precise, let $ \ffgr $ denote the full subcategory of $ \Gr $ whose objects are the free groups based
on the sets $ \un = \{1,2, \ldots, n\} $ for $ n \ge 0 $. We denote such groups by
$ \rn := {\mathbb F}\langle \un \rangle $ (where, by convention, $ \langle 0 \rangle $ is the identity group)
and write $ k\rn $ for the corresponding group algebras over $k$. The category $ \ffgr $ is a PROP, with
monoidal product $ \boxtimes $ being the free product of groups, so that $\,\rn \boxtimes \langle m \rangle =
\langle n+m \rangle $. A commutative Hopf algebra $ \mathcal H $ over $k$ defines the (strong monoidal)
covariant functor $\,\ffgr \to \ve_k \,$, $\,\rn \mapsto {\mathcal H}^{\otimes n}\,$, which we denote by
$ \underline{\mathcal H} $. The assignment $\, \mathcal H  \mapsto \underline{\mathcal H} \,$ gives an
equivalence between the category of commutative Hopf algebras over $k$ and the category of $k$-algebras
over the PROP $ \ffgr $. Dually, the category of {\it cocommutative} Hopf algebras is equivalent to
the category of $k$-algebras over the opposite PROP $ \ffgr^{\rm op} $.

Now, observe that for any commutative Hopf algebra $ {\mathcal H} $, the functor
$ \underline{\mathcal H}: \ffgr \to \ve_k $ takes values in the category of commutative algebras,
{\it i.e.}, it can be viewed as a functor $\, \underline{\mathcal H}: \ffgr \to \cAlg_k \,$.
We extend this last functor to the category $ \fgr $ of all free groups  by
taking the left Kan extension along the inclusion $ \ffgr \into \fgr $.
To be precise, let $ \fgr $ denote the category of {\it based}\, free groups whose objects are pairs
$ (\Gamma, S) $, where $ \Gamma = \langle S \rangle $ is a free group with a specified
generating set $S$, and morphisms are arbitrary group homomorphisms $ \Gamma \to \Gamma' $
(not necessarily, preserving the generating sets). We have the  natural inclusion functor
$ i: \ffgr \into \fgr $ that takes $ \rn $ to $ (\rn, \un) $. The Kan extension of
$ \underline{\mathcal H} $ along $i $ then defines a functor
$\, \fgr \to  \cAlg_k\,$ that assigns to the free group $ \frgr{S} $ on a set $S$ the commutative algebra
$ S \otimes {\mathcal H} = \otimes_{s \in S}\, {\mathcal H}_s $.
We continue to denote this functor by $ \underline{\mathcal H} $.

Let $X $ be a reduced simplicial set (or equivalently, a pointed connected  topological space).
Recall that the Kan loop group construction gives a functor $ \lgr{X}: {\Delta}^{\mathrm{op}} \to
\fgr $ that takes $ [n] \in \Delta^{\mathrm{op}} $ to the free group $ \lgr{X}_n = \langle B_n \rangle $ based
on the set $ B_n = X_{n+1}\!\setminus s_0(X_n) $. Now, given a commutative Hopf algebra $ \mathcal H $,
we consider the composition of functors
$$
{\Delta}^{\mathrm{op}} \xrightarrow{\ \lgr{X}\ } \fgr \xrightarrow{\ \underline{\mathcal H}\ } \cAlg_k \ ,
$$
which defines a simplicial commutative algebra $ \underline{\mathcal H}(\lgr{X})$.

\vspace{1ex}

\begin{definition} The {\it representation homology} of $X$ in $\mathcal H$ is defined by
\begin{equation}\la{hrep}
\mathrm{HR}_{\ast}(X, \mathcal H)\,:=\, \pi_{\ast}[\underline{\mathcal H}(\lgr{X})]\,=\,
\H_{\ast}[N(\underline{\mathcal H}(\lgr{X}))]\,\text{.}
\end{equation}
\end{definition}


Clearly, a morphism $f\,:\,X \rar Y$ of reduced simplicial sets induces a map  of graded commutative algebras $\mathrm{HR}_{\ast}(f,\mathcal H) \,:\, \mathrm{HR}_{\ast}(X,\mathcal H) \rar \mathrm{HR}_{\ast}(Y,\mathcal H)$. Thus, representation homology defines a covariant functor $\,\HR(\,\mbox{--}\,, {\mathcal H}): \sset_0 \to \gcAlg_k \,$.
The following proposition justifies the above definition of representation homology.
\bprop \la{rephom0}
Let $G $ be an affine group scheme defined over $k$ with coordinate ring  $\mathcal 
H = \O(G)$. Then, for any $ X \in \sset_0 $, there is a natural isomorphism of graded commutative algebras
\begin{equation}
\la{imhr}
\mathrm{HR}_{\ast}(X, \O(G))\,\cong\, \mathrm{HR}_{\ast}(X, G)\,\text{.}
\end{equation}
In particular, $\, \mathrm{HR}_0(X, \O(G)) \cong \pi_1(X)_G \,$, where
$\pi_1(X)$ is the fundamental group of $ X$.
\eprop
\bproof
If $\mathcal H\,=\,\mathcal O(G)$, we have natural isomorphisms
$\underline{\mathcal H}(\frgr{S}) \cong \Motimes_{s\,\in\,S} \mathcal O(G)_s \cong (\frgr{S})_{G}$ for any set $S$.
This implies that $\,\underline{\mathcal H}(\lgr{X}) \cong (\lgr{X})_{G} $ in $ s \cAlg_k $. On the other hand, by Proposition~\ref{Kanprop}, the simplicial group $\lgr{X}$ is semi-free, and hence a cofibrant object in $\sGr$. This implies that $\, (\lgr{X})_{G} \cong \L (\lgr{X})_{G} $ in $\, \Ho(s \cAlg_k)  $, which, in turn, implies the isomorphism \eqref{imhr} in homology.
The isomorphism for $ \mathrm{HR}_0(X,G) $ is the composition of \eqref{imhr} with \eqref{commpi} and the natural isomorphism of groups $ \pi_0(\lgr{X}) \cong \pi_1(X) $.
\eproof
%
%
%

Let $\Gamma $ be a discrete group, and let $ X = \mathrm{B}\Gamma$ be the classifying space ({\it i.e.}, the simplicial nerve) of $\Gamma$. As a simple
application of Proposition~\ref{rephom0}, we get
\bcor
\la{compwithdrep}
$\mathrm{HR}_{\ast}(\mathrm{B}\Gamma,  \O(G))\,\cong\, \HR_*(\Gamma, G)\,$.
In particular, $\, \mathrm{HR}_0(\mathrm{B}\Gamma,G) \cong \Gamma_{G}\,$.
\ecor
\bproof
The Kan adjunction \eqref{kanl} gives the canonical cofibrant resolution $\lgr{\overline{W}\Gamma} \stackrel{\sim}{\rar} \Gamma \,$ in $\sGr$. Since $ \Gamma $ is discrete, we have $ \overline{W}\Gamma \,=\,\mathrm{B}\Gamma $, and the result follows from Proposition~\ref{rephom0}.
\eproof

\bcor \la{wedgeprod}
For any $X, Y\,\in\,\sset_0$, there is a natural isomorphism
$$\HR_\ast(X \vee Y, G) \,\cong\, \HR_\ast(X, G) \otimes \HR_\ast(Y,G)\ .$$
\ecor

\bproof
Recall that the wedge sum is a (categorical) coproduct in $\sset_0$. Since $\lgr{}$ is a left adjoint functor, we have $\lgr{(X \vee Y)} \,\cong\, \lgr{X} \ast \lgr{Y}$. By Theorem~\ref{hopushout}, it follows that
$$
\L(\lgr{(X \vee Y)})_G  \,\cong\,\L(\lgr{X})_G \otimes \L(\lgr{Y})_G \ .
$$
The desired result is now immediate from  K\"{u}nneth's Theorem and Proposition~\ref{rephom0}.
\eproof

\subsubsection{The fundamental spectral sequence}
Now, we introduce the functor categories $\ffgr$-$\mathtt{Mod}$ and
$\mathtt{Mod}$-$\ffgr$, whose objects are all covariant (resp., contravariant) functors from
$ \ffgr $ to the category of vector spaces. We regard these objects as left and right modules over
$ \ffgr $, respectively. Both categories are abelian with sufficiently many projective and injective
objects. There is a natural bifunctor
$$
\mbox{--} \,\otimes_{\ffgr} \mbox{--}\,:\, \mathtt{Mod}\text{-}\ffgr \times \ffgr \text{-}\mathtt{Mod} \rar \ve_k
$$
which is right exact with respect to each argument, preserves sums and is left balanced in the sense of \cite{CE}.
Explicitly, this bifunctor can be defined by formula \eqref{tenspr} with $ \fin $ replaced by $ \ffgr $.

Since $\, \mbox{--} \,\otimes_{\ffgr} \mbox{--}\,$ is left balanced, the derived functors with respect to each argument are naturally isomorphic, and we denote their common value by
$\mathrm{Tor}^{\ffgr}_{\ast}(\mbox{--},\,\mbox{--})$.
Note that for any left $\ffgr$-module $ {\mathcal M} $, the functor
  $ \, \mbox{--} \,\otimes_{\ffgr} {\mathcal M} :\,\mathtt{Mod}$-$\ffgr \rar \ve_k$ is left adjoint to the functor $\underline{\Hom}({\mathcal M},\, \mbox{--})\,:\,
  \ve_k \rar \mathtt{Mod}$-$\ffgr $, where $\underline{\Hom}({\mathcal M},V)$ is the right $\ffgr$-module
$ \rn \mapsto \Hom_k({\mathcal M}(\rn), V)$ for any vector space $V$. Similarly, for any right $\ffgr$-module $ {\mathcal N}$, the functor ${\mathcal N} \otimes_{\ffgr} \mbox{--}$ is left adjoint to the functor
  $\underline{\Hom}({\mathcal N}, \, \mbox{--})\,:\,\ve_k \rar \ffgr$-$\mathtt{Mod}$. Hence, both functors  $ \mbox{--} \otimes_{\ffgr} {\mathcal M} $ and $ {\mathcal N} \otimes_{\ffgr} \mbox{--}$ commute with colimits.

To state our first theorem we need some notation. First, we recall that
if $\Gamma$ is any group, $k[\Gamma]$ is a cocommutative Hopf algebra: thus, $k[\Gamma]$ defines a
right $\ffgr$-module in  $\mathtt{Mod}$-$\ffgr$.
Now, if $X$ is a reduced simplicial set,  $ k[\lgr{X}] $ defines a simplicial right $\ffgr$-module
in $ \mathtt{sMod}$-$\ffgr$.  Applying the normalization functor
$ \underline{N}:\, \mathtt{sMod}\text{-}\ffgr \to \mathtt{Ch}_{\ge 0}(\mathtt{Mod}\text{-}\ffgr) $ to this simplicial module,
 we get a chain complex of  $\ffgr$-modules and hence an object in the derived category
$ {\mathscr D}(\mathtt{Mod}\text{-}\ffgr) $. Abusing notation, we will denote this object by $ \underline{N}({k[\lgr{X}]}) $.

\bthm
\la{rephomdtp}
For any $ X \in \sset_0 $ and any commutative Hopf algebra $\mathcal H$,
there is a natural isomorphism of graded commutative algebras
$$
\mathrm{HR}_{\ast}(X,\mathcal H) \,\cong\,
\H_{\ast}[\underline{N}({k[\lgr{X}]}) \otimes^{\L}_{\ffgr} \underline{\mathcal H}]\ .
$$
\ethm
To prove Theorem~\ref{rephomdtp} we need a simple lemma.
Recall that for $ n\ge0 $, we denote by $ k\rn $ the group algebra
of the free group based on the set $ \un = \{1,2,\ldots,n\}$. Regarding it as a cocommutative Hopf algebra, we get a right $\ffgr$-module which (to simplify the notation) we also denote by $ k \frgr{n} $.
\blemma
\la{gammangrp}
For each $ n \ge 0 $, the  $\ffgr$-module $\, k \frgr{n} $ is a projective object in $ \mathtt{Mod}\text{-}\ffgr $.
\elemma
\bproof
For a fixed $ n \ge 0 $, let  $ h^n\,:=\,k[\Hom_{\ffgr}(\,\mbox{--}\,,\,\frgr{n})] $ denote the standard right $\ffgr$-module associated to the object $ \frgr{n} \in \ffgr $. By Yoneda Lemma, there is a natural isomorphism $\Hom_{\mathtt{Mod}\text{-}\ffgr}(h^n, {\mathcal N})\,\cong\, {\mathcal N}(\rn) $ for any $ {\mathcal N} \in \mathtt{Mod}\text{-}\ffgr $. The sequence of $\ffgr$-modules
$\, 0 \to {\mathcal N}' \to {\mathcal N} \to {\mathcal N}'' \to 0 \,$ is exact in $ \mathtt{Mod}\text{-}\ffgr $ if and only if the sequence of $k$-vector spaces $\, 0 \to {\mathcal N}'(\rn) \to {\mathcal N}({\rn}) \to {\mathcal N}''(\rn) \to 0 \,$ is exact for all $ n \ge 0 $. It follows that $\, \Hom_{\mathtt{Mod}\text{-}\ffgr}(h^n,\,\mbox{--}\,):\,  \mathtt{Mod}\text{-}\fin \to \ve_k $ is an exact functor, and hence $ h^n $ is a projective object in  $ \mathtt{Mod}\text{-}\ffgr$.
On the other hand, for any $m \ge 0 $, we have
$$
h^n(\frgr{m})\,=\,k[\Hom_{\ffgr}(\frgr{m},\frgr{n})] \,\cong\, k[\frgr{n}^{\times m}] \,\cong\,[k\frgr{n}]^{\otimes m} \,=\,{k\frgr{n}}(\frgr{m})\ ,
$$
which shows that $ k\frgr{n} \cong h^n $ as right $\ffgr$-modules. This finishes the proof of the lemma.
\eproof

\bproof[Proof of Theorem~\ref{rephomdtp}]
By Lemma~\ref{gammangrp}, for any $ n \ge 0 $, ${k \frgr{n}}$ is a projective right $\ffgr$-module such that
${k\frgr{n}} \otimes_{\ffgr} \underline{\mathcal H}\,\cong\, \underline{\mathcal H}(\frgr{n})$.
Since colimits of projective modules are flat and commute with left Kan extensions,
this implies that ${k \frgr{S}}$ is a {\it flat} right $\ffgr$-module
and ${k\frgr{S}} \otimes_{\ffgr} \underline{\mathcal H}\,\cong\, \underline{\mathcal H}(\frgr{S})$ for any set $S$.
Extending the last isomorphism levelwise to simplicial sets, we get an isomorphism of simplicial vector spaces $ {k[\lgr{X}}] \otimes_{\ffgr} \underline{\mathcal H} \,\cong\, \underline{\mathcal H}(\lgr{X}) $.
Further, since each ${k[\lgr{X}_n]}$ is a flat right $\ffgr$-module, the normalized chain complex $ \underline{N}({k[\lgr{X}]}) $
is a complex of flat $\ffgr$-modules, hence we have a natural isomorphism in the derived category $ {\mathscr D}(\mathtt{Mod}\text{-}\ffgr) $:
$$
\underline{N}({k[\lgr{X}]}) \otimes^{\L}_{\ffgr} \underline{\mathcal H}\,\cong\, \underline{N}({k[\lgr{X}]}) \otimes_{\ffgr} \underline{\mathcal H} \,\cong\, N(\underline{\mathcal H}(\lgr{X}))\,\text{.}
$$
At the homology level, this induces the desired isomorphism of Theorem~\ref{rephomdtp}.
\eproof
%

%
%

For our next theorem we recall that the singular chain complex $ C_{\ast}(\Omega X; k)$ of the (Moore) loop space $\Omega X$ of a pointed topological space $X$ has a natural structure of a DG Hopf algebra. The coproduct on $ C_{\ast}(\Omega X; k) $ is induced by the Alexander-Whitney diagonal, while the product comes from the structure of a topological monoid on $\Omega X $ via the Eilenberg-Zilber map  (see, e.g., \cite[Section~26]{FHT1}). Thus, the homology
$ \H_\ast(\Omega X; k) $ of  $ \Omega X $ is a graded {\it cocommutative} Hopf algebra called the {\it Pontryagin algebra} of $X$.

Now, any graded cocommutative Hopf algebra $ {\rm H} $ defines a graded right
$\ffgr$-module $ \underline{\rm H} $ ({\it i.e.}, a contravariant functor from $ \ffgr $ to the category of graded vector spaces). For $ q \in \Z $, we let $ \underline{\H}_{\,q} $ denote the graded component of $\underline{\H}$ of degree $q$; thus,
$ \underline{\H}_{\,q}: \ffgr^{\rm op} \to \ve_k \,$ is a right $\ffgr$-module that assigns $\frgr{n} \mapsto [\H^{\otimes n}]_q $, the $q$-th graded component of the graded vector space $\H^{\otimes n}$. Note that the $\ffgr$-module
$ \underline{\H}_{\,q} $ depends on {\it all} graded components of the Hopf algebra $ \H $, and not solely on $ \H_q $.
With this notation, we can now state our second theorem, which is an analogue of \cite[Theorem 2.4]{P1} for representation homology.

\bthm
\la{rephomdtp2}
There is a natural first quadrant spectral sequence
\begin{equation}
\la{ssrephom}
E^2_{pq}\,=\, \mathrm{Tor}^{\ffgr}_p(\underline{\H}_{\,q}(\Omega X;k), \underline{\mathcal H})\ \underset{p}{\Longrightarrow}\  \mathrm{HR}_{n}(X,\mathcal H)\
\end{equation}
converging to the representation homology of $X$.
\ethm
\bproof
Recall from the proof of Theorem~\ref{rephomdtp} that $\, \underline{N}({k[\frgr{X}]}) $ is a
non-negatively graded chain complex of flat right $\ffgr$-modules. Hence,
for any left $ \ffgr$-module $ \underline{\mathcal H} $, there is a standard
`Hypertor' spectral sequence (see, e.g., \cite[Application~5.7.8]{W}):
\begin{equation*}\la{torseq}
E^2_{pq}\,=\,\mathrm{Tor}^{\ffgr}_p(\H_{q} [\underline{N}({k[\lgr{X}]})],\,\underline{\mathcal H})
\ \underset{p}{\Longrightarrow}\ \H_{p+q}\,[\underline{N}({k[\lgr{X}]}) \otimes_{\ffgr} \underline{\mathcal H}]\ .
\end{equation*}
By Theorem~\ref{rephomdtp}, the limit of this spectral sequence is isomorphic to $ \mathrm{HR}_{\ast}(X,\mathcal H) $. To prove the theorem we need only to show that $\, \H_{\ast} [\underline{N}({k[\lgr{X}]})] \cong \underline{\H}_{\,\ast}(\Omega X;k) \,$ as graded right $\ffgr$-modules.

By Kan's Theorem~\ref{kanstheorem}, $ |\lgr{X}| $ is weakly equivalent to the based loop space $ \Omega X $. In fact, both
$ |\lgr{X}| $ and $ \Omega X $ have natural structures of topological monoids, and they are known to be
weakly equivalent as an $H$-spaces (see, e.g., \cite{Berg}, Sect.~2 and Prop.~3.3(c)).
This implies, in particular, that $ \,\H_{\ast}[N(k[\lgr{X}])]\cong \H_{\ast}(\Omega X;k) $ as graded Hopf algebras, and hence
$ \,\underline{\H}_{\,\ast}[N(k[\lgr{X}])]\cong \underline{\H}_{\,\ast}(\Omega X;k)\, $ as graded  $\ffgr$-modules. Note that
$\,N(k[\lgr{X}])\,$ stands here for the normalized chain complex of
the simplicial Hopf algebra $ k[\lgr{X}] $, while $\, \underline{N}({k[\lgr{X}]})\,$ in the above spectral sequence denotes the normalized chain complex of the simplicial
$\ffgr$-module $k[\lgr{X}]$. We need to check that $\,  \H_{\ast} [\underline{N}({k[\lgr{X}]})] \cong \underline{\H}_{\,\ast} [N({k[\lgr{X}]})] \,$
as graded  $\ffgr$-modules.  Now, the
simplicial $\ffgr$-module ${k[\lgr{X}]}$ assigns to $ \frgr{m} \in \ffgr $ the simplicial vector space $\, k[\lgr{X}_\ast]^{\otimes m}  = \{ k[\lgr{X}_n]^{\otimes m} \}_{n \ge 0} $. By the Eilenberg-Zilber Theorem, the normalized chain complex of this simplicial vector space is homotopy equivalent to
$N(k[\lgr{X}])^{\otimes m}$, while, by Kunneth's formula, the homology of  $\,N(k[\lgr{X}])^{\otimes m}\,$ is naturally isomorphic to $\,\H_{\ast}[N(k[\lgr{X}])]^{\otimes m}$. This shows that $\,\H_{\ast}(\underline{N}({k[\lgr{X}]}))(\frgr{m})\,\cong\,\H_{\ast} [N({k[\lgr{X}]})]^{\otimes m} $ for any $ m \ge 0 $, completing the proof of the theorem.
\eproof

\vspace{1ex}

Theorem~\ref{rephomdtp2} has several interesting implications. First, we consider one important special case when the spectral sequence \eqref{ssrephom} collapses at $ E^2$-term.
\bcor
\la{bgamma}
Let $ \Gamma $ be a discrete group. Then, for any  affine algebraic group $ G $,  there is a natural isomorphism
$$
\HR_{\ast}({\rm B}\Gamma, G) \cong \Tor_{\ast}^{\ffgr}(k[\Gamma], \O(G))\ .
$$
In particular, $ \HR_{0}({\rm B}\Gamma, G) \cong k[\Gamma] \otimes_{\ffgr} \O(G)\,$.
\ecor
\bproof
The classifying space $ X = {\rm B}\Gamma $ is an Eilenberg-MacLane space of type $ K(\Gamma, 1)$.
Its loop space $ \Omega X $ is homotopy equivalent to $ \Gamma $, where $ \Gamma $ is considered
as a discrete topological space. Hence, $ \H_q(\Omega X; k) = 0 $ for all $ q > 0 $,
while $ \H_0(\Omega X; k) \cong k[\Gamma] $ as a Hopf algebra. Thus, for $ X = {\rm B}\Gamma $,
the spectral sequence
\eqref{ssrephom}  collapses on the $p$-axis, giving the required isomorphism.
\eproof

\begin{remark}
Combining the isomorphisms of Corollary~\ref{compwithdrep} and Corollary~\ref{bgamma}, we can
express the representation homology of $ \Gamma $ (originally defined as a non-abelian derived
functor) in terms of classical abelian homological algebra:
$$
\HR_\ast(\Gamma, G)\, \cong\, \Tor_{\ast}^{\ffgr}(k[\Gamma],\O(G)) \ .
$$
In degree 0, we have a natural isomorphism expressing the coordinate ring of the representation variety $ \Rep_G(\Gamma) $ as a functor tensor product:
\begin{equation*}\la{MS}
 \O[\Rep_G(\Gamma)] \cong  k[\Gamma] \otimes_{\ffgr} \O(G)\ .
\end{equation*}
This last isomorphism was found in \cite{KP}, and it was one of the starting points for the present paper.
\end{remark}

\vspace{1ex}

\begin{remark}
The result of Theorem~\ref{rephomdtp2} holds for any (not necessarily, monoidal) left
$\ffgr$-module. In particular, if we take a {\it
reductive} affine algebraic group $G$ and define a left $\ffgr$-module $\O(G)^G \in \ffgr\text{-}\mathtt{Mod}$ by the
formula $\,\frgr{n} \mapsto[\O(G)^{\otimes n}]^G = \O(G \times \stackrel{n}{\ldots} \times G)^G\,$, then, for any $ X \in \sset_0$,  we obtain a homology spectral sequence
\begin{equation}
\la{ssrephom2}
E^2_{pq}\,=\, \mathrm{Tor}^{\ffgr}_p(\underline{\H}_{\,q}(\Omega X;k), \,\O(G)^G) \ {\Longrightarrow}\  \mathrm{HR}_{n}(X,\,G)^G\
\end{equation}
converging to the $G$-invariant part of representation homology of $X$.
The proof of Corollary~\ref{bgamma} shows that, for $ X = {\rm B}\Gamma $, the spectral sequence \eqref{ssrephom2} collapses on the $p$-axis,   giving an isomorphism
$$
\HR_{\ast}({\rm B}\Gamma, G)^G \cong \Tor_{\ast}^{\ffgr}(k[\Gamma], \O(G)^G)\ .
$$
In degree $0$, we therefore have $\,
 \O[\Rep_G(\Gamma)]^G \cong  k[\Gamma] \otimes_{\ffgr} \O(G)^G $.
\end{remark}

\vspace{1ex}

\begin{remark}
Using Corollary~\ref{bgamma}, we can write the 5-term exact sequence associated to the spectral sequence \eqref{ssrephom} in the form
$$
\HR_2(X, G) \to
\HR_2(\pi_1(X),\, G) \to
\underline{\H}_{\,1}(\Omega X; k) \otimes_{\ffgr}  \O(G) \to
\HR_1(X, G) \to
\HR_1(\pi_1(X),\, G) \to 0
$$
If the fundamental group $ \pi_1(X) $ is f.g. virtually free (in particular, finite or f.g. free), then, by \cite[Theorem~5.1]{BRY}, $ \HR_i(\pi_1(X),\, G) $ vanishes for all $ i > 0 $, and hence in this case, we get
\[
\HR_1(X, G) \cong \underline{\H}_{\,1}(\Omega X; k) \otimes_{\ffgr}  \O(G)\ .
\]
\end{remark}

To state futher consequences of Theorem~\ref{rephomdtp2} we introduce some terminology. We will say
that a map $ f: X \to Y $ of pointed topological spaces is a {\it Pontryagin equivalence} (over $k$)
if it induces an isomorphism $ \H_{\ast}(\Omega X; k) \cong \H_{\ast}(\Omega Y;k) $
of Pontryagin algebras (or equivalently, a quasi-isomorphism
$ C_{\ast}(\Omega X; k) \stackrel{\sim}{\to} C_{\ast}(\Omega Y; k) $ of DG Hopf algebras).
The next  result is obtained by applying to \eqref{ssrephom} a
standard comparison theorem for homology spectral sequences (see \cite[Theorem~5.1.12]{W}).
\bcor
\la{rephomhlinv}
If $ f: X \to Y $ is a Pontryagin equivalence, the induced map on
representation homology $ f_*: \mathrm{HR}_{\ast}(X,\mathcal H) \stackrel{\sim}{\to}
\mathrm{HR}_{\ast}(Y,\mathcal H)$ is an isomorphism for any  Hopf algebra $ {\mathcal H} $.
\ecor

We remark that Corollary~\ref{rephomhlinv} does not say that an {\it arbitrary}
isomorphism of Hopf algebras $ \H_{\ast}(\Omega X; k) \cong \H_{\ast}(\Omega Y;k) $
gives an isomorphism $ \mathrm{HR}_{\ast}(X,\mathcal H) \cong \mathrm{HR}_{\ast}(Y,\mathcal H) $.
(Indeed, an abstract isomorphism of Pontryagin algebras need not even induce a map
on representation homology.) Still, Corollary~\ref{bgamma} shows that if both $X$ and $Y$
are aspherical spaces, then
any isomorphism of Pontryagin algebras induces an isomorphism on representation homology.

Next, we recall that the singular chain complex $ C_\ast(X; k) $ of any space $X$
is naturally a DG coalgebra with comultiplication defined by the Alexander-Whitney diagonal. Moreover,
if $X$ is path-connected, there is a quasi-isomorphism of DG coalgebras (see~\cite[Theorem~6.3]{FHT2})
\begin{equation*}
\la{fhthm}
C_{\ast}(X;k)\,\simeq\, \bB[C_{\ast}(\Omega X; k)]\ ,
\end{equation*}
where $\bB$ is the classical bar construction. Since $ \bB $ preserves quasi-isomorphisms,
any Pontryagin equivalence $ f: X \rar Y $ of path-connected spaces
is necessarily a homology equivalence, {\it i.e.}, it induces an isomorphism on singular homology
$ \H_{\ast}(X;k) \stackrel{\sim}{\to} \H_{\ast}(Y;k) $. The converse is not always true
unless $ X $ and $ Y $ are simply-connected. In the latter case, we have the following
well-known result ({\it cf.} \cite[Part I, Prop.~1.1]{Q}).
\blemma\la{equivl}
Let $\,f: X \to Y\,$ be a map of simply-connected pointed topological spaces. The following conditions are equivalent:
\begin{enumerate}
\item[$(1)$] $f$ is a rational homology equivalence: i.e. $ f_*: \H_{\ast}(X; \Q) \stackrel{\sim}{\to} \H_{\ast}(Y; \Q) $;

\item[$(2)$] $f$ is a rational Pontryagin equivalence: i.e. $ f_\ast: \H_{\ast}(\Omega X; \Q)  \stackrel{\sim}{\to} \H_{\ast}(\Omega Y; \Q) $;

\item[$(3)$] $f$ is a rational homotopy equivalence: i.e.
$\, f_\ast: \pi_{\ast}(X) \otimes_{\Z} \Q   \stackrel{\sim}{\to} \pi_{\ast}(Y) \otimes_{\Z} \Q \,$.
\end{enumerate}
\elemma
\bproof
The equivalence $ (1) \Leftrightarrow (2) $ follows a classical theorem of Adams \cite{A} that asserts that, for any simply-connected space $X$, there is a quasi-isomorphism of DG algebras:
$\, C_{\ast}(\Omega X;k) \simeq \cb[C_{\ast}(X;k)]\,$, where $\cb$ is the cobar construction.

To prove that $ (2) \Leftrightarrow (3) $ we first recall that, for any simply-connected $X$,
the $\Q$-vector space $\, L_X := \pi_\ast(\Omega X)_{\Q} \cong \pi_{\ast + 1}(X) \otimes_{\Z} \Q \,$
carries a natural bracket (called the Whitehead product) making it a graded Lie algebra\footnote{The Lie algebra $ L_X $
is called the {\it homotopy Lie algebra} of $X$.}. Thus, a map $ f: X \to Y $ is a rational homotopy equivalence
if and only if it induces an isomorphism of Lie algebras $ f_\ast: L_X \to L_Y $. Then, a classical theorem of Milnor
and Moore  (see \cite[Theorem~21.5]{FHT1}) implies that the Hurewicz homomorphism $\, \pi_\ast(\Omega X) \to \H_\ast(\Omega X; \Q)\,$ induces an isomorphism of graded Hopf algebras
$\,U L_X \stackrel{\sim}{\to} \H_{\ast}(\Omega X; \Q)\,$, where $U(L_X) $ is the universal enveloping algebra of $ L_X$. This yields the equivalence $ (2) \Leftrightarrow (3) $.
\eproof

We say that a map $f: X \to Y $ of simply-connected spaces is a {\it rational homotopy equivalence} if the equivalent conditions of Lemma~\ref{equivl} hold.

\bprop
\la{rephomhtpinv}
A  rational homotopy equivalence induces an isomorphism on representation homology.
Thus, $\HR_{\ast}(X,\mathcal H)$ depends only on the rational homotopy type of $X$.
\eprop
\bproof
By Lemma~\ref{equivl}$(2)$, a rational homotopy equivalence $ X \to Y $ induces an isomorphism
$ \H_{\ast}(\Omega X; \Q)  \stackrel{\sim}{\to} \H_{\ast}(\Omega Y; \Q) $.
Since $ {\rm char}(k) = 0 $, we have $ \Q \subseteq k $, and the Universal Coefficient Theorem implies that
$ \H_{\ast}(\Omega X; k) \cong \H_{\ast}(\Omega Y; k) $. The claim then follows from Corollary~\ref{rephomhlinv}.
\eproof

Next, we look at higher connected spaces. Recall that a space $X$ is called {\it $n$-connected}
if $X$ is path-connected and its first $n $ homotopy groups vanish, {\it i.e.} $ \pi_i(X) = 0 $ for $ 1 \leq i \leq n $.
\bprop
\la{nconnect}
Let $ X $ be an $n$-connected space for some $ n\ge 1$, and let $ \mathcal H = \O(G) $. Then
\begin{equation}
\la{hi}
\HR_q(X, G) =
\left\{
\begin{array}{cll}
k & \mbox{\rm for} & q = 0\\*[1ex]
0 & \mbox{\rm for} & 1 \leq q < n\\*[1ex]
\H_{q+1}(X; \g^*) & \mbox{\rm for} & n \leq q \leq 2n-1
\end{array}
\right.
\end{equation}
where $ \g := {\rm Lie}(G) $ is the Lie algebra of $G$ and $ \g^* $ is its $k$-linear dual.
\eprop
\bproof
If a space $X$ is $n$-connected, its homotopy Lie algebra $ L_X = \pi_*(\Omega X)_{\Q} \cong \pi_{*+1}(X)_{\Q} $ is $n$-reduced, {\it i.e.} $ (L_X)_q = 0 $ for $ 0 \leq q \leq n-1 $. Since $ \H_{\ast}(\Omega X; \Q) \cong U L_X $ and $ \Q \subseteq k $, we have  $\, \H_0(\Omega X; k) \cong k \,$, $\,\H_q(\Omega X; k) = 0 $ for $ 1 \leq q \leq n-1 $, and
\begin{equation*}\la{hiq}
 \H_q(\Omega X; k) \cong  (L_X)_q \otimes_{\Q} k \cong  \pi_{q+1}(X)_{k} \cong \H_{q+1}(X; k)\quad \mbox{\rm for}\quad n \leq q \leq 2n-1\ ,
\end{equation*}
where the last isomorphism is a consequence of the Rational Hurewicz Theorem (see, e.g., \cite{KK}).

Now, recall that for a fixed $ q \ge 0$, the right $\ffgr$-module $ \underline{\H}_{\,q}(\Omega X; k) \,$  is defined
as the functor $ \ffgr^{\rm op} \to \ve_k $, $\, \frgr{m} \mapsto [\H_\ast(\Omega X; k)^{\otimes m}]_q $. It follows from this definition that
\begin{equation*}
 \underline{\H}_{\,q}(\Omega X; k)
 =
\left\{
\begin{array}{lll}
\underline{k}  & \mbox{\rm for} & q = 0\\*[1ex]
0 & \mbox{\rm for} & 1 \leq q \leq n-1\\*[1ex]
\lin_k^* \otimes \H_{q+1}(X; k)  & \mbox{\rm for} & n \leq q \leq 2n-1
\end{array}
\right.
\end{equation*}
where $ \lin_k $ is  the linearization functor:
\begin{equation}
\la{linfun}
\lin_k:\,  \ffgr \to \ve_k\ ,\quad \frgr{m} \mapsto  \
\frgr{m}_{\rm ab} \otimes_{\Z}  k = k^{\oplus m}\ ,
\end{equation}
and $ \lin_k^*: \ffgr^{\rm op} \to \ve_k $ denotes its composition with linear duality.
Thus, for $X$ $n$-connected, the $E^2$-terms of the spectral sequence \eqref{ssrephom} can be identified as
\begin{equation}
\la{van}
E^2_{pq} \cong \left\{
\begin{array}{lll}
k & \mbox{for} & p = 0 \ ,\ q = 0\\*[1ex]
\Tor_{p}^{\ffgr}(\underline{k},\ \underline{\mathcal H}) & \mbox{for} & p > 0 \ ,\ q = 0\\*[1ex]
0 & \mbox{for} & p \ge 0 \ ,\ 1 \leq q < n  \\*[1ex]
\Tor_{p}^{\ffgr}(\lin_k^*,\, \underline{\mathcal H}) \otimes \H_{q+1}(X;k) & \mbox{for} &
p \ge 0 \ , \ n \leq q \leq n-1
\end{array} \right.
\end{equation}
By Lemma~\ref{gammangrp}, the right $\ffgr$-module $ \underline{k}  = k \frgr{0} $ is projective. Hence, $ E^2_{p,0} = 0 $ for  $ p > 0 $.
On the other hand,
$\,\lin_k^* \otimes_{\ffgr} \underline{\mathcal H} \cong \g^* $, while
$\,\Tor_{p}^{\ffgr}(\lin_k^*,\,\underline{\mathcal H}) = 0\,$ for $\, p > 0 \,$.
Hence, for $ n \leq q \leq 2n-1 $, we have
\begin{equation}
\la{van1}
E^2_{0,q}= \g^* \otimes \H_{q+1}(X;k) \cong \H_{q+1}(X; \g^*)\ ,\qquad
E^2_{pq} = 0\ \quad \mbox{for}\ p > 0 \ .
\end{equation}
The vanishing of $ E^2_{pq} $  for all $ p>0 $ in the range $\,0 \leq q \leq 2n-1\,$ shows that the spectral sequence \eqref{ssrephom} collapses on the $q$-axis for these values of $q$. Thus, we have  $\, \HR_q(X, \underline{\mathcal H}) \cong E^2_{0,q} \,$ for  $\,0 \leq q \leq 2n-1\,$. By  \eqref{van} and \eqref{van1}, these are the desired isomorphisms \eqref{hi}.
\eproof

\begin{remark}
Proposition~\ref{nconnect} shows that the representation homology of an $n$-connected space in sufficiently low degrees ($ q \leq 2n-1$) depends only on the Lie algebra $\g$. The main theorem of \cite{BRY2} implies a much stronger result: the {\it whole}  $\, \HR_*(X,G) $ is determined by $ \g $ if $X$ is (at least) $1$-connected. Thus,  for simply connected  spaces, the representation homology with coefficients in an algebraic group $G$ depends only on the connected component $G_0$ of the identity element in $G$: i.e., 
$\, \HR_*(X, G) \cong \HR_*(X, G_0) \,$. This last statement is {\it not} true in general, for non-simply connected spaces: indeed, already in the simplest example $ X = {\mathbb S}^1 $, we have $ \HR_*({\mathbb S}^1,G) \cong \O(G) $.
\end{remark}

\section{Representation homology and higher Hochschild homology}
\la{sect5}
In Section~\ref{rhs}, we defined representation homology  by analogy with
Hochschild homology, using Kan's simplicial loop group construction. In this section, we establish a direct relation between these two homology theories using another classical construction in simplicial homotopy theory due to J.~Milnor  \cite{Mi}.

\subsection{Main theorems} We begin by recalling a standard simplicial model for a (reduced) suspension $\Sigma X$ of a space $X$. The {\it  suspension functor} on pointed simplicial sets is defined by
$$\Sigma :\, \sset_{\ast} \rar \sset_0\,,\,\,\,\,\, X \mapsto C(X)/X\,,$$
where $C(X) \in \sset_{\ast}$ is the reduced cone over $X$. For a pointed simplicial set $X=\{X_n\}_{n \geq 0}$, the set of $n$-simplices in $C(X)$ is given by
$$ C(X)_n\,:=\,\{(x,m)\,:\, x \in X_{n-m}\,,\, 0 \leq m \leq n \}\,, $$
with all $(\ast, m)$ being identified to $\ast$. The face and degeneracy maps in $C(X)$ are defined by
\begin{eqnarray*}
d_i\,:\,C(X)_n \rar C(X)_{n-1}\,, & &(x,m) \mapsto \left\{ \begin{array}{ll}
                                                                   (x,m-1) & \mbox{if $\ 0 \leq i < m$}\\*[1ex]
                                                                   (d^X_{i-m}(x), m) & \mbox{ if $\ m \leq i \leq n$}
                                                                   \end{array} \right. \\*[2ex]
s_j\,:\, C(X)_n \rar C(X)_{n+1}\,, & & (x,m)\mapsto \left\{ \begin{array}{ll}
                                                                    (x,m+1) & \mbox{ if $\ 0 \leq j < m$}\\*[1ex]
                                                                   (s^X_{j-m}(x), m) & \mbox{ if $\ m \leq j \leq n$}
                                                                   \end{array} \right. \\
\end{eqnarray*}
where $\,d_1(x,1)=\ast\,$ for all $x\in X_0$.

The embedding $\, X \hookrightarrow C(X) \,$ is given by $\,x \mapsto (x,0)\,$, and
$\Sigma X $ is defined to be the corresponding quotient set.  Note that, unlike $ C(X)$, the simplicial set  $\Sigma X$ is reduced, since $ (x,0) = \ast $ in $ \Sigma X $ for all $ x \in X $ (in particular, we have
$ C(X)_0 = \{(x,0)\, :\, x \in X_0\} \sim \{\ast\} $). Now, for any pointed simplicial set $X$, there is a homotopy equivalence $ |\Sigma X| \simeq \Sigma|X| $, where $ \Sigma|X| $ is reduced suspension of the geometric realization of $X$ in the usual topological sense.

The next two theorems constitute the main result of this section.
\bthm \la{repvshh1}
For any commutative Hopf algebra $\mathcal H$ and any {\it pointed} simplicial set $X$, there is a natural isomorphism of graded commutative algebras
$$ \mathrm{HR}_{\ast}(\Sigma X, \mathcal H)\,\cong\, \HH_{\ast}(X,\mathcal H;k)\,\text{.}$$
\ethm
To state the next theorem, we recall that there is a natural way to make an arbitrary simplicial set pointed by adding to it a disjoint basepoint.
To be precise, the forgetful functor $\sset_{\ast} \rar \sset$ has a left adjoint $(\, \mbox{--}\, )_+: \sset \to \sset_\ast $ obtained by extending to simplicial sets the obvious functor $X \mapsto X \sqcup \{\ast\}$ on the category of sets. Explicitly, if $\{X_n\}_{n \geq 0}$ is a simplicial set, then $(X_+)_n\,=\,X_n \sqcup\{\ast\}$ for all $n$, and the face and degeneracy maps of $X_+$ are the (unique) basepoint-preserving extensions of the corresponding maps of $X$. Being a left adjoint, the functor $\,(\, \mbox{--}\, )_+\,$ commutes with colimits; in particular, we have
$$
|X_+| \,\cong\, |X|_+\,,
$$
where $|X|_+$ is the space obtained from $|X|$ by adjoining a basepoint.
\bthm \la{repvshh}
For any commutative Hopf algebra $\mathcal H$ and any simplicial set $X$, there is an isomorphism of graded commutative algebras
$$ \mathrm{HR}_{\ast}(\Sigma(X_+), \mathcal H)\,\cong\, \HH_{\ast}(X,\mathcal H)\,\text{.}$$
\ethm
The proofs of Theorem~\ref{repvshh1} and Theorem~\ref{repvshh} are based on a classical simplicial group model of
the spaces $\Omega\Sigma X$, which we now briefly review.

\subsection{Milnor's $FK$-construction} \la{smilnor}
For a pointed simplicial set $ K\,\in\,\sset_{\ast}$, we define $ FK\,:=\,\lgr{\Sigma K}$. Then, by Kan's Theorem~\ref{kanstheorem},  there is a homotopy equivalence of spaces
$$
|FK|\,\simeq\, \Omega \Sigma |K|\,\text{.}
$$
The following observation is due to J. Milnor \cite{Mi} (see also \cite[Theorem~V.6.15]{GJ}).
\blemma[Milnor]
\la{Miln}
For any $ K \in \sset_\ast $,  $\, FK $ is a semi-free simplicial group generated by the simplicial set $ K $ with basepoint identified with $1\,$,
i.e.
$$
FK_n\,=\, (\lgr{\Sigma K})_n \,\cong\, \frgr{K_n}/\langle s_0^n(\ast)=1 \rangle \,\cong\,\frgr{K_n \!\setminus\! s_0^n(\ast)}\,\text{.}
$$
The face and degeneracy maps are induced by the face and degeneracy maps of $K$.
\elemma
\bproof
By definition of the reduced suspension, we have $\,(x,0) = \ast\,$ for all $x \,\in\,K$ and $s_0(x,m)\,=\,(x,m+1)$ for all
$m >0$. Hence, $ (\Sigma K)_{n+1}/s_0(\Sigma K_n)\,=\,\{(x,1)\,|\, x \,\in\, K_n\}$, with $(\ast,1)$ being the basepoint. It follows that
$$
(\lgr{\Sigma K})_n\,=\, \frgr{(\Sigma K)_{n+1}/s_0(\Sigma K_n)}\,\cong\, \frgr{K_n}/(\ast=1)\,\text{.}$$
To calculate the face and degeneracy maps, we recall from Section~\ref{klg} that
$$
d_0^{\lgr{\Sigma K}}(x,1)\,=\, d_1(x,1)\,d_0(x,1)^{-1}\,=\,(d_0 x,1)\,(x,0)^{-1}\,=\,(d_0 x,1)\,,
$$
and $d_i^{\lgr{\Sigma K}}(x,1) = d_{i+1}(x,1) = (d_i x, 1)$  for $i>0$. Similarly, $s_j^{\lgr{\Sigma K}}(x,1) = s_{j+1}(x,1) = (s_j x,1) $ for all
$j \geq 0$. This proves the desired lemma.
\eproof

\subsubsection{Proofs of Theorems~\ref{repvshh1} and Theorem~\ref{repvshh}}
Recall that, for a commutative Hopf algebra $\mathcal H$, we denote by $\underline{\mathcal H} $ the functor $ \mathtt{FGr} \rar \cAlg_k $ on the category of based free groups obtained from the $\ffgr$-module
$\, \langle n\rangle \mapsto {\mathcal H}^{\otimes n}\,$ by taking its left Kan extension along the inclusion $ \ffgr \into \mathtt{FGr} $ (see Section~\ref{rhs}).
\bprop
\la{isofunc}
There is an isomorphism of functors from $\sset$ to $\,\mathtt{s}\cAlg_k\,${\rm :}
$$
\underline{\mathcal H} \circ \lgr \circ \Sigma \circ (\,\mbox{--}\,)_+ \,\cong\, (\,\mbox{--} \,\otimes \mathcal H)\,,
$$
where $\mathcal H$ in the right-hand side is regarded as a commutative $k$-algebra.
\eprop
\bproof
 By Lemma~\ref{Miln}, for any simplicial set $ X = \{X_n\}_{n\ge 0} $, there are natural isomorphisms of groups
 $\, [\lgr{\Sigma( X_+)}]_n\,\cong\, \frgr{X_n}\,$, $\,n \ge 0\,$,
with structure maps on $ \lgr{\Sigma( X_+)}$ being compatible with those of $X$. By applying the functor $ \mathcal H $, we thus get isomorphisms of simplicial commutative algebras
$$
\underline{\mathcal H}([\lgr{\Sigma( X_+)}]_\ast) \,\cong\, \underline{\mathcal H}[\frgr{X_\ast}] \,\cong\, X_\ast \otimes \mathcal H\ ,
$$
which are obviously functorial in $X$. This proves the proposition.
\eproof
Theorem~\ref{repvshh} is an immediate consequence of the above proposition.
To prove Theorem~\ref{repvshh1}, we first note that, although the unreduced cone on a space $X$ coincides with the reduced cone on $X_+$, the corresponding suspensions differ. Instead, for any pointed space $X$, there is a homotopy equivalence (see~\cite[p. 106]{M})
\begin{equation}
\la{may}
\Sigma (X_+) \simeq \Sigma X \vee \mathbb S^1\,\text{.}
\end{equation}
From this we can deduce the following
\blemma \la{unpsusp}
For a pointed topological space $X$, there is a natural isomorphism
$$\mathrm{HR}_{\ast}(\Sigma(X_+), \mathcal H) \,\cong\, \mathrm{HR}_{\ast}(\Sigma X, \mathcal H) \otimes \mathcal H \text{.}$$
\elemma
\bproof
Applying Corollary~\ref{wedgeprod} to \eqref{may}, we have $ \mathrm{HR}_{\ast}(\Sigma(X_+),\mathcal H)\cong \mathrm{HR}_{\ast}(\Sigma X,\mathcal H) \otimes \mathrm{HR}_{\ast}(\mathbb S^1,\mathcal H)$. Now, since $\mathbb S^1 \,\cong\, \Sigma (\mathrm{pt}_+)$, Theorem~\ref{repvshh} implies $\mathrm{HR}_{\ast}(\mathbb S^1,\mathcal H) \cong \HH_{\ast}(\mathrm{pt} ,\mathcal H) \cong \mathcal H$, where $ {\mathcal H} $ is concentrated in degree $0$. It follows that
$\, \mathrm{HR}_{\ast}(\Sigma(X_+),\mathcal H) \cong \mathrm{HR}_{\ast}(\Sigma X,\mathcal H) \otimes \mathcal H\,$ as desired.
\eproof
Lemma~\ref{unpsusp} shows that $\,\mathrm{HR}_{\ast}(\Sigma X,\mathcal H) \cong
\mathrm{HR}_{\ast}(\Sigma(X_+),\mathcal H) \otimes_{\mathcal H} k\,$.
Combining this last isomorphism with that of Theorem~\ref{repvshh}, we now conclude
$$\mathrm{HR}_{\ast}(\Sigma X,\mathcal H)\, \cong\, \mathrm{HR}_{\ast}(\Sigma(X_+),\mathcal H) \otimes_{\mathcal H} k \, \cong \, \HH_{\ast}(X,\mathcal H) \otimes_{\mathcal H} k \,\cong\, \HH_{\ast}(X,\mathcal H;k)\ .
$$
This proves Theorem~\ref{repvshh1}.

\subsection{Examples} We conclude this section with a few simple examples illustrating the use of Theorems~\ref{repvshh1} and~\ref{repvshh}.
More examples will be given in the next two sections.
In what follows, $G$ denotes an arbitrary affine algebraic group and $\g = {\rm Lie}(G) $ stands for its Lie algebra.

\subsubsection{Spheres}
The representation homology of the circle $ {\mathbb S}^1 $ is given
by $ \HR_0({\mathbb S}^1, G) \cong \O(G) $ and $  \HR_i({\mathbb S}^1, G) = 0 $ for $ i > 0 $. This follows, for example, from Lemma~\ref{gammangrp} and Corollary~\ref{bgamma} (since $ {\mathbb S}^1 \cong B \Z \,$). Now, for  higher dimensional spheres, we have
\bprop \la{rephomsph}
$\mathrm{HR}_{\ast}(\mathbb S^n,G)\,\cong\, \Sym_k(\g^{\ast}[n-1])\,$ for all $\,n \geq 2$.
\eprop
\bproof
Note that  $\mathbb S^n \simeq \Sigma\, \mathbb S^{n-1}$  for all $n \geq 2$. By Theorem~\ref{repvshh1}, we conclude
\[
\mathrm{HR}_{\ast}(\mathbb S^n,G) \,\cong \, \HH_{\ast}(\mathbb S^{n-1},\mathcal O(G);k)
             \, \cong \, \Sym_{\mathcal O(G)} (\Omega^1(G)[n-1]) \otimes_{\mathcal O(G)} k
             \, \cong \, \Sym_k(\g^{\ast}[n-1])\,\text{,}\\
\]
where the second isomorphism follows from~\cite[Section 5.5]{P1}.
\eproof

\subsubsection{Suspensions} We now generalize the previous example
to arbitrary suspensions.
\bprop
\la{rephomsusp}
Let $\Sigma X$ be the suspension of a pointed connected space $X$ of finite type. Then
\begin{equation}
\label{isosigma}
\mathrm{HR}_{\ast}(\Sigma X, G)\,\cong\, \Sym_k[\,\rH_\ast(X; \g^*)]\ ,
\end{equation}
where $ \rH_\ast(X; \g^*) $ stands for the reduced $($singular$)$ homology of $X$ with constant coefficients in $ \g^* $.

\noindent
Consequently, by induction,
$$
\mathrm{HR}_{\ast}(\Sigma^n X, G)\,\cong\, \Sym_k \bigl(\,\rH_\ast(X; \g^*)[n-1]\bigr)\ ,\quad \forall\,n\ge 1 \ .
$$
\eprop
\bproof
It is known  (see~\cite[Theorem 24.5]{FHT1}) that  $ \Sigma X $ is  rationally homotopy equivalent to a bouquet of spheres:
$\,\Sigma X\, {\simeq}_{\Q} \,\Mvee_{i \in I}\, \mathbb S^{n_i}\,
$, where each $\mathbb S^{n_i}$ have dimension $n_i \geq 2$.
By Proposition~\ref{rephomhtpinv}, it thus suffices to compute
$ \mathrm{HR}_{\ast}(S, G) $ for  $ S :=  \Mvee_{i \in I} \mathbb S^{n_i} $. Note that the reduced homology $\overline{\H}_{\ast}(S; k)$ of $S$ is isomorphic to $\,\oplus_{i \in I}\, k \cdot v_i \,$ with trivial coproduct, where $v_i $ is a basis element of homological degree $ \deg(v_i) = n_i$. Now, by Corollary~\ref{wedgeprod} and Proposition~\ref{rephomsph}, we have
 \begin{eqnarray*}
 \mathrm{HR}_{\ast}(\Sigma X,G) &\cong & \mathrm{HR}_{\ast}(S , G)
 \, \cong \, \Motimes_{i \in I}\, \mathrm{HR}_{\ast}(\mathbb S^{n_i},G)
 \, \cong \, \Motimes_{i \in I} \Sym_k(\g^{\ast}[n_i -1]) \\
 & \cong & \Sym_k \bigl(\,\Moplus_{i \in I} \g^{\ast}[n_i -1]\,\bigr)
 \, \cong \, \Sym_k\bigl(\Moplus_{n \geq 2} \g^{\ast} \otimes \H_n(\Sigma X;k)[n-1]\,\bigr)\\
 & \cong & \, \Sym_k\bigl(\,\Moplus_{n \geq 1} \g^{\ast} \otimes \H_n(X;k)[n]\bigr)\,\cong\,
\Sym_k\bigl[\,\g^\ast \otimes \rH_\ast(X;k)\bigr]\, \cong \, \Sym_k [\, \rH_\ast(X; \g^*)] \ ,
\end{eqnarray*}
where the last isomorphism is a consequence of the Universal Coefficient Theorem.
\eproof

As a consequence of Theorem~\ref{repvshh}, Lemma~\ref{unpsusp} and Proposition~\ref{rephomsusp}, we have the following general formula for the higher Hochschild homology of $X$ with coefficients in a commutative Hopf algebra.
\bcor \la{pirhhsusp}
For any pointed connected topological space $X$ of finite type,
 \begin{equation} \la{formulaforhh}
 \HH_{\ast}(X,\mathcal O(G))\,\cong\,\Sym_{\mathcal O(G)}\bigl[\,\rH_\ast(X;k) \otimes \Omega^1(G)\, \bigr]\,\text{.}\end{equation}
\ecor
The isomorphism \eqref{formulaforhh} is a refinement of Pirashvili's generalization of the classical
HKR Theorem which (in our notation) asserts that $\, \HH_{\ast}(X,A) \cong \underline{\H}_\ast(X;k) \otimes_{\mathfrak F} \underline{A}\,$ for any smooth commutative algebra $A$ ({\it cf.} \cite[Theorem 4.6]{P1}).
%
%
%
%

\vspace{2ex}

%
\subsubsection{Co-H-Spaces}
The result of Proposition~\ref{rephomsusp} can be seen in a more conceptual way. The key fact is
that the suspension $ \Sigma X $ of any pointed connected space $X$ is a cogroup object in the homotopy category of pointed spaces, with coproduct $\, \Sigma X \to \Sigma X \,\vee\, \Sigma X \,$
given by the natural `pinching' map (see, e.g., \cite[p. 41]{Sp}). The functor
$\,\HR_*(\,\mbox{--}\,,G): \Ho({\tt Top}_{0,*}) \to \gcAlg_{k}\,$ preserves coproducts and hence maps cogroup
objects in $ \Ho({\tt Top}_{0,*}) $ to cogroup objects in $ \gcAlg_{k} $. The latter are precisely the graded commutative Hopf algebras;
thus, the representation homology of  $ \Sigma X $ carries a natural  Hopf algebra structure for any space $X $. Since $ \Sigma X $ is $1$-connected, $\,\HR_*(\Sigma X, G)\,$
is actually a {\it connected} graded commutative Hopf algebra, and hence, by the (dual) Milnor-Moore Theorem
(see \cite[Theorem~0.2]{Fre1}), its underlying algebra structure is free: i.e., $\,\HR_*(\Sigma X, G) \cong \Sym_k V\,$
for some graded vector space $V$. As shown in the proof of Proposition~\ref{rephomsusp}, the rational equivalence
$\,\Sigma X\, {\simeq}_{\Q} \,\Mvee_{i \in I}\, \mathbb S^{n_i}\,$ implies $\, V \cong \rH_\ast(X; \g^*) \,$, and it is easy to see that \eqref{isosigma} is actually an isomorphism of graded Hopf algebras.

The above argument is similar to Berstein's `categorical' proof of the classical Bott-Samelson Theorem describing the Pontryagin algebra of the suspension $ \Sigma X $  (see \cite{Bers}). This formal argument works actually for any (simply-connected associative) co-$\H$-space, provided one replaces the homology $ \rH_\ast(X, k) $
with the so-called Berstein-Scheerer coalgebra $ B_\ast(X,k) $ of $X$ (see, e.g., \cite{Ark}).  In this way, we have the following generalization of Proposition~\ref{rephomsusp}.
\bprop
\la{coHsp}
Let $ X $ be a $1$-connected, associative co-$\H$-space. Then
$$
\HR_\ast(X, G) \cong \Sym_k \bigl[ B_\ast(X, k) \otimes \g^* \bigr]\ ,
$$
where $ B_\ast(X,k) $ is the Berstein-Scheerer coalgebra of $X$.
\eprop
We remark that $\, B_\ast(\Sigma X,k) \cong  \rH_\ast(X, k) \,$ as coalgebras (see \cite{Ark}, p. 1150),
so in the case of suspensions, Proposition~\ref{coHsp} indeed reduces to Proposition~\ref{rephomsusp}.

\section{Examples: surfaces and $3$-manifolds}
\la{sec7}

In this section, using standard topological decompositions, we compute representation homology of some classical non-simply connected spaces. Our examples include closed surfaces (both orientable and non-orientable) as well as some
three-dimensional spaces (link complements in $ \R^3 $, lens spaces and general
closed orientable $3$-manifolds). The representation homology of surfaces and link complements is given in terms of classical Hochschild homology of $ \O(G) $ (or $\O(G^n) $ for some $ n\ge 2$)  with twisted coefficients. The representation homology of
a closed 3-manifold $ M $ is expressed in terms of a differential `Tor', which gives
rise to an (Eilenberg-Moore) spectral sequence converging to $\HR_\ast(M, G)$.

\subsection{Surfaces}
\la{S7.1}
\subsubsection{The torus}
\la{S7.1.1} As a cell complex, the $2$-torus $ \bT^2 = \bS^1 \times \bS^1 $ can be
constructed as the homotopy cofibre (the mapping cone) of the map
$\,\alpha: \bS_c^1 \to \bS_a^1 \vee \bS_b^1 \,$, where the subscripts on the circles indicate
the generators of the respective fundamental groups, and the map itself is specified, up to homotopy, by
its effect on these generators:
\begin{equation}
\la{tor1}
\alpha(c) = [a,b] := a b a^{-1} b^{-1} \ .
\end{equation}
Thus $\,\bT^2 \simeq {\rm hocolim} [\, \ast \leftarrow \bS_c^1 \xrightarrow{\alpha} \bS_a^1 \vee \bS_b^1]\,$,
where the homotopy colimit is taken in the category $ {\tt Top}_{0,\ast} $ of connected pointed spaces.
Applying to this the Kan loop group functor $ \bG $ (more precisely, the composition\footnote{Note that,
being Quillen equivalences, both of these functors preserve homotopy colimits, and hence so does their composition.} of $ \bG $ with the Eilenberg  subcomplex functor $ \overline{S} $, see Section~\ref{qeq}),
we get a simplicial group model for  $ \bT^2 $:
\begin{equation}
\la{tor2}
\bG(\bT^2) \cong {\rm hocolim} [\, 1 \leftarrow \bF_1  \xrightarrow{\alpha} \bF_2  \,]\ .
\end{equation}
Here $ \bF_1 $ and $ \bF_2 $ are the free groups on the generators $c$ and $\{a, b\}$ respectively;
the map $ \alpha $ is given by \eqref{tor1}, and the homotopy colimit is taken in the category $ \sGr $ of
of simplicial groups.

Now, by Theorem~\ref{hopushout}, the derived representation functor preserves homotopy pushouts
for any algebraic group $G$. Hence, it follows from \eqref{tor2} that
\begin{equation}
\la{tor4}
\O[\DRep_G(\bT^2)] \cong {\rm hocolim} [\,k \leftarrow \O(G)  \xrightarrow{\alpha_*} \O(G\times G) \,] \ ,
\end{equation}
where the homotopy colimit is taken in $ s\cAlg_k $,
and the map $ \alpha_*: \O(G) \to \O(G \times G) $ is induced by \eqref{tor1} (explicitly,
$ \alpha_*(f)(x, y) = f([x, y]) $ for $ f \in \O(G)$). Since
$$
{\rm hocolim} \,[\,k \leftarrow \O(G)  \xrightarrow{\alpha_*} \O(G\times G)\,]\, \cong\,
\O(G\times G) \otimes^{\L}_{\O(G)} k\ ,
$$
by Proposition~\ref{rephom0}, we conclude that
\begin{equation}
\la{tor3}
\HR_*(\bT^2, G) \,\cong \, \Tor_*^{\O(G)}(\O(G \times G), \,k)\ ,
\end{equation}
where $ \O(G \times G) $ is viewed as a (right) $\O(G)$-module via the algebra map $ \alpha_* $.

By standard homological algebra (see \cite[ Theorem~2.1, p. 185]{CE}), we can  identify the Tor-groups
in \eqref{tor3} as the classical Hochschild homology of $ \O(G) $ with coefficients in the bimodule
$ \O(G\times G) $, where the right $ \O(G) $-module structure is given via the map $ \alpha_* $ and
the left module structure via the augmentation map $ \varepsilon: \O(G) \to k \,$:
\begin{equation}
\la{torhoch}
 \HR_\ast(\bT^2, G)\,\cong\,\HH_*(\O(G),\, {}_{\varepsilon}\O(G\times G)_{\alpha})\ .
\end{equation}

Alternatively, for classical (matrix) groups $ G $, we can give an explicit `small' DG algebra model for the
representation homology $  \HR_\ast(\bT^2, G) $. Specifically, let $ \m :=\Ker (\varepsilon) $ denote
the maximal (augmentation) ideal of $ \O(G)$ corresponding to the identity element $ e \in  G$.
Assume that $ \m $ is generated by a {\it regular} sequence of elements $ (r_1, r_2, \ldots, r_d) $ in $ \O(G) $,
so that $ d = \dim G $. Consider the free module $ E := \O(G)^{\oplus d} $ and define the $\O$-module
map  $\, \pi:  E \to \O(G) \,$ by $ \pi(f_1, f_2, \ldots, f_d) := \sum_{i=1}^d r_i f_i \,$.
Then, associated to $ (E, \pi) $ is the (global) Koszul complex $ K_\ast(G) := (\Lambda^*_{\O(G)}(E), \delta_K) \,$
with differential
$$
\delta_K(e_0 \wedge e_1 \wedge \ldots \wedge e_n) =
\sum_{i=0}^n \,(-1)^i\, \pi(e_i) \ e_0 \wedge \ldots \wedge \hat{e}_i \wedge \ldots \wedge e_n\ .
$$
Since $ \m $ is generated by a regular sequence, the canonical projection
$\, K_\ast(G) \onto \O(G)/\m \cong k \,$ is a quasi-isomorphism of complexes, and therefore
$ K_*(G) $ is a free resolution of $ k $ over $ \O(G) $. It follows from \eqref{tor3} that
\begin{equation}
\la{tores}
 \Tor_*^{\O(G)}(\O(G \times G), \,k) \,\cong \, \H_*[\cA(\bT^2, G)]\ ,
\end{equation}
where $ \cA(\bT^2, G) := \O(G\times G) \otimes_{\O(G)} K_*(G) $ is a commutative DG algebra with
differential $ d = \id \otimes \delta_K $. In particular,
$\,\HR_i(\bT^2, G) = 0\,$ for all  $\,i > \dim G \,$.

We conclude this example with a conjectural description
of the $G$-invariant part of representation homology $  \HR_\ast(\bT^2, G)^G $. Our conjecture
can be viewed as a multiplicative analogue of the derived Harish-Chandra conjecture proposed in \cite{BFPRW}.

Assume that $G$ is a connected reductive algebraic group of rank $ l \ge 1 $
defined over an algebraically closed field $k $ of characteristic zero.
Let $ T \subset G $ be a Cartan subgroup (i.e., a maximal torus) in $ G$,
and let $ W $ be the corresponding Weyl group. Note that, since $ T $ is commutative,
the map $ \alpha_* : \O(T) \to  \O(T \times T) $ associated to $ T $ factors through the augmentation $ \varepsilon: \O(T) \to k $.
Hence, by \eqref{tor3}, we have canonical isomorphisms
\begin{eqnarray}
\HR_\ast(\bT^2, T) &\cong& \Tor_\ast^{\O(T)}(\O(T \times T),\,k) \la{torT}\\
&\cong& \O(T \times T) \otimes \Tor_\ast^{\O(T)}(k,\,k)\nonumber\\
&\cong& \O(T \times T) \otimes \Lambda^*_k(\m_T/\m_T^2) \nonumber\\
&\cong& \O(T \times T) \otimes \Lambda^*_k(\h^*)\ , \nonumber
\end{eqnarray}
where $ \m_T := \Ker(\varepsilon) $ is the augmentation ideal,
$\, \h = (\m_T/\m_T^2)^* $ is the Lie algebra of $ T $ (i.e., a Cartan subalgebra of $ \g $), and
$ \Lambda^*_k(\h^*) $ is the (homologically) graded exterior algebra with $ \h^* $ placed in degree one.

Now, by functoriality, the natural inclusion $ T \into G $ induces a map of simplicial commutative algebras
\begin{equation}\la{hcmap}
\Phi_G(\bT^2):\ \O[\DRep_G(\bT^2)]^G \,\to \,\O[\DRep_{T}(\bT^2)]^W\ ,
\end{equation}
which is (a multiplicative analogue of) the derived Harish-Chandra homomorphism constructed in \cite{BFPRW}. Then,
the multiplicative version of the derived Harish-Chandra conjecture states

\bconj
\la{conj2}{\it
Assume that $ G $ is one of the classical groups $\, {\rm GL}_n(k)$, $\, {\rm SL}_n(k)$, $\, {\rm Sp}_{2n}(k) \,$, $ n\ge 1\,$, or any  simply-connected, semi-simple\footnote{It is known that every simply-connected reductive affine algebraic group is automatically semi-simple. This follows from two classical facts:
(1) every reductive Lie algebra is a product of a semi-simple one and an
abelian one; (2) there are no nontrivial simply-connected abelian reductive
algebraic groups. } affine algebraic group. Then the derived Harish-Chandra homomorphism \eqref{hcmap} is a weak equivalence in $ s \cAlg_k $. Hence, by \eqref{torT}, there is an isomorphism of graded commutative algebras
\begin{equation}
\la{hhcmap}
 \HR_*(\bT^2, G)^G\, \cong\, [\O(T \times T) \otimes \Lambda^*_k(\h^*)]^W\ .
\end{equation}
}
\econj

We illustrate Conjecture~\ref{conj2} for $ G = \GL_n $. Since $\,
\O(\GL_n) \cong k[x_{ij}, \det(x_{ij})^{-1}]_{1 \leq i,j \leq n} \,$,
the elements $ \{x_{ij} - \delta_{ij}\}_{1 \leq i,j \leq n} $ form a regular sequence in $ \O(\GL_n) $
generating the maximal ideal $ \m $, so we have a canonical commutative DG algebra representing
$ \HR_\ast(\bT^2, \GL_n) \,$:
$$
\cA(\bT^2, \GL_n) \cong k[x_{ij},\,y_{ij},\,\theta_{ij}; \,\det(X)^{-1},\,\det(Y)^{-1}]_{1 \leq i,j \leq n}\ .
$$
Here the variables $ x_{ij} $ and $ y_{ij} $ have homological degree $0$, $\, \theta_{ij} $ have homological
degree $1$, and $ \det(X) $ and $ \det(Y) $ denote the determinants of the generic matrices
$ X := \|x_{ij}\|$ and $ Y := \|y_{ij}\|$. The differential on $ \cA(\bT^2, \GL_n) $ can be written
in matrix terms as
$$
d\,\Theta = XYX^{-1}Y^{-1} - I_n\ ,
$$
where $ \Theta := \|\theta_{ij}\| $ and $ I_n $ is the identity $n\times n$-matrix.
The Harish-Chandra homomorphism
\begin{equation*}
\la{hhcmap1}
\Phi_{\GL_n}(\bT^2):\
\cA(\bT^2, \GL_n)^{\GL_n} \to k\left[x_1^{\pm 1},\,\ldots,\,x_n^{\pm 1},\,y_1^{\pm 1},\ldots,\,y_n^{\pm 1},\,\theta_1,\,\ldots,\,\theta_n\right]^{S_n}
\end{equation*}
is given explicitly (on generators) by the following map
$$
x_{ij} \mapsto \delta_{ij} x_i\ ,\quad y_{ij} \mapsto \delta_{ij} y_i\ , \quad
\theta_{ij} \mapsto \delta_{ij} \theta_i\ ,
$$
and the derived Harish-Chandra conjecture asserts that $ \Phi_{\GL_n}(\bT^2) $ induces an isomorphism ({\it cf.}  \eqref{hhcmap})
\begin{equation}
\la{hhcmap2}
\HR_*(\bT^2, \GL_n)^{\GL_n}\, \stackrel{\sim}{\to}\, k\left[x_1^{\pm 1},\,\ldots,\,x_n^{\pm 1},\,y_1^{\pm 1},\ldots,\,y_n^{\pm 1},\,\theta_1,\,\ldots,\,\theta_n\right]^{S_n}\ ,
\end{equation}
where $ \theta_1, \ldots, \theta_n $ have homological degree $ 1 $ and the symmetric group $ S_n $
acts diagonally by permuting the variables. Note that, in the case of $ \GL_n(k)$, unlike for other algebraic
groups, Conjecture~\ref{conj2} follows from the derived Harish-Chandra conjecture for the corresponding Lie algebra $ \gl_n(k) $
stated in \cite{BFPRW}. This is because the Harish-Chandra map $ \Phi_{\GL_n}(\bT^2) $ can be obtained by formally localizing the derived  Harish-Chandra map for the Lie algebra $ \gl_n(k) $ ({\it cf.} \cite[Sect. 4]{BFPRW}).
In particular, the evidence collected in  \cite{BFPRW} for $ \gl_n(k) $ also supports Conjecture~\ref{conj2} for $ \GL_n(k) $.
we list some of this evidence here.
\begin{enumerate}
\item[1)] Conjecture \ref{conj2} holds for $ \GL_2(k) $ and $ \GL_\infty(k)$.
This follows from Theorem~4.1 and Theorem~4.2(ii) of \cite{BFPRW}.
\item[2)] For all $n\ge 1 $, the map \eqref{hhcmap2} is degreewise surjective. This follows
from \cite[Theorem~4.2(i)]{BFPRW}.
\item[3)] For all $n\ge 1 $, $\, \HR_i(\bT^2, \GL_n)^{\GL_n} = 0 \,$ for $\,i > n\,$.
This follows from \cite[Theorem~27]{BFR}.
\item[4)] For any $ G $ as in Conjecture \ref{conj2}, the map \eqref{hhcmap2} is an isomorphism
in homological degree {\it zero}, i.e. $ \HR_0(\bT^2, G)^G\, \cong\, \O(T \times T)^W $.
This follows from a theorem of Thaddeus \cite{Th} (see also \cite{Sik}).
\end{enumerate}

\vspace{1ex}

Finally, we remark that, for $ G = {\rm GL}_n(k)\,$, $ {\rm SL}_n(k) $ and
$ {\rm Sp}_{2n}(k) $, the Harish-Chandra map  is known to be an isomorphism in homological degree $0$:
$\, \HR_0(\bT^N, G)^G \cong \O(T^N)^W $ for all tori $ \bT^N $, $\, N \ge 2 $ (see \cite{Sik}).
However, by results of \cite[Sect.~5.2]{BFPRW},
the above isomorphism does not extend to higher homological degrees when $\, N \ge 3 $. In other words,
the derived Harish-Chandra homomorphism $ \Phi_{\GL_n}(\bT^N) $ is not a week equivalence
for higher dimensional tori  $ \bT^N $, $\, N \ge 3 $.

\subsubsection{Riemann surfaces}
\la{S7.1.2}
The above computation of representation homology of the $2$-torus naturally generalizes to Riemann surfaces
of an arbitrary genus. To be precise, let $ \Sigma_g $ denote a closed connected orientable surface of
genus $ g \ge 1 $. As a 2-dimensional cell complex, $ \Sigma_g $ can be described as the homotopy cofibre of
the map $\,\alpha^g:\, \bS^1_c \to \vee_{i=1}^g\left(\bS^1_{a_i} \vee \bS^1_{b_i}\right) \,$ defined by
\begin{equation}
\la{riem1}
\alpha^g(c) = [a_1, b_1]\,[a_2, b_2]\,\ldots \,[a_g, b_g] \ ,
\end{equation}
where $a_1, b_1, \ldots, a_g, b_g $ denote the $a$- and $b$-cycles on $ \Sigma_g $
generating the fundamental group $ \pi_1(\Sigma_g, \ast) $.
This gives the simplicial group model $\, \bG(\Sigma_g) \cong
{\rm hocolim} [\, 1 \leftarrow \bF_1  \xrightarrow{\alpha^g} \bF_{2g} \,]\,$ of $ \Sigma_g$,
which, in turn, implies
$$
\O[\DRep_G(\Sigma_g)] \cong {\rm hocolim} \,[\,k \leftarrow \O(G)  \xrightarrow{\alpha^g_*} \O(G^{2g})\,]\, \cong\,
\O(G^{2g}) \otimes^{\L}_{\O(G)} k\ ,
$$
where the map $\, \alpha_*^g: \O(G) \to \O(G^{2g}) $ is defined by
$$
\alpha^g_*(f)(x_1, y_1, \ldots, x_g, y_g) := f([x_1, y_1]\,[x_2, y_2]\,\ldots \,[x_g, y_g])\ ,\quad f \in \O(G)\ .
$$
By Proposition~\ref{rephom0}, we conclude
\begin{equation}
\la{riem3}
\HR_*(\Sigma_g, G) \,\cong \, \Tor_*^{\O(G)}(\O(G^{2g}), \,k) \,\cong\,
\HH_*(\O(G),\, {}_{\varepsilon}\O(G^{2g})_{\alpha})
\end{equation}
where $ {}_{\varepsilon}\O(G^{2g})_{\alpha} $ is the bimodule with left and
and right $\O(G)$-module structure given by the maps $ \varepsilon $ and $\alpha^g_*$,
respectively.

In case when $ \m \subset \O(G) $ is generated by a regular sequence, we can also
express the representation homology of $ \Sigma_g $ as the homology of the commutative DG algebra
$ \cA(\Sigma_g, G) := \O(G^{2g}) \otimes_{\O(G)} K_*(G) $, where $ K_*(G) $ is the global Koszul complex
constructed in Section~\ref{S7.1.1}:
$$
\HR_*(\Sigma_g, G) \,\cong \, \H_\ast[\,\cA(\Sigma_g, G)\,]\ .
$$

Like in the torus case, for a reductive group $G$ with a Cartan subgroup $ T $, there is an algebra map
induced by the derived Harish-Chandra homomorpism $ \Phi_G(\Sigma_G) \,$:
$$
\HR_*(\Sigma_g, G)^G \,\to\, [\O(T^{2g}) \otimes \Lambda^*_k(\h^*)]^W
$$
where $W$ operates diagonally on the target. However, in contrast to the torus case, this map seems
far from being an isomorphism in general. In fact, for $ g \ge 2 $, it is conjectured in \cite{BRY} that $\, \HR_i(\Sigma_g, G) = 0 \,$ if $ i > \dim\,{\mathcal Z}(G) $, where
$ {\mathcal Z}(G) $ denotes the center of $G$; in particular, this implies that
$\, \HR_i(\Sigma_g, G) = 0 \,$ for all $ i > 0 $ if $G$ is semisimple.

\subsection{$3$-Manifolds}
\subsubsection{Link complements in $ \R^3 $}
\la{S7.2}
By a link $ L $ in $ \R^3 $
we mean a smooth (oriented) embedding of the disjoint union
$\, \bS^1 \sqcup \ldots \sqcup \bS^1 \,$ of (a finite number of)
copies of $ \bS^1 $ into $ \R^3 $. The link complement
$ X := \R^3\!\setminus\! L $ is then defined to be the complement of an (open) tubular
neighborhood of the image of $L$ in $ \R^3 $.
To describe a simplicial group model for
$X$ we recall two classical facts from geometric topology  ({\it cf.} \cite{Bir}).
First, by a well-known theorem of J. W. Alexander, every link $L$ in $ \R^3$ can be obtained geometrically as the
closure of a braid $ \beta $ in $ \R^3 $ (we write $ L = \hat{\beta} $
to indicate this relation). Second, for each $n\ge 1$, the braids on $n$
strands  in $ \R^3$ form a group $B_n$ (the Artin braid group),
which admits a faithful representation by automorphisms of the free group
$ \bF_n $ (the Artin representation).  Specifically, the group $B_n$ is generated by $ n-1 $ elements
(``flips'') $\,\sigma_1, \,\sigma_2,\, \ldots, \,\sigma_{n-1} \,$ subject to
the relations
$$
\sigma_i\,\sigma_j = \sigma_j\,\sigma_i\quad (\mbox{if}\ \, |i-j| > 1)\ ,\qquad
\sigma_i\,\sigma_{j}\,\sigma_i = \sigma_{j}\,\sigma_{i}\,\sigma_{j}
\quad (\mbox{if}\ \, |i-j| = 1)\ ,
$$
and in terms of these generators, the Artin representation $ B_n \to \Aut(\bF_n) $ is given by
\begin{equation}
\la{artinact}
\sigma_i\,: \left\{
\begin{array}{lll}
x_i & \mapsto & x_i \,x_{i+1}\, x_i^{-1} \\
x_{i+1} & \mapsto  & x_i\\
x_j & \mapsto & x_j \quad (j \not= i, i+1)
\end{array}\right.
\end{equation}
To simplify the notation we will identify $B_n$ with its
image in $ \Aut(\bF_n) $ under \eqref{artinact}.

The next proposition can be viewed as a refinement of a classical theorem of Artin and Birman
\cite[Theorem~2.2]{Bir} describing the fundamental group  of the link complement $ \R^3\!\setminus\! L $ in terms of the Artin representation (see Remark below).
\bprop
\la{prop1}
Let $ L = \hat{\beta} $ be a link in $ \R^3 $ given by the closure of a braid $ \beta \in B_n $. Then
\begin{equation}
\la{grl}
\bG(\R^3\!\setminus\! L) \,\cong \,{\rm hocolim} \,[\,\bF_n \xleftarrow{\,(\beta,\,\id)\,}
\bF_n \amalg \bF_n \xrightarrow{\,(\id, \,\id)\,} \bF_n\,]\ ,
\end{equation}
where $ \beta $ acts on $ \bF_n $ via the Artin representation \eqref{artinact}.
\eprop
\vspace*{1ex}
\begin{remark}
Note that the homotopy pushout in \eqref{grl} coincides with the homotopy coequalizer $ \,\L\text{coeq} [\!\xymatrix{ \bF_n \ar@<0.4ex>[r]^{\id} \ar@<-0.4ex>[r]_{\beta}  & \bF_n }\!] \,$ of  the two endomorphisms $ \id $ and $ \beta $ of $ \bF_n $. Hence, \eqref{grl} implies
$$
\pi_1( \R^3\!\setminus\! L, \ast) \cong \pi_0[\bG(\R^3\!\setminus\! L)] \cong
\text{coeq} [\!\xymatrix{ \bF_n \ar@<0.4ex>[r]^{\id} \ar@<-0.4ex>[r]_{\beta}  & \bF_n }\!] \cong \langle x_1, \, \ldots\,,\, x_n\ |\ \beta(x_1) = x_1,\,
\ldots,\, \beta(x_n) = x_n \rangle\ ,
$$
which is the Artin-Birman presentation of the link group $ \pi(L) := \pi_1( \R^3\!\setminus\! L, \ast)  $.
\end{remark}
\bproof The proof is based on a simple van Kampen type argument ({\it cf.} \cite{Bir}).
Let's place the $n$-braid $ \beta $ in a regular position in the region $ x < 0 $ in $ \R^3 $, so that
its starting points $ \{p_1, p_2, \ldots, p_n\} $ and end points  $ \{q_1, q_2, \ldots, q_n\} $
are located on the $z$-axis with coordinates $ q_1 <  q_{2}< \ldots  < q_n < p_n < p_{n-1} < \ldots < p_1 $.
The link $L$ is the closure of $\beta$ obtained by joining the points $ p_i $ to $ q_i $ ($ i=1,2, \ldots, n $)
by simple arcs in the region $ x > 0 $, as shown in the picture
\[
\includegraphics[scale=0.3]{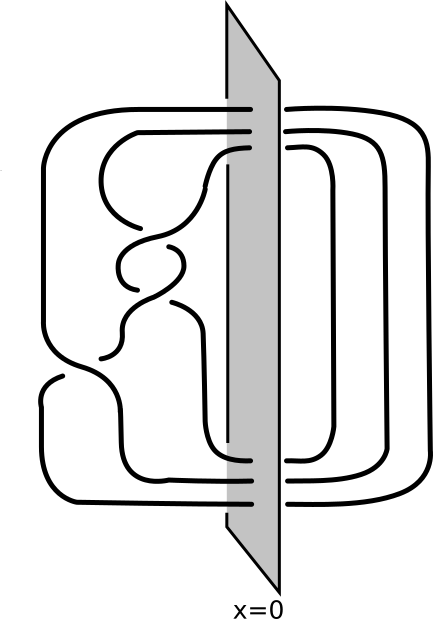}
\]
Now, let $ X := \R^3\!\setminus\! L $ denote the complement of $L$. Define
$$
X_{\ge 0} := \{(x,y,z) \in X\,:\, x\ge 0\}\ ,\quad
X_{\le 0} := \{(x,y,z) \in X\,:\, x\le 0\}\ ,\quad X_0 := X_{\ge 0} \,\cap\,X_{\le 0}\ .
$$
with a (common) basepoint $ \ast $ in $ X_0 $. It is easy to see that $ X_{\ge 0} $ is homeomorphic to the cylinder
over $\,\R^2\!\setminus\!\{p_1, \ldots, p_n\}\,$, which is, in turn, homotopic to
$ \bD^2\!\setminus\!\{p_1, \ldots, p_n\} $, where $ \bD^2 $ is  a two-dimensional disk
in (the $yz$-plane) $ \R^2 $ encompassing the points $ \{p_1,  \ldots, p_n\}  $.
Similarly, we have $\, X_{\le 0} \cong (\R^2\!\setminus\!\{q_1, \ldots, q_n\}) \times [0,1] \simeq  \bD^2\!\setminus\!\{q_1,  \ldots, q_n\} \,$, and
$$
X_0 \,\simeq\, \bD^2\!\setminus\!\{p_1, \ldots, p_n, q_1,  \ldots, q_n\}
\,\simeq\, \bD^2\!\setminus\!\{p_1, \ldots, p_n\}\,\vee\, \bD^2\!\setminus\!\{q_1, \ldots, q_n\}\ .
$$
Under these identifications, the natural inclusions $\,X_{\le 0} \hookleftarrow X_0 \into X_{\ge 0} \,$ can be identified with
\begin{equation}
\la{ddiag}
\bD^2\!\setminus\!\{q_1, \ldots, q_n\} \xleftarrow{(f_\beta, \,\id)}
\bD^2\!\setminus\!\{p_1, \ldots, p_n\}\,\vee\, \bD^2\!\setminus\!\{q_1, \ldots, q_n\} \xrightarrow{(\id, \,f_e)}
\bD^2\!\setminus\!\{p_1, \ldots, p_n\}\ ,
\end{equation}
where the map $ f_\beta $ is determined (uniquely up to homotopy) by the braid $ \beta $ and the map
$ f_e $ is determined by the trivial braid connecting the points $p_i$ and $ q_i $. Thus, we can represent
$X$ in $ \Ho({\tt Top}_{0,*})\,$ as the homotopy pushout of the diagram \eqref{ddiag}.

Next, recall that $ B_n $ can be identified with the mapping class group of
$\, \bD^2\!\setminus\! \{p_1, \ldots, p_n\} \,$ comprising (the isotopy classes of) orientation-preserving
homeomorphisms that fix pointwise the boundary of $ \bD^2 $. As a mapping class group, $B_n$ acts
naturally on the fundamental group $ \pi_1(\bD^2\!\setminus\! \{p_1, \ldots, p_n\}, \ast) $ and the latter
can be identified with the free group $ \bF_n $ on generators $ x_1, \ldots, x_n $ represented by
small loops in $ \bD^2\!\setminus\! \{p_1, \ldots, p_n\} $ around the points $ p_i $. It is well-known
(see \cite{Bir}) that the action of $ B_n $ on $ \bF_n $ arising from this construction is precisely the Artin
representation \eqref{artinact}. Now, using the map $ f_e $ we identify $ \bD^2\!\setminus\!\{q_1, \ldots, q_n\} $
with $ \bD^2\!\setminus\!\{p_1, \ldots, p_n\} $ in \eqref{ddiag} and apply the loop group functor to this diagram of spaces. As a result, we get the equivalence
\eqref{grl}, which completes the proof of the proposition.
\eproof
To state our main theorem we introduce some notation. First, observe that,
for any algebraic group $G$, the Artin representation $ B_n \into \Aut(\bF_n) $ induces naturally
a braid group action $ B_n \to \Aut[\O(G^n)] $, which we denote by $ \beta \mapsto \beta_\ast $.
On the standard generators, this action is defined by
$$
(\sigma_i)_\ast:\ \O(G^n) \to \O(G^n)\ ,\quad
f(g_1,\ldots, g_i, g_{i+1}, \ldots, g_n) \mapsto f(g_1,\ldots, g_i g_{i+1} g_{i}^{-1}, g_{i}, \ldots, g_n)\ .
$$
Now, for a braid $ \beta \in B_n $, we let $ O(G^n)_{\beta} $ denote the $ \O(G^n)$-bimodule whose
underlying vector space is $ \O(G^n) = \O(G)^{\otimes n} $, the left action of $ \O(G^n) $ is given by  multiplication, while the right action is twisted by the automorphism $ \beta_* $.
\bthm
\la{hrlink}
Let $\, L = \hat{\beta}\, $ be a link in $ \R^3 $ given by the closure of a braid $ \beta \in B_n $. Then
$$
\O[\DRep_G(\R^3 \!\setminus\! L)] \cong \O(G^n) \otimes_{\O(G^{2n})}^{\L} \O(G^n)_{\beta}\ .
$$
Consequently,
\begin{equation}
\la{HRRL}
\HR_\ast(\R^3\!\setminus\! L, G) \,\cong \,\HH_\ast(\O(G^n),\,\O(G^n)_{\beta})\ .
\end{equation}
\ethm
\bproof
By Proposition~\ref{prop1} and Theorem~\ref{hopushout}, we have
\begin{eqnarray*}
\O[\DRep_G(\R^3\!\setminus\! L)] &\cong & {\rm hocolim}\,[\,\O(G^n) \xleftarrow{\,(\beta_*,\,\id)\,} \O(G^n) \otimes_k \O(G^n) \xrightarrow{\,(\id,\,\id)\,} \O(G^n)\,]\\
                                  & \cong & {\rm hocolim}\,[\,\O(G^n) \xleftarrow{\,(\beta_*,\,\id)\,} \O(G^{2n}) \xrightarrow{\,(\id,\,\id)\,} \O(G^n)\,]\\
                                  & \cong & \O(G^n) \otimes_{\O(G^{2n})}^{\L} \O(G^n)_\beta \ .
\end{eqnarray*}

This completes the proof of the theorem.
\eproof

\begin{remark}
Theorem~\ref{hrlink} exhibits an interesting analogy between the representation homology of link complements in $ \R^3 $
and their Legendrian contact homology in the sense of L.~Ng (see \cite{Ng1}). This analogy is explained in the recent paper
\cite{BEY}, where a new algebraic construction of link contact homology is given. Roughly speaking, in terminology of
\cite{BEY}, $\, \O[\DRep_G(\R^3\!\setminus\! L)] \,$ represents the algebraic `homotopy closure' of the braid
$ \beta \in B_n $ in the category of simplicial commutative algebras, while the Legendrian contact homology of $ \R^3\!\setminus\! L $ can be computed from a certain  DG category $ \mathscr{A}_{L} $
that represents the homotopy braid closure of $ \beta \in B_n $
in the category of (small pointed) DG $k$-categories.
\end{remark}

\subsubsection{Link complements in $ \bS^3 $}
Note that Theorem~\ref{hrlink} computes the representation homology
of the topological {\it  space} $\,\R^3\!\setminus\! L\,$, not of the {\it link
group} $ \pi(L) $, which is the fundamental group of
$ \R^3\!\setminus\! L $. Even when $ L $ is a knot in $ \R^3 $
(i.e., a link with one component), the representation homologies $ \HR_\ast(\R^3\!\setminus\! L, G) $
and $ \HR_\ast(\pi(L), G) $ differ, because $ \R^3\!\setminus\! L $ is not a
$ K(\pi, 1)$-space ({\it cf.} Example~\ref{unknot} below).
In knot theory, one is usually interested in representation varieties
of the knot group $ \pi(L) $, so it is important to understand the relation between
 $ \HR_\ast(\R^3\!\setminus\! L, G) $ and $ \HR_\ast(\pi(L), G) $.
A natural way to approach this problem is to consider $ L  $ as a link in $ \bS^3 $ by  adding to $ \R^3 $ one point at infinity. If $ L \subset \R^3 \subset \bS^3 $ is a knot,
by Papakyriakopoulos' Sphere Theorem, the complement $ \bS^3\!\setminus\! L $ is an aspherical space, and $\, \pi_1(\bS^3\!\setminus\! L, \ast) \cong  \pi_1(\R^3\!\setminus\! L, \ast) = \pi(L) \,$.
Hence, for any knot $L$,  $\,\HR_\ast(\pi(L), G) \cong \HR_\ast(\bS^3\!\setminus\! L, G)\,$,
so it suffices to clarify the relation between $ \HR_\ast(\R^3\!\setminus\! L, G) $ and $ \HR_\ast(\bS^3\!\setminus\! L, G) $.

To this end, we observe that the natural inclusion $ \R^3\!\setminus\! L \into \bS^3 \!\setminus\! L $ fits into the cofibration sequence
$\, \bS^2 \stackrel{i}{\into} \R^3\!\setminus\! L \into \bS^3 \!\setminus\! L \,$,
so that
\begin{equation}
\la{S3L}
\bS^3 \!\setminus\! L \cong {\rm hocolim}\,[\,\ast \,\leftarrow\, \bS^2\, \xrightarrow{i}\, \R^3 \!\setminus\! L\,]\ ,
\end{equation}
where $ \bS^2 \subset \R^3  $ is chosen in such a way that it encloses $ L $ in $ \R^3 $.
Applying the Kan functor to \eqref{S3L}, we get
\begin{equation}
\la{S3L1}
\bG(\bS^3 \!\setminus\! L) \cong {\rm hocolim}\,[\,1 \,\leftarrow\, \bG(\bS^2)\,
\xrightarrow{i_*}\, \bG(\R^3 \!\setminus\! L)\,]\ .
\end{equation}
To describe the induced map $ i_* $, we note that
$\,\bS^2 \cong \Sigma \,\bS^1 \cong {\rm hocolim}\,[\,\ast \leftarrow \bS^1 \rightarrow \ast\,]\,$, hence
\begin{equation}
\la{S3L2}
\bG(\bS^2) \cong {\rm hocolim}\,[\,1 \,\leftarrow\, \bF_1\,
\rightarrow 1\,]\ .
\end{equation}
Now, if we identify $ \bG(\R^3 \!\setminus\! L) $ as in Proposition~\ref{prop1}, then $ i_* $  is determined  by the  morphism of diagrams
\begin{equation}
\la{diagmap}
\begin{diagram}[small, tight]
1     & \lTo & \bF_1              & \rTo & 1 \\
\dTo  &      & \dTo_{}       &      & \dTo \\
\bF_n & \lTo^{\ (\beta, \id)\ } & \bF_n \amalg \bF_n & \rTo^{\ (\id, \id)\ } & \bF_n
\end{diagram}
\end{equation}
where the map in the middle is given  (on free generators) by
$\, x \mapsto (x_1 \,x_2\, \ldots\, x_n)\,(y_1\,y_2\,\ldots\,y_n)^{-1}$.
Note that the left square in \eqref{diagmap} commutes because
the product $ x_1 x_2 \ldots x_n \in \bF_n $ stays fixed under the Artin representation for any $ \beta \in B_n $.

The map $ i_*: \bG(\bS^2) \to \bG(\R^3 \!\setminus\! L) $ induces a map of simplicial commutative algebras
\begin{equation}
\la{drepi}
i_*:\ \O[\DRep_G(\bS^2)] \to \O[\DRep_G(\R^3 \!\setminus\! L)]\ ,
\end{equation}
which (to simplify the notation) we denote by the same symbol.
By Theorem~\ref{hrlink},
$$
\O[\DRep_G(\R^3 \!\setminus\! L)] \,\cong \,\O(G^n) \otimes_{\O(G^{2n})}^{\L} \O(G^n)_{\beta}\ .
$$
On the other hand, by \eqref{S3L2},
$$
\O[\DRep_G(\bS^2)] \,\cong\, k \otimes_{\O(G)}^{\L} k \, \cong \,
\Lambda^\ast_k(\m/\m^2) \cong \Lambda^\ast_k(\g^*)\ ,
$$
where $ \Lambda^\ast_k(\g^*) $ denotes the graded exterior algebra of $ \g^*$, with $ \g^* $ being in degree one,
equipped with trivial differential. With these identifications,
the map \eqref{drepi} is induced by the algebra homomorphism
\begin{equation}
\la{alghom}
\O(G) \to \O(G^{2n}) \ ,\quad f(x) \mapsto   f((x_1 \,x_2\, \ldots\, x_n)\,(y_1\,y_2\,\ldots\,y_n)^{-1})\ .
\end{equation}
Now, we can regard $\, \O(G^n) \otimes_{\O(G^{2n})}^{\L} \O(G^n)_{\beta}\,$ as a DG module over the DG algebra $   k \otimes_{\O(G)}^{\L} k \cong \Lambda^\ast(\g^*)  $. As a consequence
of  \eqref{S3L1} and Theorem~\ref{hrlink}, we have then the following
\bthm
\la{HRSL}
Let $\, L = \hat{\beta}\, $ be a link in $ \bS^3 $ given by the closure of a braid $ \beta \in B_n $. Then
$$
\O[\DRep_G(\bS^3 \!\setminus\! L)] \,\cong\, k \otimes^{\L}_{ \Lambda^\ast(\g^*)} [\O(G^n) \otimes_{\O(G^{2n})}^{\L} \O(G^n)_{\beta}]\ .
$$
Consequenlty, there is a natural spectral sequence
\begin{equation*}
E^2_{\ast ,\ast} = \Tor_\ast^{\Lambda^\ast(\g^*)}(k,\,\HH_\ast(\O(G^n),\,\O(G^n)_{\beta})\ \Longrightarrow\ \HR_\ast(\bS^3 \!\setminus\! L,\,G)\ ,
\end{equation*}
converging to the representation homology of $\,\bS^3 \!\setminus\! L\,$.
\ethm
\begin{example}
\la{unknot}
Let $ L = \bigcirc $ be the unknot in $ \R^3 $. We can represent $ L $ by the trivial braid $ \beta = 1 \in B_1 $.
In this case, Theorem~\ref{hrlink} combined with the
classical Hochschild-Kostant-Rosenberg Theorem gives
\[
\HR_\ast(\R^3 \!\setminus\! \bigcirc,\,G) \cong \HH_\ast(\O(G), \O(G)) \cong \Omega^*(G) \ ,
\]
where $ \Omega^*(G) $ is the de Rham algebra of (algebraic) differential forms on the group $G$. On the other hand, $\, \HR_\ast(\bS^3 \!\setminus\! \bigcirc,\,G) \cong \HR_\ast(\pi(\bigcirc),\,G) \cong \O(G)\,$, since
$ \pi(\bigcirc) \cong \Z \,$. This simple example illustrates the fact that representation homology does depend on the higher homotopy structure of a space: in particular, it distinguishes the link complements in $ \R^3 $ and $ \bS^3 $, even though their fundamental groups are the same.

\end{example}

\subsubsection{Lens spaces}
\la{S7.3} Recall that, for coprime integers $p$ and $q$, the lens
space $ L(p,q) $ of type $(p,q)$ is defined as the quotient $\, \bS^3/\Z_p \,$ of the $3$-sphere $\,\bS^3 $ viewed as the unit sphere in $ \c^2 \,$ modulo the (free) action of the cyclic group $ \Z_p $ given by
$ (z,w) \mapsto (e^{2\pi i/p}\, z,\, e^{2\pi i q/p}\, w)\,$.  This definition shows that $ L(p,q) $ is a compact connected $3$-manifold, whose universal
cover is $ \bS^3 $ and the fundamental group is $ \Z_p $. Special cases include
$ L(1,0) \cong \bS^3 $, $ L(0,1) \cong \bS^1 \times \bS^2 $ and $ L(2,1) \cong \R \bP^3 $.

To compute the representation homology of  $L(p,q)$ we will use a well-known topological construction  of these spaces via Dehn surgery  in $ \bS^3 $ (see, e.g., \cite[Chap.~3B]{Rolf}). Recall that if $ K \subset \bS^3 $ is a
knot in $ \bS^3 $ and $ p, q $ are two integer numbers, the $ p/q $ {\it Dehn surgery} on $ K $ is a $3$-dimensional
space obtained by removing from $ \bS^3 $ the interior $ \mathring{N}(K) $ of a regular tubular neighborhood $ N(K)$, which is a $3$-dimensional solid torus $\, \bS^1 \times \bD^2 $, and then gluing $ \bS^1 \times \bD^2 $
back to $ \bS^3\! \setminus \!\mathring{N}(K) $ in such a way that the meridional curve of $ \bS^1 \times \bD^2 $
is identified with a $(p,q)$-curve on the boundary of $\,  \bS^3\! \setminus \!\mathring{N}(K) $.
For the trivial knot $ K \subset \bS^3 $, it is easy to see that the $p/q$ Dehn surgery on $ K $ gives precisely
the lens space $ L(p,q) $. In this case, the knot complement $ \bS^3\! \setminus \!\mathring{N}(K) $ is
homeomorphic to the solid torus $ \bS^1 \times \bD^2 $, so the space $ L(p,q) $ can be obtained by gluing together two solid tori along their boundary.

To describe this in more concrete terms, we consider the solid torus  $ \bS^1 \times \bD^2 $ as a subset in $ \c^2 \,$:
$$
\bS^1 \times \bD^2 = \{(z,w) \in \c^2\, : \, |z| = 1\,,\, |w| \leq 1 \}\ .
$$
We identify  $\,\bT^2 = \bS^1 \times \bS^1\,$ as the boundary of $ \bS^1 \times \bD^2 $ in $ \c^2 $ and denote by $ i: \bT^2 \into \bS^1 \times \bD^2 $ the natural inclusion.

Now, for the given pair $ (p,q) $ of coprime numbers, we choose $ m,n \in \Z \,$,
so that $\,mq-np = 1\,$, and define the `gluing' map
$\, \gamma: \,\bT^2 \to \bS^1 \times \bD^2\,$ by\footnote{Note that the map $ \gamma $ factors as $\, \gamma = i\circ \gamma_A \,$, where
$\, \gamma_A: \bT^2 \to \bT^2 \,$ is a homeomorphism of the torus defined by  formula \eqref{gluemap}
that depends on the matrix $ A = \left(\!\begin{array}{ll}
m & p\\
n & q
\end{array}\!\right) \in \SL_2(\Z)$. The assignment $ A \mapsto \gamma_A $  induces
an isomorphism from $ \SL_2(\Z) $ to the mapping class group of  $ \bT^2 $. Formula \eqref{Lpush1} shows that $ L(p.q) $ can be obtained as a
join of two circles, $\, \bS^1 \star_{\gamma} \bS^1 \,$,  twisted by $ \gamma_A $.}
\begin{equation}
\la{gluemap}
\gamma(z,\,w) := (z^m w^p,\,z^n w^q)\ .
\end{equation}
Then the $p/q$ Dehn surgery construction of $L(p,q) $ can  be described as the pushout in $ {\tt Top}_{0,\ast} $:
\begin{equation}
\la{Lpush}
L(p,q) \cong \colim\,[\,\bS^1 \times \bD^2 \,\stackrel{\,i\,}{\hookleftarrow}\, \bT^2 \,\stackrel{\,\gamma\,}{\into}\, \bS^1 \times \bD^2\,]\ .
\end{equation}
Since $ i $ is a cofibration in  $ {\tt Top}_{0,\ast} $, we can replace the colimit in \eqref{Lpush} by a homotopy colimit and then replace the
diagram of solid tori by a homotopy equivalent diagram of circles:
\begin{equation}\la{Lpush1}
L(p,q) \,\cong\,
{\rm hocolim}\,[\,\bS^1  \stackrel{\pi}{\twoheadleftarrow}\, \bT^2 \,\stackrel{\,\bar{\gamma}\,}{\to}\, \bS^1 \,]\  .
\end{equation}
In this diagram, the map $\pi$ is given by the canonical projection $ (z,w) \mapsto z $ and  $\bar{\gamma} $
is the composition $ \pi\circ\gamma $ defined by $ (z,w) \mapsto z^m w^p $. Now, applying the Kan loop group functor, we get a simplicial group model for $ L(p,q) $:
\begin{equation}
\la{sgrL}
\bG[L(p,q)] \,\cong \, {\rm hocolim}\,[\,\bF_1 \xleftarrow{\pi} \bG(\bT^2) \xrightarrow{\gamma} \bF_1\,]\ .
\end{equation}
Recall (see \eqref{tor2}) that $ \bG(\bT^2) $ is given in $ \sGr $ by the homotopy cofibre of  the commutator map $ \alpha: \bF_1 \to \bF_2 $,
$\, c \mapsto [a,b] \,$, where  $ a $ and $b$ are the generators  of $ \bF_2 $ corresponding to the meridian and longitude in $ \bT^2 $.
In terms of these generators, the maps $ \pi $ and $ \gamma $ in \eqref{sgrL} are  induced by
\begin{equation}
\la{maps1}
\pi: \, \bF_2 \to \bF_1\ ,\ (a,b) \mapsto (z,1)\ ,\qquad
\gamma: \, \bF_2 \to \bF_1\ ,\ (a,b) \mapsto (z^m, z^p)\ ,
\end{equation}
where  $ z  $ is a generator of $ \bF_1 $.

Now,  assume that $ G $ admits a global Koszul resolution
$  K_\ast(G) $  described in Section~\ref{S7.1}. Then, we have
an explicit DG algebra model for $ \HR_\ast(\bT^2, G) $
given by $ \cA_*(\bT^2, G) = \O(G \times G) \otimes_{\O(G)} K_\ast(G) $.
Applying  to \eqref{sgrL} the derived representation functor, we get
\begin{equation}
\la{sgrL1}
\O[\DRep_G(L(p,q))] \,\cong\, {\rm hocolim}\,[\,\O(G) \xleftarrow{\pi_*} \cA_\ast(\bT^2, G) \xrightarrow{\gamma_*} \O(G)\,]\ ,
\end{equation}
The maps $\pi_*$ and $\gamma_*$ in \eqref{sgrL1} are  determined by
\eqref{maps1};  on the degree $0$ component of the DG algebra $\cA_*(\bT^2, G) $, they are given  by
\begin{equation}
\la{maps3}
\pi_*: \, \O(G\times G) \to \O(G)\ ,\quad f(x,y) \mapsto f(z,e)\ ,
\end{equation}
\begin{equation}
\la{maps4}
\gamma_*:\,\O(G \times G) \to \O(G)\ ,\quad f(x,y) \mapsto f(z^m, z^p)\ .
\end{equation}
Using  $\pi_*$ and $\gamma_*$, we can make $ \O(G) $ into  (left and right) DG modules over the DG algebra  $\cA_*(\bT^2, G) $, which we denote by $ \O(G)_\pi $ and $ \O(G)_\gamma $ respectively.
With this notation, we have the following result that completes our calculation.
\bthm
\la{lensspaces}
The representation homology of a $3$-dimensional lens space $ L(p,q)  $ is given by
$$
\HR_\ast(L(p,q), G)\, \cong \, \boldsymbol{\mathrm{Tor}}_\ast^{\cA_*}(\cO(G)_\pi,\,\cO(G)_\gamma)\ ,
$$
where $ \boldsymbol{\mathrm{Tor}}_\ast^{\cA_*} $ denotes the differential \,{\rm Tor} taken over the DG algebra $ \cA_* = \cA_*(\bT^2, G) $.
In particular, there is an Eilenberg-Moore homology spectral sequence
$$
E^2_{\ast, \ast} = \Tor_*^{\HR_\ast(\bT^2, G)} (\O(G)_\pi,\,\O(G)_\gamma)\ \Longrightarrow\
\HR_\ast(L(p,q), G)
$$
converging to the representation homology of $ L(p,q) $.
\ethm
\subsubsection{Closed $3$-manifolds} The above construction of lens spaces generalizes to arbitrary closed $3$-manifolds. Specifically, it is well known that every closed connected
orientable $3$-manifold $ M $ admits a Heegaard decomposition $\,H_g \,\cup_{\gamma}\, H_g\,$ which can be written
as
\begin{equation}
\la{HeeG}
M \,\cong\, \colim\,[\,H_g \,\stackrel{i}{\hookleftarrow}\, \Sigma_g\, \stackrel{\gamma}{\hookrightarrow}\, H_g\,]\ ,
\end{equation}
where $ H_g $ is a handlebody of genus $ g \ge 0 \,$, $ i $ is the natural
inclusion identifying $ \Sigma_g = \partial H_g \,$ and $ \gamma $ is
a gluing map defined as the composition $\, \Sigma_g \xrightarrow{\gamma_A} \Sigma_g \stackrel{i}{\hookrightarrow} H_g \,$, where $ \gamma_A $ is an (orientation-preserving) diffeomorphism of $\Sigma_g $
representing an element in the mapping class group $ {\mathcal M}(\Sigma_g) := \pi_0({\rm Diff}^{+}\, \Sigma_g) $. In particular, for
$ g = 1 $, the Heegaard diagram \eqref{HeeG} becomes \eqref{Lpush}; in fact, the lens spaces can be characterized as
(closed) $3$-manifolds that admit Heegaard decompositions of genus $ 1 $.

Since $ H_g $ is homotopy equivalent as a cell complex to the bouquet of $g$ circles $ \vee_{i=1}^g \bS^1 $, we can represent
the homotopy type of $ M $ by
$$
M\, \cong \,{\rm hocolim}\,\left[\,\Mvee_{i=1}^g \bS^1 \,\leftarrow\, \Sigma_g \,\rightarrow\, \Mvee_{i=1}^g \bS^1\,\right]\ .
$$
This gives the simplicial group model $\,\bG(M) \cong {\rm hocolim}\,[\,\bF_g \xleftarrow{\pi} \bG(\Sigma_g) \xrightarrow{\gamma} \bF_g\,]\,$, and hence
$$
\O[\DRep_G(M)]\,\cong\, {\rm hocolim}\,[\,\O(G^g)\,\xleftarrow{\pi_*} \,\cA_\ast(\Sigma_g, G) \,\xrightarrow{\gamma_*}\,
\O(G^g)\,]\,\cong\, \O(G^g) \,\otimes^{\L}_{\cA_\ast}\, \O(G^g)\ ,
$$
where $ \cA_\ast = \cA_\ast(\Sigma_g, G) $ is an explicit DG algebra model for the
representation homology $ \HR_*(\Sigma_g, G) $ (see Section~\ref{S7.1.2}).
As a result, we have the following generalization of Theorem~\ref{lensspaces} to $3$-manifolds of
higher genus.
\bthm
\la{3mnld}
Let $ M $ a closed connected orientable $3$-manifold. Assume that $ M $ has a
Heegaard decomposition \eqref{HeeG}  determined by an element
$ \gamma \in {\mathcal M}(\Sigma_g) $ in the mapping class group of $ \Sigma_g $. Then the representation homology
of $ M $ is given by
$$
\HR_\ast(M, G)\, \cong \, \boldsymbol{\mathrm{Tor}}_\ast^{\cA_*}(\cO(G^g)_\pi,\,\cO(G^g)_\gamma)\ ,
$$
where $ \boldsymbol{\mathrm{Tor}}_\ast^{\cA_*} $ is the differential \,{\rm Tor} taken over the DG algebra $ \cA_* = \cA_*(\Sigma_g, G) $.
In particular, there is an Eilenberg-Moore homology spectral sequence
$$
E^2_{\ast, \ast} = \Tor_*^{\HR_\ast(\Sigma_g, G)} (\O(G^g)_\pi,\,\O(G^g)_\gamma)\ \Longrightarrow\
\HR_\ast(M, G)
$$
converging to the representation homology of $M$.
\ethm
\begin{remark}
If $G$ is a complex {\it semisimple} group and $ g \ge 2 $, it is conjectured in \cite{BRY}
({\it cf.} \cite[Conjecture 1.3]{BRY}) that $\,\HR_i(\Sigma_g, G) = 0 \,$ for all $ i > 0 $.
This conjecture implies, in particular, that the spectral sequence of Theorem~\ref{3mnld}
degenerates for $3$-manifolds of Heegaard genus $ g \ge 2 $, giving an isomorphism
$$
\HR_*(M, G) \cong \Tor_*^{\cA_G(\Sigma_g)} (\O(G^g)_\pi,\,\O(G^g)_\gamma)
$$
where $ \Tor_\ast $ is the ordinary `Tor' taken over $ \cA_G(\Sigma_g) := \O[\Rep_G(\Sigma_g)]\,$,
the coordinate ring of the classical representation scheme $ \Rep_G(\Sigma_g)$.
\end{remark}

\section{Representation cohomology and a non-abelian Dennis trace map}
\la{sec8}
In this section, we define representation homology and cohomology with
coefficients in an arbitrary bifunctor on the category of finitely generated free groups $ \mathfrak{G} $. Following the analogy with topological Hochschild homology, we construct a natural trace map relating representation homology to the stable homology of automorphism
groups $ \Aut(\mathbb F_n) $ with twisted coefficients.

\subsection{Representation cohomology}
\subsubsection{(Co)homology of small categories}
Let $\mathscr C$ be a small category.  By a $\Cc$-{\it bimodule}, we mean a bifunctor $D\,:\,\Cc^{\mathrm{op}} \times \Cc \rar \ve_k$, which is contravariant in the first argument and covariant in the second. We write
$\, \Bimod(\Cc) \,$ for the category of $\Cc$-bimodules.
For any $ D \in \Bimod(\Cc) $, one can define
the (Hochschild-Mitchell) homology $\HH_\ast(\Cc,D)$ and cohomology $\HH^\ast(\Cc,D)$ of $ \Cc $ with coefficients
in $ D $. For a precise definition and basic properties of these classical (co)homology theories we refer to \cite{Mit, BW, GNT} (a good summary can also be found in \cite[Appendix C]{L}). Here, we only recall that $ \HH_\ast(\Cc,\,\mbox{--}\,)$
and $ \HH^\ast(\Cc,\,\mbox{--}\,) $ are functors (covariant and contravariant, respectively)
on the category of $\Cc$-bimodules, such that  $\, \{\HH_n(\Cc,\,\mbox{--}\,)\}_{n \ge 0} \,$ and
$ \{\HH^n(\Cc, \,\mbox{--}\,)\}_{n \ge 0} $ are universal $\delta$-sequences, with $\HH_0(\Cc,\,D)$ and
$\HH^0(\Cc, \,D)$ being canonically isomorphic to the coend $ \int^{c \in \Cc} D(c,c) $ and the end
$ \int_{c \in \Cc} D(c,c) $ of the bifunctor $ D $. Moreover,  the (co)homology theories  $ \HH_\ast(\Cc,\,D)$
and $ \HH^\ast(\Cc,\,D) $ have good functorial properties with respect to the first argument: in particular, any functor
$\, F: \Cc' \rar \Cc \,$ between small categories induces a natural map on homology
$\, F_\ast: \HH_\ast(\Cc', F^\ast D) \rar \HH_\ast(\Cc, D)\,$,
where $ F^*: \Bimod(\Cc) \to \Bimod(\Cc') $ is the restriction functor on bimodules defined by $ F^* D := D \circ (F^{\rm op} \times F) $.

\subsubsection{Representation cohomology}
To express representation homology in terms of Hochschild-Mitchell homology, we need to slightly
extend the above classical setting. Specifically, we will consider chain complexes of $\Cc$-bimodules,
which are simply bifunctors $\,D: \Cc^{\mathrm{op}} \times \Cc \rar {\tt Ch}_{\ge 0}(k) \,$ with values
in the category of chain complexes of $k$-vector spaces, and define
$\, \bHH_\ast(\Cc, D) \,$ and $\, \bHH^\ast(\Cc, D) \,$ to be the Hochschild-Mitchell {\it hyperhomology} and
the Hochschild-Mitchell  {\it hypercohomology} of $D$, respectively. Now, given two chain complexes
of right and left $\Cc$-modules, say $\, M: \Cc^{\rm op} \to  {\tt Ch}_{\ge 0}(k)\,$ and
$\, N: \Cc \to  {\tt Ch}_{\ge 0}(k)\,$, we define the chain complex of $\Cc$-bimodules $\,
M \boxtimes_k N: \Cc^{\mathrm{op}} \times \Cc \to {\tt Ch}_{\ge 0}(k)\,$ by assigning to
$ (c,c') \in  {\rm Ob}(\Cc^{\mathrm{op}} \times \Cc) $ the tensor product $\,M(c) \otimes_k N(c') \,$
of the corresponding chain complexes. With this notation, we have
\blemma
\la{HRHH}
For any $ X \in \sset_0 $ and any commutative Hopf algebra $ {\mathcal H} $, there is a natural
isomorphism
\begin{equation}
\la{hrhh}
\HR_\ast(X, {\mathcal H}) \,\cong\, \bHH_\ast(\mathfrak{G}, \,\uN(k[\lgr{X}]) \boxtimes_k \ucH\,)\ .
\end{equation}
\elemma
\bproof
For any small category $ \Cc $ and any right (resp., left) $\Cc$-modules $\, M\,$ and  $\, N \,$ with values in
$ {\tt Ch}_{\ge 0}(k) $, where $ k $ is a commutative ring, there is a natural (Grothendieck) spectral sequence (see, e.g., \cite[(C.10.1)]{L}):
$$
E^2_{pq} = \bHH_{p}(\Cc,\, \H_q[M \boxtimes_k^{\L} N])\, \Longrightarrow\, \H_{p+q}[M\otimes^{\L}_{\Cc} N]\ .
$$
When $ k $ is a field, this spectral sequence degenerates giving an isomorphism
$\,\bHH_{\ast}(\Cc,\, M \boxtimes_k N) \cong \H_{\ast}[M\otimes^{\L}_{\Cc} N]\,$. In
our situation, we have
$$
\bHH_\ast(\mathfrak{G}, \,\uN(k[\lgr{X}]) \boxtimes_k \ucH\,) \,\cong\, \H_\ast[\uN(k[\lgr{X}])
\otimes^{\L}_{\ffgr} \ucH]\ ,
$$
which in composition with the isomorphism of Theorem~\ref{dtp} gives \eqref{hrhh}.
\eproof
\begin{example}
\la{groups}
In the case when $ X = {\rm B}\Gamma $ for a discrete group $ \Gamma $ and $  {\mathcal H} = \O(G) $,
formula \eqref{hrhh} reads
$$
\HR_\ast(\Gamma, G) \,\cong\, \HH_\ast(\mathfrak{G}, \,k[\Gamma] \boxtimes_k \O(G)\,)\ .
$$
\end{example}

\noindent
Lemma~\ref{HRHH} motivates the following definition.

\begin{definition}
The {\it representation cohomology} of $ X $ in $ \cH $ is defined by
$$
\HR^\ast(X, \cH) \,:=\,\bHH^\ast(\mathfrak{G}, \,\uN(k[\lgr{X}]) \boxtimes_k \ucH\,)
$$
More generally, for any $\ffgr$-bimodule $\,D: \Cc^{\mathrm{op}} \times \Cc \rar {\tt Ch}_{\ge 0}(k) \,$,
we define the {\it representation homology} and the {\it representation cohomology} of $D$ by
$$
\HR_\ast(D) := \bHH_\ast(\ffgr,\,D)\ ,\quad \HR^\ast(D) := \bHH^{\ast}(\ffgr,\,D)\ .
$$
\end{definition}

In the case when $ D $ is an ordinary $\ffgr$-bimodule (with values in $ \ve_k $),
this definition says that the representation (co)homology of $D$ is just the classical Hochschild-Mitchell (co)homology of $ D $.
\begin{example}
\la{cohgr}
For an affine algebraic group $G$, consider the $\ffgr$-bimodule $\,D := \lin_k^* \,\boxtimes \,\O(G)
\,$, where $\, \lin_k^* \,$ is the dual linearization functor
$\, \ffgr^{\rm op} \to \ve_k\,$,$\, \langle n \rangle \mapsto \Hom_{\Z}(\langle n \rangle_{\rm ab},\,k)\,$.
In this case, one can show that there are natural isomorphisms
$$
\HR^i (D) \cong \H^{i+1}(G, k)\ ,\quad \forall\,i> 0 \ ,
$$
where $ \H^{i+1}(G, k) $ stands for the classical cohomology of the affine algebraic group
with coefficients in the trivial (rational) representation.
\end{example}
\vspace{1.5ex}

\subsubsection{Relation to topological Hochschild homology}
\la{RelTHH}
For an arbtirary (associative unital) ring $ R $, denote by $ F(R) $ the full subcategory of $R$-$\Mod$
whose objects are the free modules $ R^n $, $\, n\ge 0$. For any $R$-bimodule $N$, consider
the bifunctor $\,\Hom(I, N):\, F(R)^{\rm op} \times F(R) \to \Mod(\Z)\,$ defined by
$(X,Y) \mapsto \Hom_R(X, N \otimes_R Y)) $. Then, a theorem of Pirashvili and Waldhausen \cite{PW} asserts
that the Hochschild-Mitchell homology $\, \HH_\ast(F(R),\,\Hom(I, N)) \,$
is naturally isomorphic to the  {\it topological Hochschild homology} $ \,{\rm THH}_\ast(R, N) $ of the ring $ R $
with coefficients in the bimodule $ N $. It is therefore natural to {\it define} the
topological Hochschild homology of $R$ with coefficients in an arbitrary bifunctor
$\,B:\, F(R)^{\rm op} \times F(R) \to \Mod(\Z)\,$ by ({\it cf.} \cite[Chap.~13]{L})
$$
{\rm THH}_\ast(R,B) := \HH_\ast(F(R),\,B)\ .
$$
For $ R = \Z $, the category $ F(\Z) $ is equivalent to the category $ \mathfrak{G}_{\mathrm{ab}} $
of finitely generated free abelian groups, which (as our notation suggests) is the
abelianization of the category $ \ffgr $. The abelianization functor $\, \alpha: \ffgr \to \ffgr_{\rm ab} \,$
induces a natural map $\,\HR_\ast(\alpha^* B) \rar \mathrm{THH}_\ast(\Z, B)\,$ for any $ \ffgr_{\rm ab}$-bimodule
$B \in \Bimod(\ffgr_{\rm ab})$, and conversely, for any $\ffgr$-bimodule $ D \in \Bimod(\ffgr) $, associated to the
functor $ \alpha $, there is an Andr\'{e}-type spectral sequence (see \cite[Theorem~1.20]{GNT}):
$$
E_{pq}^2 = {\rm THH}_p(\Z,\,\L_q(\alpha^{\rm op} \times \alpha)_\ast D)\,\Longrightarrow\,\HR_{p+q}(D)\ ,
$$
converging to the representation homology of $D$.

Thus, representation homology may be viewed as a non-abelian analogue of topological Hochschild homology, and it is
natural to ask for `non-abelian' analogues of various constructions known for topological Hochschild homology.
In the next section, we outline one such construction which may be thought of as a non-abelian version
of the Dennis trace map.

\subsection{Non-abelian Dennis trace map}
Recall\footnote{See~\cite[Sect.~13.1.8]{L} for the case $ B = \Hom(I, N) $ and
\cite{FFPS} for an extension to arbitrary $ F(R)$-bimodules.} that the classical Dennis trace maps the stable homology of the general linear groups of a ring $R$ to topological Hochschild homology of $ R \,$:
\begin{equation}
\la{DTRcl}
\mathrm{DTr}_\infty(R,B):\ \H_\ast(\GL_\infty(R),\, B_\infty) \,\rar\, {\rm THH}_\ast(R, B)\ ,
\end{equation}
where $ B $ is an arbitrary bimodule over $ F(R) $.  We generalize this map to the non-abelian setting.

Let $ \Aut_n := \Aut(\bF_n) $ denote the automorphism group of the free group on generators
$ x_1, \ldots, x_n $. We will regard $ \Aut_n $ as the automorphism group $ \Aut_{\ffgr}(\langle n \rangle) $
of the object $ \langle n \rangle $ in the category $ \ffgr $.
There are obvious inclusions $\, \Aut_n \hookrightarrow \Aut_{n+1}\,$ defined
by $\,g \mapsto \tilde{g}\,$, where $\,\tilde{g}(x_i) := g(x_i)$ for $i \leq n $
and $\tilde{g}(x_{n+1}) = x_{n+1} $. We set $\,\Aut_\infty := \varinjlim\, \Aut_n\,$.

Now, consider an arbitrary bimodule $ D $ on the category $ \ffgr $,
i.e. a bifunctor $\,D: \ffgr^{\rm op} \times \ffgr \to \ve_k \,$. For each $ n \geq 1 $, let
$\,D_n := D(\langle n \rangle, \langle n \rangle)$ and define the linear maps
\begin{equation}\la{indmaps}
p^{\ast} \!\circ i_\ast:\ D_n \to D(\langle n \rangle, \langle n+1 \rangle)
\to D_{n+1}\ ,
\end{equation}
where $\, i_\ast := D(\id, i_n) \,$ and $\,p^\ast := D(p_n,\id)\,$ are induced by the natural
inclusion $\, i: \langle n \rangle \into \langle n+1 \rangle\,$ and the natural projection
$\, p: \langle n+1 \rangle \onto \langle n \rangle\,$, respectively. Put
$$
D_\infty:= \varinjlim\,D_n \ ,
$$
where the inductive limit is taken with respect to the linear maps \eqref{indmaps}.

Next, observe that each $D_n$ carries a natural $\Aut_n$-module structure: namely,
$\,\Aut_n \rar \Aut(D_n) $,$\ g \mapsto g^\ast \!\circ g_\ast$,
where $g^\ast:= D(g^{-1},\id)$ and $g_\ast:=D(\id,g)$.
Moreover, for all $\,g\in \Aut_n\,$, there is a commutative diagram
$$
\begin{diagram}[small]
 D_n & \rTo^{g^\ast g_\ast} & D_n\\
  \dTo^{p^{\ast}  i_\ast} & & \dTo_{p^{\ast}  i_\ast}\\
 D_{n+1} & \rTo^{\tilde{g}^\ast  \tilde{g}_\ast} & D_{n+1}\\
 \end{diagram}
 $$
where $ \tilde{g} \in \Aut_{n+1} $ is the image of $ g $ under the natural inclusion $\, \Aut_n \hookrightarrow \Aut_{n+1}\,$ defined above.
As a consequence, the $k$-vector space $\,D_\infty\,$ carries a natural (inductive) $\,\Aut_\infty$-module structure.
Thus, we can consider the homology groups $\H_\ast(\Aut_n, D_n)$ for all $n \geq 1$ and $\H_\ast(\Aut_\infty, D_\infty)$.
Since homology commutes with direct limits, we can identify
\begin{equation}
\la{limi}
\H_\ast(\Aut_\infty, D_\infty)\,\cong\, \varinjlim\, \H_\ast(\Aut_n, D_n)\ .
\end{equation}

Next, we construct natural maps relating $ \H_\ast(\Aut_\infty, D_\infty) $
to the representation homology $ \HR_\ast(D)\,$.
Regarding each group $\,\Aut_n\,$ as the category  $\,\underline{\Aut}_n\,$ with a single object,
we consider the inclusion functors
$$
\gamma_n:\ \underline{\Aut}_n \,\to\, \mathfrak{G}\ ,
$$
identifying the single object of $ \underline{\Aut}_n $ with  $ \langle n \rangle \in \mbox{Ob}(\ffgr) $.
Since $\, D_n = \gamma_n^\ast D \,$, these functors induce natural maps
\begin{equation}
\la{HHgn}
(\gamma_n)_\ast:\, \HH_\ast(\underline{\Aut}_n,\,D_n) \to \HH_\ast(\mathfrak{G},\,D) =: \HR_\ast(D)\ .
\end{equation}
On the other hand, the Hochschild-Mitchell homology of the category
$ \underline{\Aut}_n $ coincides with the usual group homology of $ \Aut_n\,$:
\begin{equation}
\la{HHiso}
\HH_\ast(\underline{\Aut}_n,\,D_n) \cong \H_\ast(\Aut_n,\,D_n)\ .
\end{equation}
Indeed, since $ \underline{\Aut}_n $ is a category with one object, its
Hochschild-Mitchell complex $\,C^{\rm HM}_\ast(\underline{\Aut}_n,\,D_n)\,$
is isomorphic to the usual Hochschild complex $\,C_\ast(k[\Aut_n],\,D_n)\,$
of the group algebra of $ \Aut_n $, so that
$$
\HH_\ast(\underline{\Aut}_n,\,D_n) \cong \HH_\ast(k[\Aut_n],\,D_n) \ ,
$$
while $\,\HH_\ast(k[\Aut_n],\,D_n) \cong \H_\ast(\Aut_n,\,D_n)\,$
via the classical Mac Lane Isomorphism (see, e.g., \cite[Prop.~7.4.2]{L}).
Thus, combining \eqref{HHgn} and \eqref{HHiso}, for all $ n \ge 0\,$, we get canonical linear maps
\begin{equation}
\la{DTRn}
\mathrm{DTr}^{\mathfrak{G}}_n(D):\ \H_\ast(\Aut_n, D_n) \rar \HR_\ast(D) \ .
\end{equation}
As in \cite[13.1.8]{L}, it is easy to check that these maps are compatible
when passing from $n$ to $n+1$. Hence, we can stabilize \eqref{DTRn} by passing to the inductive limit as
$ n \to \infty\,$. With identification \eqref{limi}, the resulting stable map reads
\begin{equation}
\la{DTRinf}
\mathrm{DTr}^{\mathfrak{G}}_\infty(D):\ \H_\ast(\Aut_\infty, D_\infty) \,\rar \,\HR_\ast(D)\ .
\end{equation}
This is a non-abelian analogue of the  Dennis trace map \eqref{DTRcl}. As in the classical case,
it is natural to ask: When is \eqref{DTRinf} an isomorphism?
Motivated by a theorem of Scorichenko (see \cite{FFPS}), we propose a conjectural answer:

\bconj
\la{conj3}
{\it The map \eqref{DTRinf} is an isomorphism if $ D $ is a polynomial\footnote{Strictly speaking,
the notion of a polynomial functor is usually defined in the literature
for functors $ T: {\mathscr C} \to {\mathscr A}\, $ between two additive categories
(see, e.g., \cite[E.13.1.3]{L}). However, the definition of the $n$-th cross-effect, in terms of which
one defines polynomial degrees, makes sense for functors $T$ whose domain $ {\mathscr C} $ is
any pointed cocartesian monoidal category; in particular, it applies to the category $ \ffgr $ (see \cite{HPV}).} bifunctor.}
\econj

\vspace{1ex}

\noindent
We conclude this section a few remarks related to Conjecture~\ref{conj3}.

\vspace{1ex}

\begin{remark} 
A famous theorem of Galatius \cite{Ga} asserts that natural maps from the symmetric group $ S_n $
to $ \Aut_n $ (defined by permuting the generators) induce isomorphisms
$$
\H_i(\Aut_n,\,\Z) \,\cong\,\H_i(S_n,\, \Z)\ ,\quad \forall\,  n > 2i+1 \ .
$$
This implies that $\,\H_i(\Aut_\infty, \,A) = 0 \,$ for all $ i > 0 $, where $A$ is any constant
$k$-module provided $k$ has characteristic $0$ (which we always assume in this paper).
Conjecture~\ref{conj3} implies \cite[Theorem 1]{DV}, which says that
$\,\H_i(\Aut_\infty, \,D_\infty) = 0 \,$ for all $ i > 0 $, when $ D $ is a polynomial bifunctor,
{\it constant} with respect to its contravariant argument. Indeed, for such bifunctors, we have
$ \HR_i(D) = \HH_i(\ffgr, D) = 0 $ for $ i > 0 $ because $ \ffgr $ has a terminal object.
\end{remark}

\vspace{1ex}

\begin{remark}
 The direct analogue of Conjecture~\ref{conj3} is false in the abelian case. Indeed, if $ B $
is a constant bifunctor on $ F(R) $, then $ {\rm THH}_\ast(R, B) $ vanishes in positive degrees
(since $ F(R) $ has terminal object), but $ \H_\ast(\GL_\infty(R),\, B) $ may be highly
nontrivial (see \cite{FFPS}). The correct version of Conjecture~\ref{conj3}
replaces the stable group homology with Waldhausen's stable $K$-theory. In the non-abelian case, one can also state
a version of Conjecture~\ref{conj3} for the stable $K$-theory of automorphism groups $ \Aut_n $ instead of
group homology; however, we expect that the two theories are actually isomorphic.
We briefly outline an argument behind this expectation.

Let $E_\infty$ denote the commutator subgroup of $\Aut_\infty$. It is known that $E_\infty$ is a perfect normal
subgroup, hence we can form the `plus construction'
$$
\Psi:\ B \Aut_\infty \rar B\Aut_\infty^+ \ .
$$
Let $F\Psi$ denote the homotopy fiber of the map $\Psi$. We have a canonical group homomorphism
$\,\pi_1(F\Psi) \rar \pi_1(B\Aut_\infty)\cong \Aut_\infty$ that equips any $\Aut_\infty$-module with
a $\pi_1(F\Psi)$-action. In particular, the $\Aut_\infty$-module $D_\infty$ arising from a
$ \ffgr$-bimodule $ D $ may be viewed as a $\pi_1(F\Psi)$-module, and hence defines a local system on $F\Psi$.
The {\it stable $K$-theory} $\, K_\ast^{s}(\Aut_\infty, D_\infty) \,$ is then
defined to be $\, \H_\ast(F\Psi,\, D_\infty)\,$, the homology of $F\Psi$ with coefficients
in the local system $D_\infty$.
Now, consider the Serre spectral sequence associated to the homotopy fibration
$\,F\Psi \rar B\Aut_\infty \rar B\Aut_\infty^+\,$:
$$
E^2_{pq}\,=\,\H_p(B\Aut_\infty^+, \,\H_q(F\Psi, D_\infty)) \,\implies\, \H_n(B\Aut_\infty,\, D_\infty)\ .
$$
If $\Aut_\infty$ acts {\it trivially} on $ K_q^{s}(\Aut_\infty,\, D_\infty) = \H_q(F\Psi, \, D_\infty)$
(as it happens in the classical case, see \cite[13.3.2]{L}), then, since
$\, B\Aut_\infty \rar B\Aut_\infty^+\,$ is a homology equivalence for trivial coefficients,
the above spectral sequence becomes
$$
E^2_{pq}\,=\,\H_p(B\Aut_\infty,\, K_q^s(\Aut_\infty,\, D_\infty)) \,\implies\, \H_n(B\Aut_\infty,\, D_\infty)\ .$$
However, by Galatius' Theorem \cite{Ga}, we know that $\H_p(\Aut_\infty, A)=0$ for $p >0$ for any constant coefficients
over $k$. Hence, the above spectral sequence must collapse on the $p$-axis, giving the desired isomorphism $K_{\ast}^{s}(\Aut_\infty, D_\infty)\,\cong\, \H_\ast(\Aut_\infty,D_\infty)$.
\end{remark}

\vspace{1ex}

\begin{remark}A
s explained in Section~\ref{RelTHH}, the relation between topological Hochschild homology and functor homology of module categories is based on the Pirashvili-Waldhausen Theorem \cite{PW}.
Schwede \cite{S} generalized this result to arbitrary algebraic theories by associating to an algebraic PROP $\, \mathfrak{P} $ a ring spectrum $  \mathfrak{P}^s $ and identifying
the functor homology over $ \mathfrak{P} $ with topological Hochschild homology over $  \mathfrak{P}^s $ (see \cite[Theorem~6.7]{S}). In the case $ \mathfrak{P} = \mathfrak{G} $, Schwede's construction
provides a topological (spectral) interpretation of representation homology which may be useful for Conjecture~\ref{conj3}.
\end{remark}

\appendix

\section{Model approximations and derived adjunctions}
\la{App}
In this Appendix, we collect basic definitions and prove some results in abstract homotopy theory concerning derived functors. We  work in the framework of homotopical categories in the sense of Dwyer, Hirschhorn, Kan and Smith \cite{DHKS}. Apart from the original reference \cite{DHKS}, a good introduction to the subject can be found in \cite{R} and a short summary in \cite{Shul}. The
main results of this Appendix --- Theorem~\ref{ThA2} and Theorem~\ref{ThA3} ---
arise from our attempt to abstract Theorem~\ref{ThAdj} on derived representation adjunctions.
We believe that these two theorems as well as Lemma~\ref{LeA1} are of independent interest.

\subsection{Homotopical categories}
\la{App1}
A {\it homotopical category}\, is a category $\C$ equipped with a class of morphisms $ \W $ (called weak equivalences) that contains all identities
of $ \C $ and satisfies the following  {\it  2-of-6 property}: for every composable triple of morphisms $ f,g,h \in {\rm Mor}(\C)$, if $ g  f $ and $    h g $ are in $ \W $, then so are $ f, g, h $ and $h g f $.  The 2-of-6 property formally implies, but is stronger than, the usual 2-of-3 property.
The class of weak equivalences thus forms a subcategory which contains all objects and all isomorphisms of $ \C $. Since the isomorphisms satisfy the 2-of-6 property, any category can be viewed as a homotopical category by taking $ \W $ to be the class of all isomorphisms (in \cite{DHKS}, such homotopical categories called minimal). Furthermore, by forgetting the fibrations and cofibrations, any model category becomes a homotopical category: i.e., the class of weak equivalences in any model category satisfies the  2-of-6 property  (see \cite[Prop. 9.2]{DHKS}). This is a consequence of the well-known fact that in a model category, the class $ \W $ of weak equivalences is {\it saturated}: i.e., it comprises all the arrows of $ \C$ that become isomorphisms in the localized category $ \C[ \W^{-1}] \,$. Since the isomorphisms
in $ \C[\W^{-1}] \,$ satisfy the 2-of-6 property, it follows immediately that the weak equivalences in a saturated category satisfy the 2-of-6 property. Unless stated otherwise, we will assume all our homotopical categories to be saturated. If $ \C $ is a homotopical category, the category $\, \Ho(\C) := \C[ \W^{-1}] \,$ is called the {\it homotopy category} of $ \C $: it comes with the canonical functor
$ \gamma_{\C}: \C \to \Ho(\C) $ called the localization of $ \C $. It is often convenient to regard $ \Ho(\C) $ as a homotopical category itself by taking $ \W $ to be the class of isomorphisms: in other words, to think of $ \Ho(\C)$ as a minimal homotopical category.

\subsection{Derived functors and deformation retracts}
\la{App2}
If $ \C$ and $ \D $ are homotopical categories, a functor $ F: \C \to \D $  is called {\it homotopical} if it preserves weak equivalences. Such a functor induces a unique functor between the homotopy categories
of $ \C $ and $ \D $ which we will denote by  $ \bar{F}: \Ho(\C) \to \Ho(\D) $. In practice,
many important functors are not homotopical and hence do not descend to  homotopy categories.
A standard way to deal with this problem is to replace --- or `approximate' --- non-homotopical functors with their derived functors which usually come in two kinds: `left' and `right'.
We will focus on left derived functors  with understanding that all results 
apply {\it mutatis mutandis} to the right derived functors as well.

Following \cite{Q1}, we define a {\it total left derived functor} $ \L F: \Ho(\C) \to \Ho(\D) $
of a functor $ F: \C \to \D $ to be the right Kan extension of $\,\gamma_{\D} \circ F: \,\C \to \D \to \Ho(\D) \,$ along localization $ \gamma_{\C}: \C \to \Ho(\C)\, $:
$$
\begin{diagram}[small]
\C                 &  \rTo^{F}         & \D \\
\dTo^{\gamma_{\C}} &  \Uparrow                &  \dTo_{\gamma_{\D}}\\
\Ho(\C)            & \rDashto_{\L F}  &\Ho(\D)                \\
\end{diagram}
$$
By the universal property of localization, $\L F $ is uniquely determined by the
homotopical functor $ \bL F = \L F \circ \gamma_{\C} $ defined on the category $ \C $.
This last functor can be characterized as a universal homotopical functor
$ \bL F: \, \C \to \Ho(\D) $ which comes together with a natural transformation (called the
comparison map) $ \varepsilon: \bL F \to \gamma_{\D} \circ F $ that is terminal among
all  natural transformations from homotopical functors to $ \gamma_{\D} \circ F $. When they
exist, both functors $ \L F$ and $ \bL F $ are determined by $F$ uniquely up to unique isomorphism. Following \cite{Shul},
we will refer to $ \bL F $ as a {\it left derived functor} of $ F $ and $ \L F $ as the corresponding total left derived functor.

%
%

It was observed in \cite{Mal} that a stronger universal property for derived functors ---
namely, that of an absolute Kan extension --- is often very useful\footnote{In the additive setting,
absolute derived functors between triangulated categories first appeared in the work of P.
Deligne  under the name ``founcteurs d\'eriv\'e partout d\'efini'' (see \cite{D}). }. To be precise,
a total left derived functor $ \L F: \Ho(\C) \to \Ho(\D) $ is called {\it absolute} if for any
functor $ H: \Ho(\D) \to {\mathcal E} $, the right Kan extension of the composition
$ H \circ \gamma_{\D} \circ F: \C \to \D \to \Ho(\D) \to {\mathcal E} $
along $ \gamma_{\C}: \C \to \Ho(\C) $ coincides with $H \circ \L F $:
$$
\begin{diagram}[small]
\C                 &  \rTo^{F}         & \D \\
\dTo^{\gamma_{\C}} &  \Uparrow                &  \dTo_{\gamma_{\D}}\\
\Ho(\C)            & \rTo_{\L F}  &\Ho(\D)                \\
                   &  \rdTo_{\ {\rm Ran}_{\gamma_{\C}}(H \gamma_{\D} F)\ }       & \dTo_{H}\\
                   &                 & \E
\end{diagram}
$$

A fundamental theorem of \cite{Q1} asserts that any left Quillen functor $F: \C \to \D $
between model categories has a total left derived functor $ \L F: \Ho(\C) \to \Ho(\D) $, which
can be obtained as the composition $ F \circ Q $, where $ Q $ is the cofibrant replacement
functor on $\C$; moreover, as noticed in \cite{Mal}, such a left derived functor
is automatically absolute. This construction of derived functors was axiomatized and extended
to homotopical categories in \cite{DHKS}. We briefly recall the main definitions.
If $ \C $ is a homotopical category, a {\it left deformation retract} of $ \C $ is
a full subcategory $ i: \C_{Q} \into \C $ given together with a homotopical functor
$ Q: \C \to \C_{Q} $ and natural weak equivalence $q: \, i \circ Q \to \id_{\C} $.
It is easy to see that, for any left deformation retract of $ \C $, the inclusion
functor $ i: \C_{Q} \into \C $ induces an equivalence of categories
$ \Ho(\C_{Q}) \simeq \Ho(\C) $ with inverse induced by $ Q $. Now, we say that a functor
$ F: \C \to \D $ between two homotopical categories is {\it left deformable} if there is
a left deformation retract $ \C_Q $ of the domain category such that the restriction of $F$
to $ \C_Q $ is homotopical. For example, if $ \C $ and $ \D $ are model categories,
any left Quillen functor $ F: \C \to \D $ is canonically left deformable: for the corresponding
deformation retract $ \C_Q $, we can always take the subcategory of cofibrant objects in
$ \C $, with $Q: \C \to \C_Q $ being the cofibrant replacement functor.
\bprop[\cite{DHKS}]
\la{leftdef}
A left deformable functor $F: \C \to \D $ has a left derived functor
$\bL F: \C \to \Ho(\D) $ given by $\, \bL F = \gamma_\D \circ F \circ Q
\,$ with comparison map $ \varepsilon = \gamma F q: \bL F \to \gamma \circ F $.
The corresponding total left derived functor $ \L F: \Ho(\C) \to \Ho(\D) $
is absolute in the sense of \cite{Mal}.
\eprop
The first statement of Proposition~\ref{leftdef} is proved in
Sections~41.2-5 of \cite{DHKS} (see also \cite[Theorem~2.2.8]{R}).
The second statement is verified in (the proof of) Proposition~2.2.13
of \cite{R}.

It is well-known that, for any composable pair $(F_1, F_2)$ of left Quillen functors, the derived
functor $ \L(F_1 \circ F_2)$ of their composition coincides with $ \L F_1 \circ \L F_2 $.
In the more general context of homotopical categories, this is not the case even when both functors
$F_1$ and $ F_2 $ are left deformable. To guarantee this property one needs to impose an extra
condition on deformation retracts of the functors involved. Following \cite{DHKS}, we say that
a composable pair $(F_1, F_2) $ of left deformable functors
$ \C \xrightarrow{F_1} \D \xrightarrow{F_2} \E $ is {\it left deformable} if $ F_1 $ maps
the left deformation retract $ \C_Q $, on which it is homotopical, into the left
deformation retract $ \D_Q  $, on which $ F_2 $ is homotopical: i.e,,
$ F_1(\C_Q) \subseteq \D_Q $. With this definition, we have
\bprop[\cite{DHKS}, 42.4]
\la{comdef}
For any left deformable pair $ (F_1, F_2) $, there is a canonical isomorphism of total
left derived functors $\,\L(F_1 \circ F_2) \cong \L F_1 \circ \L F_2\,$.
\eprop

\subsection{Derived adjunctions}
\la{App4}
We now turn to the important question when an adjunction between two homotopical categories
induces a derived adjunction between the corresponding homotopy categories. We begin by
stating the main result of \cite{Mal} ({\it cf.} \cite[2.2.15]{R}):
\bthm[\cite{Mal}]
\la{dermal}
Let $ F:\,\C \rightleftarrows \D \,: G $ be a pair of adjoint functors
between homotopical categories. Assume that $ F $ has a total left derived functor
$ \L F $, $G$ has a total right derived functor $ \bR G $, and both derived
functors are absolute. Then $ \L F $ and $ \bR G $ are adjoint to each other:
\begin{equation}
\la{LFRG}
\L F:\,\Ho(\C) \rightleftarrows \Ho(\D) \,:\bR G
\end{equation}
\ethm
Following \cite{DHKS}, let us call an adjunction $ F:\,\C \rightleftarrows \D \,: G $
{\it deformable} if $F$ is left deformable and $G$ is right deformable. As an immediate
consequence of Theorem~\ref{dermal} and Proposition~\ref{leftdef}, we get
\bcor[\cite{DHKS}, 44.2]
\la{cordef}
If $ F:\,\C \rightleftarrows \D \,: G $ is a deformable adjunction,
then both total derived functors $ \L F $ and $ \bR G $ exist and form an adjoint pair \eqref{LFRG}.
\ecor
This result is one of the key observations of \cite{DHKS}, which, in particular,
formally implies Quillen's Adjunction Theorem for model categories \cite{Q1}. Unfortunately,
the assumption that a pair of adjoint functors is deformable is rather restrictive and does not always hold in practice.
In what follows we propose a different --- somewhat roundabout --- way to produce derived adjunctions using model approximations of homotopical categories.

We begin with the following simple lemma which can be viewed as a partial converse of Theorem~\ref{dermal}
\blemma
\la{LeA1}
Let $ F:\,\C \rightleftarrows \D \,: G $ be a pair of adjoint functors
between homotopical categories. Assume:

$(1)$ $\,F $ has an absolute total left derived functor $ \L F: \Ho(\C) \to \Ho(\D) $,

$(2)$ $\, \L F $ has a right adjoint functor $ \tilde{G}: \Ho(\D) \to \Ho(\C) $.

\vspace{1ex}

\noindent
Then $  	\tilde{G} $ is an absolute total right derived functor of $ G $: i.e., $\bR G $ exists and
$ \bR G \cong 	\tilde{G} $.
\elemma

\bproof
Let us spell out the universal mapping property of the absolute total right derived functor $ \bR G\,$: for any functors $E: \Ho(\C) \to \E $ and $H: \Ho(\D) \to \E $,
%
%
there is a natural (in $E$ and $H$)  bijection:
\begin{equation}
\la{hom1}
\Hom(E \circ \bR G, \, H) \,\cong\, \Hom(E \circ \gamma_{\C} \circ G,\, H \circ \gamma_{\D})\ ,
\end{equation}
where by $ \Hom$'s we denote the sets of natural transformations between the corresponding functors.
To prove the lemma it suffices to check that $ \tilde{G} $ satisfies this property.

First,
$ \tilde{G} $ being right adjoint to $ \L F $ implies that $ \tilde{G}^* = (\,\mbox{--}\,) \circ \tilde{G} $ is left adjoint to $ \L F^* =  (\,\mbox{--}\,) \circ \L F $ on the functor category
$ {\rm Fun}(\Ho(\C), \, \E)\,$, so that there is a natural bijection
\begin{equation}
\la{hom2}
\Hom(E \circ \tilde{G}, \, H) \,\cong\, \Hom(E,\, H \circ \L F)\ .
\end{equation}
Second, the universal mapping property of $ \L F $ being an absolute left derived functor of $F$ gives
\begin{equation}
\la{hom3}
\Hom(E, \, H \circ \L F) \,\cong\, \Hom(E \circ \gamma_{\C},\, H \circ \gamma_{\D} \circ F)\ .
\end{equation}
Third, $F$ being left adjoint to $G$ implies that $ F^* =  (\,\mbox{--}\,) \circ F $ is right adjoint
to $ G^*  =  (\,\mbox{--}\,) \circ G $, hence
\begin{equation}
\la{hom4}
\Hom(E \circ \gamma_{\C}, \, H \circ \gamma_{\D} \circ F) \,\cong\, \Hom(E \circ \gamma_{\C} \circ G,\,
H \circ \gamma_{\D})\ .
\end{equation}
Combining now \eqref{hom2}, \eqref{hom3}, \eqref{hom4} and comparing the result with \eqref{hom1}, we see that $ \tilde{G} $ satisfies the same universal mapping property as $ \bR G $. Whence
$ \bR G = \tilde{G} $.
\eproof
\begin{remark}
There is a dual version of Lemma~\ref{LeA1}: if the absolute total right
derived functor $ \bR G $ for the right adjoint in the pair
$ F:\,\C \rightleftarrows \D \,: G $ exists and
has a left adjoint $ \tilde{F}: \Ho(\C) \to \Ho(\D) $, then this left adjoint $ \tilde{F} $
is the absolute total left derived functor of $F$.
\end{remark}

\subsection{Model approximations}
Next we recall the notion of a model approximation introduced in \cite{CS}.
This notion plays an important role in abstract homotopy theory allowing one
to define homotopy colimits of arbitrary diagrams in model categories. We will use it,
however, for a different purpose: to construct derived adjunctions between homotopical categories.

\begin{definition}[\cite{CS}]
\la{DefA1}
A {\it left model approximation} of a homotopical category $ \C $ is a model
category $ \M $ given together with a pair of adjoint functors $\, l: \M \rightleftarrows \C: r \,$ such that

$(1)$ $r$ is homotopical, i.e. $ r(\W_{\C}) \subseteq \W_{\M} $;

$(2)$ $l$ is homotopical on cofibrant objects of $ \M $;

$(3)$ $(l,r) $ is an `almost Quillen equivalence' in the sense:
for any $ A \in {\rm Ob}(\C) $ and any cofibrant $ X \in {\rm Ob}(\M)  \,$, if
$ f: X \to r(A) $ is a weak equivalence in $ \M $ then the adjoint map
$ f^{\#}: l(X) \to A $ is a weak equivalence in $ \C $.
\end{definition}

The intuition behind this definition is that --- from the homotopy-theoretical point of view ---
being a model category or having a model approximation should not make much difference. Our
Theorem~\ref{ThA2} below illustrates this principle in the case of derived adjunctions.

We will need one more definition ({\it cf.} \cite[Def. 5.8]{CS}). If $ F: \C \to \D $ is a
functor between homotopical categories, we say that a left model approximation
$\, l: \M \rightleftarrows \C: r \,$ is {\it good} for $F$ if the restriction
$ F \circ l: \M \to \C \to \D $ is homotopical on cofibrant objects of $ \M $.
In this case, it follows from property $(3)$ of Definition~\ref{DefA1} that
$\, Q_{\C} := l \circ Q \circ r:\, \C \to \C \,$ provides a left deformation
for $ F $, where $ Q $ is the cofibrant replacement functor on $ \M $.
Thus, if $F$ admits a good left model approximation, then $ F $ is a left deformable functor
and hence, by Proposition~\ref{leftdef}, has an absolute total left derived functor
$ \L F: \Ho(\C) \to \Ho(\D) $. This applies, in particular, to the functor $ l: \M \to \C $
itself (since we can take the identity adjunction on $ \M $ as a good model
approximation for $l$). Now, since $l$ is left deformable and  $r$
is homotopical,  by Corollary~\ref{cordef}, the adjunction  $\, l: \M \rightleftarrows \C: r \,$
induces the adjunction of derived functors
\begin{equation}
\la{Llradj}
\L l:\, \Ho(\M) \rightleftarrows \Ho(\C)\,: \bar{r}
\end{equation}
%


The next lemma clarifies the properties of the derived functors \eqref{Llradj}; it is
essentially a reformulation of \cite[Proposition~5.5]{CS}.
\blemma
\la{LeA2}
Let $\, l: \M \rightleftarrows \C: r \,$ be a left model approximation of
a homotopical category $ \C $. The functor $\, \bar{r}: \Ho(\C) \to \Ho(\M) \,$
induced by $r$ is fully faithful, and the counit morphism $\,
\L l \circ \bar{r}  \xrightarrow{\!\sim}  \id_{\Ho(\C)} $ associated with
\eqref{Llradj} is an isomorphism.
\elemma

\bproof To simplify the notation we write $ \bar{X} \in \Ho(\C) $ for the image
of $ X \in \Ob(\C) $ under the localization functor $ \gamma: \C \to \Ho(\C)$, and
similarly for $ \M $. We need to prove that, for any  $ X, Y \in \Ob(\C) $, the map
\begin{equation*}
\bar{r}_{X,Y}:\,\Hom_{\Ho(\C)}(\bar{X},\,\bar{Y})\,\to\,
\Hom_{\Ho(\M)}(\overline{r(X)},\,\overline{r(Y)})
\end{equation*}
is bijective. For this, we will explicitly construct the inverse map.

Let $ Q,  R: \M \to \M $ denote the cofibrant and the fibrant
replacement functors in $\M$, respectively. Since $ \M $ is a model category,
any morphism $ \bar{f}: \overline{r(X)} \to \overline{r(Y)} $ in $ \Ho(\M) $ can be
represented by a morphism $ f: Q r(X) \to  R Q r(Y) $ in $ \M $. Moreover, we have
the following natural diagram in $ \C $:
\begin{equation}
\la{diagf}
X \xleftarrow{\sim} l Q r(X) \xrightarrow{l(f)} l R Q r (Y) \xleftarrow{\sim} l Q r(Y)
\xrightarrow{\sim} Y
\end{equation}
The first and the last maps in \eqref{diagf} are the adjoints of the cofibrant
resolutions $ Qr(X) \sonto r(X) $ and $ Qr(Y) \sonto r(Y) $ in $ \M $, hence,
by property $(3)$ of Definition~\ref{DefA1}, they are weak equivalences in $ \C $.
The third map is obtained by applying the functor $ l $ to the fibrant
resolution $ R Q r(Y) \xrightarrow{\sim} Q r(Y)  $ of the (cofibrant) object
$ Q r(Y) $ in $ \M $, hence it is also a weak equivalence, by Definition~\ref{DefA1}$(2)$.
Now, applying the localization functor $ \gamma: \C \to \Ho(\C) $ transforms
the weak equivalences in \eqref{diagf} into isomorphisms, and  by inverting these isomorphisms,
we can define a (unique) morphism $ \bar{\psi}_{X,Y}(\bar{f}): \bar{X} \to \bar{Y} $ in $ \Ho(\C) $,
which depends only on $ \bar{f} $. It is straightforward to check that the map
$ \bar{\psi}_{X,Y} $ given by this construction is inverse to $ \bar{r}_{X,Y} $.
This proves the first claim of the lemma. The second claim is equivalent to the first
by abstract properties of adjunctions (see, e.g., \cite[Prop. I.1.3]{GZ}).
\eproof
%


We are now in position to state the main result of this Appendix.
\begin{theorem}
\la{ThA2}
Let $ F:\,\C \rightleftarrows \D \,: G $ be a pair of adjoint functors between homotopical categories.
Assume that $ \C$ admits a left model approximation $\,l: \M \rightleftarrows \C: r\,$  together with
adjoint functors $ \hat{F}: \M \rightleftarrows \D \,: \hat{G} $, such that

$\,(i)\,$ $ (\hat{F}, \hat{G}) $ is a deformable adjunction,

$\,(ii)\,$ $ (\hat{F}, r)$ is a left deformable pair, and there is a natural weak equivalence  $\, \hat{F} \circ r \xrightarrow{\sim} F\,$.

$\,(iii)\,$ $ \im(\bR \hat{G}) \subseteq \im(\bar{r})$.

\vspace*{1ex}

\noindent
Then $ F $ and $G $ have total (left and right) derived functors given by 
\begin{equation}
\la{LlRG}
\L F = \L \hat{F} \circ \bar{r}\ , \quad
\bR G = \L l \circ \bR \hat{G}\ .
\end{equation}
The derived functors $\L F $ and $ \bR G $ are both absolute and  adjoint to each other:
$$
\L F:\,\Ho(\C) \rightleftarrows \Ho(\D) \,:\bR G
$$
\end{theorem}
\bproof
First, note that, by $ (i)$ and Proposition~\ref{leftdef}, the derived functors $\, \bL \hat{F} $ and $ \R \hat{G} $ exist, and by Corollary~\ref{cordef}, the corresponding total derived functors $ \L \hat{F} $ and $ \bR \hat{G} $ are adjoint to each other.
By $ (ii)$, the functor $ \bL F := \bL (\hat{F} \circ r) $ satisfies the universal property of a left derived functor of $F$ provided we define the comparison map $ \varepsilon: \bL F \to \gamma \circ F $ to be the composition $ \varepsilon := \gamma(\varphi) \circ \hat{\varepsilon} $, where $ \varphi: \hat{F} \circ r \to F $ is the natural weak equivalence of $(ii)$ and $ \hat{\varepsilon} $ is the comparison map for the derived functor $  \bL (\hat{F} \circ r) $. By Proposition~\ref{comdef}, the total left derived functor $  \L (\hat{F} \circ r) $ is absolute and isomorphic to  $ \L \hat{F} \circ \bar{r} $. Hence $ \L F = \L \hat{F} \circ \bar{r}  $ is an absolute total left derived functor of $F$.

Now, by $(iii)$, we can factor $ \bR \hat{G} $ as
a composition: $\,\Ho(\D) \xrightarrow{\bar{G}_0} \bar{\C} \stackrel{\bar{i}}{\into} \Ho(\M) \,$, where $ \bar{\C} := \im(\bar{r}) $ denotes the essential image of  $ \bar{r} $ in $ \Ho(\M) $ and
$ \bar{i} $ is the inclusion functor. By Lemma~\ref{LeA2}, we can also factor $ \bar{r} = \bar{i} \circ \bar{r}_0 $, where
$\, \bar{r}_0: \Ho(\C) \stackrel{\sim}{\to} \bar{\C} \,$ is an equivalence, with
quasi-inverse $\, \bar{l}_0 := \L l \circ \bar{i}:\, \bar{\C} \to \Ho(\C) \,$. Combining these two factorizations, we can write $ \L l \circ \bR \hat{G} = \L l \circ \bar{i} \circ \bar{G}_0 =
\bar{l}_0 \circ \bar{G}_0 $. Then, for any objects $ X \in \Ho(\C) $ and  $ A \in \Ho(\D) $, we have
\begin{eqnarray*}
\Hom_{\Ho(\C)}(X, \,(\L l \circ \bR \hat{G})(A)) & = & \Hom_{\Ho(\C)}(X, \,\bar{l}_0(\bar{G}_0(A))) \\
& \cong & \Hom_{\bar{\C}}( \bar{r}_0(X), \, \bar{G}_0(A))\\
& \cong & \Hom_{\Ho(\M)}( \bar{i}(\bar{r}_0(X)), \, \bar{i}(\bar{G}_0(A)))\\
& \cong & \Hom_{\Ho(\M)}(\bar{r}(X), \, \bR \hat{G}(A))\\
& \cong & \Hom_{\Ho(\D)}(\L\hat{F}(\bar{r}(X)), \, A) \\
& \cong & \Hom_{\Ho(\D)}(\L F(X), \, A) \ .
\end{eqnarray*}
This shows that $ \L l \circ \bR \hat{G} $ is right adjoint to $\L F$, which is an absolute left derived functor. Hence, by Lemma~\ref{LeA1}, we conclude that $ \bR G $ exists and $ \bR G \cong  \L l \circ \bR \hat{G} $.
\eproof
\begin{remark}
\la{R1}
{\bf 1.} Under the assumptions of Theorem~\ref{ThA2}, there is
a natural isomorphism of functors
\begin{equation}
\la{Lleft}
\L \hat{F} \cong \L F \circ \L l \ ,
\end{equation}
which is {\it a priori}\, a stronger condition than $ \L F \cong \L \hat{F} \circ \bar{r} $.
Indeed, by Theorem~\ref{ThA2}, the functor $ \L F $ has a right adjoint $ \bR G $ which can be written, using
the notation introduced in the proof, as $ \bR G = \bar{l}_0 \circ \bar{G}_0 $.
Since $ \bar{l}_0 $ is an equivalence with quasi-inverse $ \bar{r}_0 $, this implies
$$
\bar{r} \circ \bR G = \bar{r} \circ \bar{l}_0 \circ \bar{G}_0 = \bar{i} \circ
\bar{r}_0 \circ \bar{l}_0 \circ \bar{G}_0 \cong \bar{i} \circ \bar{G}_0 = \bR \hat{G}
$$
Thus, we have an isomorphism of functors $ \bR \hat{G} \cong \bar{r} \circ \bR G $,
where each functor has a left adjoint. By adjunction, this gives \eqref{Lleft}.

{\bf 2.}
The main assumption of Theorem~\ref{ThA2} --- namely, the condition that the adjunction
$ \hat{F}: \M \rightleftarrows \D \,: \hat{G} $ is defined on the whole model category $ \M $ ---
can be weakened. The  proof shows that it suffices
to assume that $ \hat{F} $ exists on a full subcategory $ \M' $ of $ \M $, which is closed under the weak equivalences in $ \M $ and whose image in $ \Ho(\M) $ contains $ \im(\bar{r}) $.
\end{remark}

\subsection{Homotopy colimits}
Recall that any adjunction $ F: \C \rightleftarrows \D : G $ extends formally to an adjunction
$ F^I: \C^I \rightleftarrows \D^I : G^I $ of the diagram categories
$ \C^I := {\rm Fun}(I,\C) $ and $ \D^I := {\rm Fun}(I,\D) $
for any small category $I$. The corresponding functors $ F^I $ and $ G^I $ are given by compositions $ F^I(X) = F \circ X $ and $ G^I(Y) = G \circ Y $, where $ X \in \Ob(\C^I) $ and $Y \in \Ob(\D^I) $.
If $ C $ is a homotopical category, the diagram category $ \C^I $ has a natural homotopical
structure in which a morphism of $I$-diagrams $\, \varphi: X \to X' \,$ is a weak equivalence
if $ \varphi_i: X(i) \xrightarrow{\sim} X'(i) $ is a weak equivalence in $\C$ for every object $ i \in \Ob(I) $.
Moreover, as observed in \cite{DHKS}, if the functor $F: \C \to \D $ is  left deformable,
then so is $ F^I: \C^I \to \D^I $: in fact, if $ Q: \C \to \C_Q $ is a left deformation
retract for $ F $, then $ Q^I: \C^I \to \C_Q^I $ is a left deformation retract for $ F^I$. By Proposition~\ref{leftdef}, this implies that
for any left deformable functor $F: \C \to \D $, the functor $ F^I: \C^I \to \D^I $
has an absolute total left derived functor $ \L F^I: \Ho(\C^I) \to \Ho(\D^I) $ induced
by $ \bL F^I = \gamma_{\D^I} \, F^I \, Q^I $. Informally
speaking, the left derived functor of $ F^I $ is just the left derived functor of $F$ applied
objectwise.

Now, for a small category $ I $, let $\, \diag^{\C}_I: \C \to \C^I \,$ denote the
diagonal functor which assigns to an object $ A \in \Ob(\C) $ the constant
diagram $ \diag^{\C}_I(A):\, I \to \C $, $\,i \mapsto A \,$. Recall that
the colimit $ \colim^{\C}_I: \C^I \to \C $ is the left adjoint functor of $ \diag^{\C}_I $.
If $ \C $ is a homotopical category, we define the {\it homotopy colimit}
$\, \bL\colim^{\C}_I: \C^I \to \Ho(\C) $ to be the left derived functor of
$ \colim^{\C}_I $, and following our convention, write $  \L\colim^{\C}_I: \Ho(\C^I) \to \Ho(\C)$
for the corresponding total left derived functor. By Proposition~\ref{leftdef}, $ \bL\colim^{\C}_I $
exists if $ \colim^{\C}_I $ exists and is left deformable; in that case, since the diagonal
functor is homotopical, we have a deformable adjunction $\,\colim_I^{\C}:\, \C^I \rightleftarrows \C\,: \diag^{\C}_I $, and hence, by Corollary~\ref{cordef}, the derived adjunction
$$
\L\colim_I^{\C}:\, \Ho(\C^I) \rightleftarrows \Ho(\C): \, \overline{\diag}_{I}^{\, \C}
$$
After these preliminary remarks, we can state our second main theorem.
\begin{theorem}
\la{ThA3}
Let $ F:\,\C \rightleftarrows \D \,: G $ be a pair of adjoint functors
satisfying the conditions of Theorem~\ref{ThA2}. Assume, in addition, that for a small category $I$,
the functors $ \colim_I^{\C} $ and $ \colim_I^{\D} $ exist and are left deformable.
Then there is a natural isomorphism of functors
\begin{equation}
\la{hocomm}
\L F \circ \L \colim_I^{\C} \,\cong \, \L \colim_I^{\D} \circ \L(F^I)
\end{equation}
In other words, the  functor $\L F$ preserves homotopy colimits.
\end{theorem}
Theorem~\ref{ThA3} follows readily from Theorem~\ref{ThA2} and the main results of \cite{CS}
concerning homotopy colimits. For reader's convenience, we will summarize these results below,
before proving Theorem~\ref{ThA3}. We start with a simple lemma, which is probably well known
to experts, but since we could not find a reference, we provide a quick proof.
\begin{lemma}
\la{LeA4}
Let  $ \hat{F}:\,\M \rightleftarrows \D \,: \hat{G} $ be a deformable adjunction between homotopical
categories. Assume that, for a small category $I$, the functors $ \colim_I^{\M} $ and $ \colim_I^{\D} $ exist and are left deformable. Then there is a natural isomorphism
\begin{equation}
\la{hocomm1}
\L \hat{F} \circ \L \colim_I^{\M} \,\cong \, \L \colim_I^{\D} \circ \L(\hat{F}^I)
\end{equation}
\end{lemma}
\bproof
Since $ \hat{G} $ is right deformable, so is $ \hat{G}^I $, and there is a right
deformation functor on $\D$, say $\, R: \D \to \D \,$, such that
$\, \R \hat{G} = \gamma_{\M} \circ \hat{G} \circ R \,$ and
$\, \R(\hat{G}^I) = \gamma_{\M^{I}} \circ \hat{G}^I \circ R^I \,$ are
the right derived functors of $ \hat{G} $ and $ \hat{G}^I $, respectively.
Now, since $ \diag_I $ is homotopical,  we have obvious isomorphisms:
\begin{eqnarray*}
\overline{\diag}_{I}^{\, \M} \circ \R \hat{G} & = & \overline{\diag}_{I}^{\, \M} \circ  \gamma_{\M} \circ \hat{G} \circ R \\
& \cong & \gamma_{\M^{I}} \circ \diag_{I}^{\M} \circ \hat{G} \circ R \\
& \cong &
\gamma_{\M^{I}} \circ \hat{G}^I \circ \diag_{I}^{\D} \circ R \\
& \cong &
\gamma_{\M^{I}} \circ \hat{G}^I \circ R^I \circ \diag_{I}^{\D} \\
& \cong & \R(\hat{G}^I) \circ \diag_{I}^{\D} \ ,
\end{eqnarray*}
which induce an isomorphism of the total right derived functors
\begin{equation}
\la{obis}
\overline{\diag}_{I}^{\, \M} \circ \bR \hat{G} \,\cong \,
\bR(\hat{G}^I) \circ \overline{\diag}_{I}^{\,\D} \ .
\end{equation}
By lemma's assumptions, each functor in \eqref{obis} has a left adjoint, hence
\eqref{obis} implies \eqref{hocomm1}.
\eproof

Now, we briefly review the results of \cite{CS} needed for the proof of our Theorem~\ref{ThA3}.
We warn the reader that our notation differs from that of \cite{CS} but 
this should not cause confusion. For a small category $ I $, we denote by $ \Delta I $
the simplex category of $ I $, i.e. the category of simplices
$\, \Delta \downarrow N(I) \,$ of the nerve of $I$, and write $ \M^{\Delta I}_b $
for the full subcategory\footnote{In \cite{CS}, the category $ \M^{\Delta I}_b $ is denoted
$\,Fun^b({\bf N}(I), \M)  \,$.} of $ \M^{\Delta I} $ consisting of {\it bounded}
$\Delta I$-diagrams in a category $ \M $. Recall that a functor $ X: \Delta I \to \M $ is
{\it bounded} if it maps every degeneracy map $ s^i: s^i \sigma \to \sigma $ in $ \Delta I$
to an isomorphism in $ \M $; thus, modulo isomorphisms, a bounded functor is determined
by its values on nondegenerate simplices in $ \Delta I $.
The simplex category comes together with a forgetful functor
$ \tau: \Delta I \to I $ which takes an $n$-simplex $ \sigma $ in $ \Delta I $, i.e.
a chain $ \sigma = (i_0 \leftarrow i_1 \leftarrow \ldots \leftarrow i_n) $ of $n$
composable maps in $ I $, to its target $ \tau(\sigma) = i_0 $. This forgetful
functor yields the restriction functor $ \tau^*: \M^I \to \M^{\Delta I} $ whose image
is in $ \M^{\Delta I}_b $ (in fact, it is easy to check that $ \im(\tau^*) $ consists of
bounded functors $ X: \Delta I \to M $ which, in addition to inverting
all the degeneracy maps in $ \Delta I $, also invert all the boundary maps
$ d^i: d^i \sigma \to \sigma $ with $ i > 0 $). Now, if the category $ \M $ is closed under
colimits, the functor $\,\tau^*:  \M^I \to \M^{\Delta I}_b \,$ has a left adjoint
$\,\tau_*: \M^{\Delta I}_b \to \M^I \,$, which is given by restricting to
$ \M^{\Delta I}_b $ the left Kan extension  $\,{\rm Lan}_{\tau}: \M^{\Delta I} \to \M^I \,$ taken along $ \tau: \Delta I \to I \,$.  In this way, for any cocomplete category $ \M $, 
we get the adjunction
\begin{equation}
\la{BK1}
\tau_* :\, \M^{\Delta I}_b \, \rightleftarrows\,  \M^I \,: \, \tau^*
\end{equation}
In $ \M $ is a model category,  \eqref{BK1} is called the
{\it Bousfield-Kan approximation} of $ \M^I $. More generally,
in $\,l: \M \rightleftarrows \C: r\,$ is a left model approximation of
a homotopical category $ \C $, the composition of adjunctions
\begin{equation}
\la{BK2}
l^I \circ \tau_* :\, \M^{\Delta I}_b \, \rightleftarrows\,  \M^I  \, \rightleftarrows\, \C^I \,: \, \tau^* \circ r^I
\end{equation}
is called the {\it Bousfield-Kan approximation} of $ \C^I $. Now, the main results
of \cite{CS} can be encapsulated into the following theorem.
\bthm[\cite{CS}, Theorem~11.2 and Theorem 11.3]
\la{ThA4}
Let $I$ be a small category.

$(1)$ For any model category $ \M $, the category $ \M^{\Delta I}_b $ has a
model structure, where the weak equivalences (resp., fibrations) are the objectwise
weak equivalences (resp., fibrations) of bounded $ \Delta I$-diagrams in $\M$.

$(2)$ For any left model approximation $\,l: \M \rightleftarrows \C: r\,$, the Bousfield-Kan
approximation \eqref{BK2} is a left model approximation of $ \C^I $. In particular,
\eqref{BK1} is a left model approximation of $ \M^I $.

$(3)$ If $\,\C$ is closed under colimits and admits a
left model approximation $\,l: \M \rightleftarrows \C: r\,$, the corresponding
Bousfield-Kan approximation \eqref{BK2} is good for
$ \colim_I^{\C}: \C^I \to \C $. In particular, the functor $ \colim_I^{\C} $ is left deformable
and its left derived functor $($the homotopy colimit$)$
$ \bL \colim_I^{\C} $ exists.
\ethm

We now explain how to construct homotopy colimits using the Bousfield-Kan model
structure on $ \M^{\Delta I}_b $. To this end, we need another important observation
of \cite{CS} that, for any model category $ \M $, the functor
$ \colim^{\M}_{\Delta I}:  \M^{\Delta I}_b \to \M $, obtained by restricting the usual
colimit to bounded diagrams, is homotopical on cofibrant objects in $ \M^{\Delta I}_b $, and hence has a left derived functor (see \cite[Cor.~13.4 and Prop.~14.2]{CS}). Following \cite{CS},
we denote this  derived functor by\footnote{It is important to note that
the functor $ \ocolim^{\M}_{\Delta I} $ is {\it not}
equivalent, in general, to the usual homotopy colimit $ \bL \colim^{\M}_{\Delta I} $
restricted to $  \M_b^{\Delta I} $ (see \cite[Remark~14.3]{CS}).}
\begin{equation}
\la{ocolim}
\ocolim^{\M}_{\Delta I}:\, \M_b^{\Delta I} \to \Ho(\M)
\end{equation}
In terms of \eqref{ocolim},
the homotopy colimit functor on arbitrary $I$-diagrams
$ \bL\colim_{I}^{\M}: \M^I \to \Ho(\M) $ is given by
\begin{equation}
\la{ocolim2}
\bL\colim_{I}^{\M} \, \cong \,\ocolim^{\M}_{\Delta I} \circ \tau^*
\end{equation}
where $ \tau^*: \M^I \to \M_b^{\Delta I} $ is the restriction functor in
the Bousfield-Kan approximation \eqref{BK1}.
More generally, for a left model approximation $\,l: \M \rightleftarrows \C: r\,$,
the homotopy colimit $ \bL\colim_{I}^{\C}: \C^I \to \Ho(\C) $ is given by
the composition
$$
\C^I \xrightarrow{r^I} \M^I \xrightarrow{\tau^*} \M_b^{\Delta I}
\xrightarrow{\ocolim^{\M}_{\Delta I}} \Ho(\M) \xrightarrow{\L l} \Ho(\C)
$$
that is
\begin{equation}
\la{ocolim3}
\bL\colim_{I}^{\C} \, \cong \,\L l \circ \ocolim^{\M}_{\Delta I} \circ \tau^* \circ r^I
\end{equation}
Combining the isomorphisms \eqref{ocolim2} and \eqref{ocolim3} and passing to
total derived functors, we arrive at the following result
which we will use in the proof of Theorem~\ref{ThA3}.
\bcor
\la{lastcor}
Assume that a homotopical category $ \C $ admits a left model
approximation $ l: \M \rightleftarrows \C: r $ and is closed under colimits. Then, for any
small category $I$, $\, \L\colim_{I}^{\C} $ exists and
\begin{equation}
\la{ocolim4}
\L\colim_{I}^{\C} \, \cong \,\L l \circ \L\colim_{I}^{\M} \circ \bar{r}^I
\end{equation}
\ecor
Finally, we turn to

\bproof[Proof of Theorem~\ref{ThA3}]
By Lemma~\ref{LeA4}, we have a natural isomorphism \eqref{hocomm1} which yields
by restriction:
\begin{equation}
\la{hocomm3}
\L \hat{F} \circ \L \colim_I^{\M} \circ \bar{r}^I \,\cong \,
\L \colim_I^{\D} \circ \L(\hat{F}^I) \circ \bar{r}^I
\end{equation}
Now, by Theorem~\ref{ThA2} (see  \eqref{Lleft}) and
Corollary~\ref{lastcor}, the composition of functors in the left-hand side of \eqref{hocomm3} is isomorphic to
\begin{equation*}
\L F \circ \L l \circ \L \colim_I^{\M} \circ \bar{r}^I \cong \L F \circ \L\colim_{I}^{\C}
\end{equation*}
On the other hand, by condition $(ii)$ of Theorem~\ref{ThA2}, the pair of functors
$(\hat{F}^I, r^I) $ is left deformable, and there is a natural weak equivalence
$ \hat{F}^I \circ r^I = (\hat{F} \circ r)^I \xrightarrow{\sim} F^I $, inducing an
isomorphism
$$
\L(\hat{F}^I) \circ \bar{r}^I \cong \L(\hat{F}^I \circ r^I) \cong \L(F^I)\ .
$$
Hence, the right-hand side of \eqref{hocomm3} is isomorphic to $ \L \colim_I^{\D} \circ \L(F^I) $.
Combining \eqref{hocomm3} with these two isomorphisms gives \eqref{hocomm}.
\eproof
%


\begin{remark}
The assumption of Theorem~\ref{ThA3} that $ \colim_I^{\C} $ is a left deformable functor
is superfluous. Indeed, thanks to Theorem~\ref{ThA4}$(3)$, it suffices only to assume the existence of
$ \colim_I^{\C} $.
\end{remark}

\end{document}